\magnification=1000
\hsize=11.7cm
\vsize=18.9cm
\lineskip2pt \lineskiplimit2pt
\nopagenumbers

\hoffset=-1truein
\voffset=-1truein

\advance\voffset by 4truecm
\advance\hoffset by 4.5truecm

\newif\ifentete

\headline{\ifentete\ifodd	\count0 
      \rlap{\head}\hfill\tenrm\llap{\the\count0}\relax
    \else
        \tenrm\rlap{\the\count0}\hfill\llap{\head} \relax
    \fi\else
\global\entetetrue\fi}

\def\entete#1{\entetefalse\gdef\head{#1}} 
\entete{}

\input amssym.def
\input amssym.tex

\def\-{\hbox{-}}
\def\.{{\cdot}}
\def\O{{\cal O}}

\def\F{{\cal F}}

\def\L{{\cal L}}
\def\M{{\cal M}}
\def\N{{\cal N}}
\def\P{{\cal P}}
\def\Q{{\cal Q}}
\def\G{{\cal G}}
\def\T{{\cal T}}

\def\B{{\cal B}}

\def\X{{\cal X}}

\def\Z{{\cal Z}}

\def\C{{\cal C}}

\def\ab{\frak a\frak b}
\def\Ab{\frak A\frak b}

\def\Gr{\frak G\frak r}

\def\Fct{\frak F\frak c\frak t}

\def\Ker{\frak K\frak e\frak r}

\def\id{\frak i\frak d}
\def\int{\frak i\frak n\frak t}

\def\qq{\quad{\rm and}\quad}

\def\tr{\frak t\frak r}
\def\mod{\frak m\frak o\frak d}

\def\too{\longrightarrow}

 3
 2
\font\large=cmr10  scaled \magstep 2
 2
\font\larti=cmti10  scaled \magstep 2
 1
\font\larsy=cmsy10  scaled \magstep 2

\font\cds=cmr7
\font\cdt=cmti7
\font\cdy=cmsy7
\font\cdi=cmmi7

\count0=1

\centerline{\large The perfect {\larsy F}-locality from the basic {\larsy F}-locality}
\medskip
\centerline{\large over a Frobenius {\larti P}-category {\larsy F}}
\medskip
\centerline{\bf Lluis Puig }
\medskip
\noindent 
\centerline{\cds 6 Av Bizet, 94340 Joinville-le-Pont, France}

\medskip
\noindent
{\bf Abstract:} {\cds   Let {\cdt p} be a prime, {\cdt P} a finite {\cdt p\-}group, {\cdy F} a Frobenius 
{\cdt P\-}category and  $\scriptstyle \F^{^{\rm sc}}$ the full subcategory of {\cdy F} over  the set of   {\cdt {\cdy F}-selfcentralizing\/} subgroups of~{\cdt P}. Recently, we have understood an easy way to obtain the {\cdt perfect $\scriptstyle \F^{^{\rm sc}}\-$locality} $\scriptstyle \P^{^{\rm sc}}$  from the {\cdt basic $\scriptstyle \F^{^{\rm sc}}\-$locality}~$\scriptstyle \L^{^{{\rm b},sc}}\!$: it depends on a suitable filtration of  the {\cdt basic $\scriptstyle \F\-$locality} $\scriptstyle \L^{^{\rm b}}\!$ and on a vanishing cohomology result, given with more generality in [11, Appendix].}

\bigskip
\noindent
{\bf £1. Introduction }
\medskip
£1.1. Let $p$ be a prime and $P$ a finite $p\-$group. After our introduction of the Frobenius $P\-$categories $\F$ [7] and the question of Dave Benson [1] on the existence  of  a suitable category  $\P^{^{\rm sc}}\!$ --- called {\it linking system\/} in [2] and {\it perfect $\F^{^{\rm sc}}\-$locality\/} in  [8,~Chap.~17] --- extending the {\it full\/} subcategory  $\F^{^{\rm sc}}$ of~$\F$ over the set of   {\it $\F\-$selfcentralizing\/} subgroups of~$P$~[8,~Chap.~3], the existence and the uniqueness of 
$\P^{^{\rm sc}}\!$ has concentrate some effort.

\medskip
£1.2.  In [2]  Carles Broto, Ran Levi and Bob Oliver formulate the existence and the uniqueness of the category $\P^{^{\rm sc}}$ in terms of the annulation of an {\it obstruction $3\-$cohomology element\/} and of the vanishing of a  {\it $2\-$cohomology group\/}, respectively. They actually state a sufficient condition for the va-nishing of the corresponding {\it $n\-$cohomology groups\/}.

\medskip
£1.3.  In [3] Andrew Chermak has proved the existence and the uniqueness of $\P^{^{\rm sc}}\!$ {\it via\/} his {\it objective partial groups\/}, but his proof depends on the so-called  {\it Classification of the finite simple groups\/} and on some  results by U. Meierfrankenfeld and B. Stellmacher. In [6] Bob Oliver, following some of Chermak's methods, has also proved for $n\ge 2$ the vanishing of the {\it $n\-$cohomology groups\/} mentioned above. In [5] George Glauberman and Justin Lynd remove the use of the {\it Classification of the finite simple groups\/} in [6]{\footnote{\dag}{\cds 
Although they need a {\cdi partial} {\cdi classification} for ${\scriptstyle p = 2}$.}} .

\medskip
£1.4. Independently, with direct methods which already employ the {\it basic $\F\-$locality\/} $\L^{^{\rm b}}$ [8,~Chap.~22],  in [9] and [10]{\footnote{\dag\dag}{\cds In [10] we give a full correction of the uniqueness argument for {\cdy P ${\scriptstyle ^{\rm sc}\!}$} in [9].}} we prove the existence and the uniqueness of an extension 
$\P$ of $\F$ --- called {\it perfect $\F\-$locality\/} in [8,~Chap.~17] --- which contains $\P^{^{\rm sc}}\!$ as the {\it full\/} subcategory over the set of   {\it $\F\-$selfcentralizing\/} subgroups of~$P$~[8,~Chap.~3].

\medskip
£1.5. But recently, we have understood an easier way to obtain $\P^{^{\rm sc}}\!$ from the {\it full\/} subcategory~$\L^{^{{\rm b},sc}}$ of $\L^{^{\rm b}}$  over the set of   {\it $\F\-$selfcentralizing\/} subgroups of~$P$~[8,~Chap.~3]. Denoting by $\Z^{^{\rm sc}}\,\colon \L^{^{{\rm b}, sc}}\to \Ab$ the obvious {\it contravariant\/} functor from $\L^{^{{\rm b},sc}}$ to the category $\Ab$ of finite Abelian groups, mapping any {\it $\F\-$selfcentralizing\/} subgroup $Q$ of $P$ on its center $Z(Q)\,,$ it is easy to see that we have a quotient category $\widetilde{\L^{^{{\rm b}, sc}}} =\L^{^{{\rm b}, sc}}/\Z^{^{\rm sc}}$ and that the structural functor $\pi^{^{\rm sc}}\,\colon \L^{^{{\rm b}, sc}}\to \F^{^{\rm sc}}$ factorizes  through a functor $\widetilde{\pi^{^{\rm sc}}}\,\colon \widetilde{\L^{^{{\rm b}, sc}}}\to \F^{^{\rm sc}}\,.$

\medskip
£1.6. The point is that $\widetilde{\pi^{^{\rm sc}}}$ admits an {\it essentially unique section functor\/}
$\widetilde{\sigma^{^{\rm sc}}}\,\colon \F^{^{\rm sc}}\to \widetilde{\L^{^{{\rm b}, sc}}}\,,$ and then 
$\P^{^{\rm sc}}\!$ is just the {\it converse image\/} in $\L^{^{{\rm b}, sc}}$ of the {\it image\/} 
$\widetilde{\sigma^{^{\rm sc}}} (\F^{^{\rm sc}})$ of $\F^{^{\rm sc}}$ in 
$ \widetilde{\L^{^{{\rm b}, sc}}}\,;$ since in [9,~Theorem~7.2] we prove that any perfect 
$\F^{^{\rm sc}}\-$locality  $\P^{^{\rm sc}}$ can be extended to a unique  perfect 
$\F\-$locality~$\P\,,$ this proves the existence of $\P\,.$ Moreover,
in [9,~8.5 and~Theorem~8.10] we prove that there is an {\it $\F\-$locality functor\/} $\sigma$
from any perfect $\F\-$locality~$\P$ to $\L^{^{\rm b}}\,;$ then, it is easy to check that $\sigma$
induces a  functor 
$\widetilde{\sigma^{^{\rm sc}}}\,\colon \F^{^{\rm sc}}\to \widetilde{\L^{^{{\rm b}, sc}}}$ which
is a section of $\widetilde{\pi^{^{\rm sc}}}\,,$ proving the uniqueness of $\P^{^{\rm sc}}\!$ and therefore, by [9,~Theorem~7.2], the uniqueness of $\P\,.$

\medskip
£1.7. In Section £2 we recall all the definitions and state properly all the quoted results. The existence and the essential uniqueness of the {\it section funtor\/} $\widetilde{\sigma^{^{\rm sc}}}$ mentioned above depend on a suitable {\it filtration\/} of $\L^{^{\rm b}}$ and on a 
{\it vanishing cohomology result\/}; this {\it filtration\/} has been already introduced in [9,~8.3 and~Corollary~8.4], but it seems necessary to give here a complete account  in Section~£3. The {\it vanishing cohomology result\/} we need here is given
in [11, Appendix] in a more general framework. In Section~£4 we give explicit proofs of all the results announced in~£1.6 above and, in particular, an independent proof of the existence of the functor $\sigma\,\colon \P\to \L^{^{\rm b}}$ mentioned above.

\bigskip
\noindent
{\bf £2.  Definitions and quoted results\/}

\medskip
£2.1. Denote by $\Ab$ the category of  Abelian groups and  by  $\frak i\Gr$ the category formed by the finite groups and by the injective group  homomorphisms. Recall that, for any category $\frak C\,,$ we denote by $\frak C^\circ$  the {\it opposite\/} category and, for any $\frak C\-$object $A\,,$ by $\frak C_A$  (or $(\frak C)_A$ to avoid confusion)  the category of ``$\frak C\-$morphisms to $A$'' [8, 1.7]; moreover, for any pair of objects $A$ and $B\,,$
$\frak C (B,A)$ denote the set of $\frak C\-$morphisms from $A$ to~$B$ and we set 
$\frak C (A) = \frak C (A,A)$ for short.

\medskip
£2.2.  For any finite subgroup $G$ and any $p\-$subgroup $P$ of $G\,,$  denote by $\F_{G,P}$ and  $\T_{G,P}$ the respective categories where the objects are all the  subgroups of $P$ and, for
two of them $Q$ and $R\,,$ the respective sets of morphisms $\F_{G,P} (Q,R)$
and $\T_{G,P} (Q,R)$  are formed by the group homomorphisms from $R$ to $Q$ respectively induced by the conjugation by elements of $G\,,$ and by the set $T_G (R,Q)$ of such elements,
the {\it compositions\/} being the obvious ones.

\medskip
£2.3. For a finite $p\-$group $P\,,$ a {\it Frobenius  $P\-$category\/} $\F$ is a subcategory of $\frak i\Gr$ containing $\F_{\!P} = \F_{P,P}$ where the objects are all the  subgroups of $P$
and the morphisms fulfill the following three conditions [8, 2.8 and Proposition~2.11]
\smallskip
\noindent
£2.3.1\quad {\it For any subgroup $Q$ of $P\,,$ the inclusion functor $(\F)_Q\to 
(\frak i\Gr)_Q$ is full.\/}
\smallskip
\noindent
£2.3.2\quad {\it $\F_P (P)$ is a Sylow $p\-$subgroup of $\F (P)\,.$\/}
\smallskip
\noindent
We say that a subgroup $Q$ of $P$ is {\it fully centralized\/} in $\F$ if for any $\F\-$morphism
$\xi\,\colon Q\.C_P(Q)\to P$ we have $\xi \big(C_P (Q)\big) = C_P\big(\xi (Q)\big)\,;$
similarly, replacing  in this condition the centralizer by the normalizer, we say that $Q$ is {\it fully normalized\/}.
\smallskip
\noindent
£2.3.3\quad {\it For any  subgroup $Q$ of $P$ fully centralized in $\F\,,$
any $\F\-$morphism $\varphi\,\colon Q\to P$ and any subgroup $R$ of $N_P\big(\varphi(Q)\big)$ such that  $\varphi (Q)\i
R$ and that $\F_P(Q)$ contains the action of $\F_R \big(\varphi(Q)\big)$ over $Q$ via $\varphi\,,$ there exists an $\F\-$morphism 
$\zeta\,\colon R\to P$ fulfilling $\zeta\big(\varphi (u)\big) = u$ for any $u\in Q\,.$\/}
\smallskip
\noindent
 We denote by $\tilde\F$ -- called the {\it exterior quotient\/} of~$\F$ --- the quotient of $\F$ by the {\it inner\/} automorphisms of the $\F\-$objects  [8,~1.3] and by $\epsilon_\F : \F\to \tilde\F$
 the {\it canonical functor\/}. Note that, with the notation above, if $P$ is Sylow $p\-$subgroup of $G$ then $\F_{G,P}$ is a {\it Frobenius  $P\-$category\/}; often, we write $\F_G$ instead of $\F_{G,P}\,.$

\medskip 
£2.4. Then, a ({\it divisible\/}) {\it $\F\-$locality\/}{\footnote{\dag}{\cds Here we only consider {\cdt divisible
 {\cdy F}-localities} in the sense of [8,~Chap.~17].}} is a triple $(\tau,\L,\pi)$ formed by a {\it finite\/} category $\L\,,$ a {\it surjective\/} functor $\pi\,\colon \L\to \F$ and a functor 
$\tau\,\colon \T_P\to \L$ from the {\it transporter\/} category $\T_P= \T_{P,P}$ of $P\,,$ fulfilling the following two conditions [8,~17.3 and 17.8]
\smallskip
\noindent
£2.4.1\quad {\it The composition $\pi\circ\tau$ coincides with the composition of the canonical functor
defined by the conjugation $\kappa_P\,\colon \T_P\to \F_P$ with the inclusion $\F_P\i \F\,.$\/}
\smallskip
\noindent
We denote by $\tilde\kappa_P\,\colon \T_P\to \tilde\F_{\!P}$ the composition of $\kappa_P$ with $\epsilon_{\F_{\!P}}$
above.
\smallskip
\noindent
£2.4.2\quad {\it  For any pair of subgroups $Q$ and $R$ of $P\,,$ ${\rm Ker}(\pi_{_R})$ acts regularly on the fibers of the following maps  determined by $\pi$}
$$\pi_{_{Q,R}} : \L (Q,R)\too \F(Q,R)$$
\smallskip
\noindent
 Analogously, for any pair of subgroups $Q$ and $R$ of $P\,,$ we denote by 
$$\tau_{_{Q,R}} : \T_P (Q,R)\too \L (Q,R)
\eqno £2.4.3\phantom{.}$$
the map determined by $\tau\,,$ and whenever  $R\i Q$ we set 
$i_{_R}^{^Q}  = \tau_{_{Q,R}}(1)\,;$ if $R = Q$ then we write $Q$ once  for short.

\medskip
£2.5. We say that an {\it $\F\-$locality\/} $(\tau,\L,\pi)\,,$ or $\L$ for short, is {\it coherent\/}
if it fulfills the following  condition [8,~17.9]
\smallskip
\noindent
 (\Q)\quad {\it  For any pair of subgroups $Q$ and $R$ of $P\,,$ any $x\in \L (Q,R)$ and any $v\in R\,,$
we have $x\.\tau_{_R} (v) = \tau_{_Q} \big(\pi_{_{Q,R}}(x) (v)\big)\.x\,.$ }
\smallskip
\noindent
In this case,, if $Q'$ and $R'$ are subgroups of $P\,,$ and we have the inclusions $R\i Q$ and $R'\i Q'\,,$ denoting by $\L (Q',Q)_{R',R}$ the set of $y\in \L (Q',Q)$ such that
$\big(\pi_{_{Q',Q}} (y)\big)(R) \i R'\,,$ we get a {\it restriction\/} map (possibly empty!)
$$r_{_{R',R}}^{_{Q',Q}} : \L (Q',Q)_{R',R}\too \L (R',R)
\eqno £2.5.1\phantom{.}$$
fulfilling $y\. i_{_R}^{^Q} = i_{_{R'}}^{^{Q'}}\.r_{_{R',R}}^{_{Q',Q}} (y)$ for any 
$y\in \L (Q',Q)_{R',R}\,.$ Note that, with the notation above, if $P$ is Sylow $p\-$subgroup of $G$ then $\T_G = \T_{G,P}$ endowed with the obvious functors
$$\tau^{_G} : \T_P\too \T_G\qq \pi^{_G} :\T_G\too \F_G = \F_{G,P}
\eqno £2.5.2\phantom{.}$$
becomes a {\it coherent $\F_G\-$locality\/}.
Moreover, we say that a {\it coherent $\F\-$locality\/} $(\tau,\L,\pi)$ is {\it $p\-$coherent\/} 
(resp. {\it ab-coherent\/}) when 
${\rm Ker}(\pi_{_Q})$ is a  finite  $p\-$group (resp. a finite abelian group) for any subgroup $Q$ of $P\,.$

\medskip
£2.6. Recall that  the $\F\-${\it hyperfocal subgroup\/} is the subgroup~$H_{\!\F}$
of~$P$ generated by the union of the sets $\big\{u^{-1}\sigma (u)\big\}_{u\in Q}$
where $Q$ runs over the set of subgroups of~$P$ and $\sigma$ over the set of
{\it $p'\-$elements\/} of~$\F (Q)\,.$ We say that an {\it $\F\-$locality\/} $(\hat\tau,\P,\hat\pi)$
is {\it perfect\/} if $\P$ is coherent and, for any subgroup $Q$ of $P$ {\it fully centralized\/} 
in~$\F\,,$ the {\it $C_\F (Q)\-$hyperfocal subgroup\/} $H_{C_\F (Q)}$ coincides with 
${\rm Ker} (\hat\pi_Q)$ [8, 17.13]; actually, this is equivalent to say that $\P (Q)\,,$ endowed with
$$\hat\tau_Q : \T_{N_P (Q)} \too \P (Q) \qq \hat\pi_Q : \P (Q)\too \F (Q)
\eqno £2.6.1,$$
is an {\it $\F\-$localizer\/} of $Q$ [8, 18.5 and Theorem 18.6], for  any  subgroup $Q$ of $P$ 
{\it fully centralized\/} in~$\F$.

\medskip
£2.7. Further, for any {\it $\F\-$locality\/} $(\tau,\L,\pi)$ we get a {\it contravariant\/} functor from $\L$  to the category $\Gr$ of finite groups [8,~17.8.2]
$$\Ker (\pi) : \L\too \Gr
\eqno £2.7.1\phantom{.}$$
sending any subgroup $Q$ of $P$ to ${\rm Ker}(\pi_{_Q})$ and any $\L\-$morphism 
$x\,\colon R\to Q$ to the group homomorphism 
$$\Ker (\pi)_x : {\rm Ker}(\pi_{_Q})\too {\rm Ker}(\pi_{_R})
\eqno £2.7.2\phantom{.}$$ 
fulfilling $u\. x = x\. \big(\Ker (\pi)_x (u)\big)$ for any $u\in {\rm Ker}(\pi_{_Q})\,.$ 
If $\L$ is {\it ab-coherent\/} then the functor $\Ker (\pi)$ factorizes through
the {\it exterior quotient\/} $\tilde\F$, inducing a functor
$$\widetilde\Ker (\pi) = \tilde{\frak k}_\L : \tilde\F\too \Gr
\eqno £2.7.3;$$
 indeed, in this case it follows from [8, Proposition 17.10]
that, for any subgrpup $Q$ of~$P\,,$ $\tau_{_Q} (Q)$ centralizes ${\rm Ker}(\pi_{_Q})$
and therefore, for any $u\in {\rm Ker}(\pi_{_Q})$ and any $v\in Q\,,$ we have
$$\tau_{_Q} (v)\. u = u.\tau_{_Q} (v) = \tau_{_Q} (v)\.\big(\Ker (\pi)_{\tau_{_Q} (v)} (u)\big)
\eqno £2.7.4,$$
so that $\Ker (\pi)_{\tau_{_Q} (v)}= {\rm id}_{{\rm Ker}(\pi_{_Q})}\,;$ the same argument holds
for $w\in{ \rm Ker}(\pi_{_Q})\,.$

\medskip
£2.8.  If $\L'$ is a second  $\F\-$locality with {\it structural functors\/} $\tau'$ 
and~$\pi'\,,$  we call  {\it $\F\-$locality functor\/} from $\L$ to $\L'$ any functor 
$\frak l\colon \L\to \L'$ fulfilling 
$$\tau' = \frak l\circ\tau\qq \pi'\circ\frak l = \pi
\eqno £2.8.1;$$
the composition of two $\F\-$locality functors is obviously an $\F\-$locality functor.
It is easily checked that
any $\F\-$locality functor  $\frak l\colon \L\to \L'$  determines a {\it natural map\/}
$$\chi_{\frak l} : \Ker (\pi)\too \Ker (\pi')
\eqno £2.8.2;$$
conversely, it is quite clear that any subfunctor $\frak h$ of $\Ker (\pi)$ determines a 
{\it quotient  $\F\-$locality\/} $\L/\frak h$ defined by the quotient sets
$$(\L/\frak h)(Q,R) = \L (Q,R)/\frak h (R)
\eqno £2.8.3,$$
 for any pair of subgroups $Q$ and $R$ of $P\,,$ and by the corresponding induced composition; moreover, $\L/\frak h$ is {\it coherent\/} whenever $\L$ is it.

 \medskip
 £2.9.  We say that two  $\F\-$locality functors 
 $\frak l$ and  $\bar\frak l$ from 
$\L$ to $\L'$ are {\it naturally $\F\-$isomorphic\/} if we have a {\it natural isomorphism\/} 
$\lambda\colon \frak l \cong \bar\frak l$ fulfilling $\pi'* \lambda = \id_\pi$ and 
$\lambda * \tau = \id_{\tau'}\,;$
 in this case,  
$\lambda_Q$ belongs to ${\rm Ker} (\pi'_{_Q})$ for any  subgroup $Q$  of $P$ and, since 
$\frak l (i_{_Q}^{^P}) = i'^{^P}_{_Q}\! = \bar\frak l  (i_{_Q}^{^P}) \,,$  $\lambda$ is uniquely determined by~$\lambda_P\,;$ indeed, we have
$$\lambda_P\.i'^{^P}_{_Q}\!  = i'^{^P}_{_Q}\.\lambda_Q
\eqno £2.9.1.$$
Once again, the composition of a natural $\F\-$isomorphism with an 
$\F\-$locality functor or with another such a natural $\F\-$isomorphism
is  a natural $\F\-$iso-morphism. 

\medskip
£2.10. Moreover, from two $\F\-$localities $(\tau,\L,\pi)$ and $(\tau',\L',\pi)\,,$ we can construct a third  $\F\-$locality $\L'' = \L \times_{\F} \L'$ from the corresponding  category defined by the {\it pull-back\/} of sets
$$\L''(Q,R) = \L (Q,R)\times_{\F (Q,R)}\L' (Q,R)
\eqno £2.10.1\phantom{.}$$
with the obvious composition and with the structural maps
$$\T_P (Q,R)\buildrel \tau''_{_{Q,R}}\over{\hbox to 30pt{\rightarrowfill}} \L''(Q,R)\buildrel \pi''_{_{Q,R}}\over{\hbox to 30pt{\rightarrowfill}}\F(Q,R)
\eqno £2.10.2\phantom{.}$$
respectively induced by $\tau$ and $\tau'\,,$ and by $\pi$ and $\pi'\,.$ Note that we have obvious 
$\F\-$locality functors 
$$\L\longleftarrow \L \times_{\F} \L'\too \L'
\eqno £2.10.3\phantom{.}$$
and that $ \L \times_{\F} \L'$ is {\it coherent\/} if $\L$ and $\L'$ are so.

\medskip
£2.11. In order to define the {\it basic $\F\-$locality\/}, we have to consider the  {\it $\F\-$basic $P\times P\-$sets\/}; recall that an {\it $\F\-$basic $P\times P\-$set\/} $\Omega$  is a finite nonempty $P\times P\-$set $\Omega$ fulfilling the following three conditions 
[8,~21.4 and~21.5], where  $\Omega^\circ$ denotes the $P\times P\-$set obtained from 
$\Omega$ by exchanging both factors, for any subgroup $Q$ of $P$ we denote by $\iota_Q^P$ the corresponding inclusion map, and for any $\varphi\in \F (P,Q)$ we set$$\Delta_{\varphi} (Q)  =\{(\varphi (u),u)\}_{u\in Q}
\eqno £2.11.1.$$
\smallskip
\noindent
£2.11.2 \quad {\it  We have $\Omega^\circ \cong \Omega\,,$ $\{1\}\times P$ acts  freely on 
$\Omega$ and $\vert\Omega\vert/\vert P\vert \not\equiv 0 \bmod{p}\,.$\/}
\smallskip
\noindent
£2.11.3 \quad {\it  For any subgroup $Q$ of $P$ and any $\varphi\in \F (P,Q)$  we have a 
$Q\times P\-$set isomorphism
$${\rm Res}_{\varphi\times {\rm id}_P} (\Omega)\cong 
{\rm Res}_{\iota_Q^P\times {\rm id}_P} (\Omega)$$\/}
£2.11.4 \quad {\it Any $P\times P\-$orbit in $\Omega$ is isomorphic to $(P\times P)/\Delta_\varphi (Q)$
for a suitable subgroup $Q$ of $P$ and some $\varphi\in \F (P,Q)\,.$\/}
\smallskip
\noindent
Moreover,  we say that an $\F\-$basic $P\times P\-$set $\Omega$ is {\it thick\/} if the multiplicity of the indecomposable $P\times P\-$set  $(P\times P)\big/\Delta_\varphi (Q)$  is at least two for any subgroup $Q$ of $P$ and any $\varphi\in \F (P,Q)$ [8,~21.4].
 
 \medskip 
 £2.12. The existence of a {\it thick $\F\-$basic $P\times P\-$set\/}  is guaranteed by
 [8,~Proposition~21.12]; we choose one of them $\Omega$  and denote by $G$ the group of 
$\{1\}\times P\-$set automorphisms  of~${\rm Res}_{\{1\}\times P}(\Omega)\,;$ it is clear that we have an injective~map from $P\times \{1\}$ into $G$ and we identify this image with the 
$p\-$group $P$ itself so that, from now on, $P$ is contained in $G$ and acts freely on $\Omega\,.$ Then, it follows from the conditions above
that we have
$$\F_{G,P} = \F
\eqno £2.12.1\phantom{.}$$
and it is quite clear that, as in~£2.5.2, $\T_G = \T_{G,P}$ becomes a {\it coherent $\F\-$locality.\/}

\medskip
£2.13. For any subgroup $Q$ of $P\,,$ it is clear that  the centralizer $C_G (Q)$ coincides with the group of {\it $Q\times P\-$set automorphisms\/} of ${\rm Res}_{Q\times P}(\Omega)\,;$ moreover, since any $Q\times P\-$orbit in $\Omega$ is isomorphic to the $Q\times P\-$set 
$(Q\times P)\big/\Delta_\eta (T)\,,$ for a suitable subgroup $T$ of $P$ such that 
$\F (Q,T)\not= \emptyset$ and some $\eta\in \F (Q,T)$ (cf.~condition~£2.11.3), and since we have  [8,~22.3]
$${\rm Aut}_{Q\times P}\big((Q\times P)\big/\Delta_\eta (T)\big)\cong \bar N_{Q\times P}\big(\Delta_\eta (T)\big)
\eqno £2.13.1,$$ 
denoting by $k_\eta$ the {\it multiplicity\/} of $(Q\times P)\big/\Delta_\eta (T)$ in $\Omega$ and 
by $\frak S_{k_\eta}$ the cor-responding {\it $k_{\psi}\-$symmetric\/} group, we actually get  obvious group isomorphisms
$$ C_G (Q)\cong \prod_{T\in \C_P}\,\prod_{\eta\in \frak O_Q^T} \bar N_{Q\times P}\big(\Delta_\eta (T)\big)\wr  \frak S_{k_{\eta}}
\eqno £2.13.2\phantom{.}$$
where $\wr$ denotes the {\it wreath product\/}, $\C_P$ is a set of representatives  for the set of $P\-$conjugacy classes of subgroups $T$ of $P$ and, for any $T\in \C_P\,,$\break
\eject
\noindent
$\frak O_Q^T\i \F (Q,T)$ is a (possibly empty) set of representatives for the quotient set $Q\backslash\F (Q,T)/N_P (T)\,.$  For short, let us set 
$$\frak O_Q = \bigsqcup_{T\in \C_P}\frak O_Q^T
\eqno £2.13.3;$$
this set actually indexes the {\it  set of isomorphic classes of transitive $Q\times P\-$sets\/};
to avoid confusion, we note by $(T,\eta)$ the element $\eta$ in~$\frak O_Q^T\,.$

\medskip
£2.14. Then, it follows from [8,~Proposition~22.11] that {\it  the correspondence sending $Q$ to the minimal normal subgroup $\frak S_G (Q)$ of $C_G (Q)$ containing the image of 
${\displaystyle \prod_{(T,\eta)\in \frak O_Q}} \frak S_{k_{\eta}}$ for any  isomorphism 
{\rm £2.13.2} induces a functor \/}
$$\frak S_G  : \T_G\too \frak i\Gr
\eqno £2.14.1;$$ 
it is actually a subfunctor of $\Ker (\pi^{_G})$ (cf.~£2.5.2) and therefore determines a 
{\it coherent $\F\-$locality\/} $\L^{^{\rm b}} = \T_G/\frak S_G$ (cf.~£2.8) --- called the {\it basic $\F\-$locality\/}
[8,~Chap.~22] --- which, according to [9,~Corollary~4.11], does not depend on the choice of the 
{\it thick $\F\-$basic $P\times P\-$set\/} $\Omega\,.$ Moreover, denoting by 
$$\tau^{_{\rm b}} : \T_P\too \L^{^{\rm b}}\qq \pi^{_{\rm b}} : \L^{^{\rm b}}\too \F
\eqno £2.14.2\phantom{.}$$
the structural functors,  it follows from [8,~Proposition~22.7] that, for any subgroup $Q$ of $P\,,$ isomorphisms in £2.13.2 induce a 
{\it canonical isomorphism\/}
$$\big(\Ker (\pi^{_{\rm b}})\big)(Q) \cong \prod_{(T,\eta)\in \frak O_Q}
\ab\Big(\bar N_{Q\times P}\big(\Delta_\eta (T)\big)\Big) 
\eqno £2.14.3\phantom{.}$$
where $\ab\,\colon \Gr\to \Ab$ denotes the obvious functor mapping any finite group $H$~on its maximal Abelian quotient~$H/[H,H]\,;$ in particular, note that $\L^{^{\rm b}}$ is 
{\it $p\-$coherent\/} (cf.~£2.5).

\medskip
£2.15. Moreover, any $\L^{^{\rm b}}\-$morphism $x\,\colon R\to Q$ can be lifted to an element
$\hat x\in G$ fullfilling
$\hat x\circ R\circ \hat x^{-1}\i Q$ in the group of bijections of $\Omega\,;$ in particular, we also have
$$\hat x^{-1}\circ C_G (Q)\circ \hat x\i C_G (R)
\eqno £2.15.1\phantom{.}$$
and, considering isomorphisms £2.13.2 for both $C_G(Q)$ and $C_G (R)$, it is clear that
the conjugation by $\hat x^{-1}$ sends the factor determined by $T\in \C_P$ and by 
$\eta\in \frak O_Q^T$ in  some factors determined by $U\in \C_P$ and by 
$\theta\in \frak O_R^U$ in such a way that, setting $\varphi = \pi^{_{\rm b} }(x)$, there exists an injective  $R\times P\-$set homomorphism
$$f :(R\times P)\big/\Delta_\theta (U) \too 
{\rm Res}_{\varphi\times {\rm id}_P}\big((Q\times P)\big/\Delta_\eta (T)\big)
\eqno £2.15.2\phantom{.}$$
or, equivalently, we have 
$$\Delta_{\varphi\circ\theta} (U) = \big(\varphi (R)\times P\big)\cap {}^{(u,n)}\Delta_\eta (T)
\eqno £2.15.3$$
for suitable $u\in Q$ and $n\in P\,.$
\eject

\medskip
£2.16. More precisely, the $\L^{^{\rm b}}\-$morphism $x\,\colon R\to Q$ determines the group homomorphism
$$\big(\Ker (\pi^{_{\rm b}})\big)(x) : \big(\Ker (\pi^{_{\rm b}})\big)(Q)\too 
\big(\Ker (\pi^{_{\rm b}})\big)(R)£2.17
\eqno £2.16.1;$$
considering isomorphisms £2.14.3 for both  $ \big(\Ker (\pi^{_{\rm b}})\big)(Q)$ and 
$\big(\Ker (\pi^{_{\rm b}})\big)(R)$, it makes sense to introduce the projection  in 
$\ab\Big(\bar N_{R\times P}\big(\Delta_\theta (U)\big)\Big) $  of the restriction 
of~$\big(\Ker (\pi^{_{\rm b}})\big)(x)$ to £2.17
$\ab\Big(\bar N_{Q\times P}\big(\Delta_\eta (T)\big)\Big)$ --- noted 
$\big(\Ker (\pi^{_{\rm b}})\big)(x)_{_{(T,\eta)}}^{^{(U,\theta)}}$  --- for any $(T,\eta)\in \frak O_Q$ and any $(U,\theta)\in \frak O_R\,;$ according to £2.15 above, $\big(\Ker (\pi^{_{\rm b}})\big)(x)_{_{(T,\eta)}}^{^{(U,\theta)}}\not= 0$ forces
$$\Delta_{\varphi\circ\theta} (U) = \big(\varphi (R)\times P\big)\cap {}^{(u,n)}\Delta_\eta (T)
\eqno £2.16.2$$
for suitable $u\in Q$ and $n\in P\,.$

 \medskip
£2.17. In this case, in [8, Proposition 22.17] we describe 
$\big(\Ker (\pi^{_{\rm b}})\big)(x)_{_{(T,\eta)}}^{^{(U,\theta)}}$ as follows. Consider the  set of  injective $R\times P\-$set homomorphisms as in~£2.15.2 above; it is clear that 
$\bar N_{R\times P}\big(\Delta_\theta (U)\big) \times \bar N_{Q\times P}\big(\Delta_\eta (T)\big)$ acts on this set by left- and right-hand composition and, denoting by 
$\frak I_{_{(T,\eta)}}^{^{(U,\theta)}}(\varphi)$ a   set of representatives for the set of  
$\bar N_{R\times P}\big(\Delta_\theta (U)\big)\times \bar N_{Q\times P}\big(\Delta_\eta (T)\big)\-$orbits,  for any $f\in \frak I_{_{(T,\eta)}}^{^{(U,\theta)}}(\varphi)$ consider the stabilizer
$\bar N_{Q\times P}\big(\Delta_\eta (T)\big)_{{\rm Im}(f)}$ of ${\rm Im}(f)$ in 
$\bar N_{Q\times P}\big(\Delta_\eta (T)\big)\,,$  so that we get an inclusion and an obvious group homomorphism
$$\eqalign{\varepsilon_f \,\colon \bar N_{Q\times P}\big(\Delta_\eta (T)\big)_{{\rm Im}(f)}&\too 
\bar N_{Q\times P}\big(\Delta_\eta (T)\big)\cr
\delta_f : \bar N_{Q\times P}\big(\Delta_\eta (T)\big)_{{\rm Im}(f)}&\too 
\bar N_{R\times P}\big(\Delta_\theta (U)\big)\cr}
\eqno £2.17.1.$$
fulfilling $\bar a\. f = f\. \delta_f (\bar a)$ for any 
$\bar a\in \bar N_{Q\times P}\big(\Delta_\eta (T)\big)_{{\rm Im}(f)}\,.$
Then, denoting by $\ab^\frak c\,\colon \frak i\Gr\to \Ab$ the {\it contravariant\/} functor mapping any finite group $H$ on its maximal Abelian quotient~$H/[H,H]$ and any injective group homomorphism on the group homomorphism  induced by the {\it transfert\/}, it follows from [8,~Proposition~22.17] that  for any $(T,\eta)\in \frak O_Q$ and any $(U,\theta)\in \frak O_R$ fulfilling condition~£2.16.2 for suitable $u\in Q$ and $n\in P$ we have
$$\big(\Ker (\pi^{_{\rm b}})\big)(x)_{_{(T,\eta)}}^{^{(U,\theta)}} = 
\sum_{f\in \frak I_{_{(T,\eta)}}^{^{(U,\theta)}}(\varphi)}\ab(\delta_f)\circ \ab^\frak c (\varepsilon_f)
\eqno £2.17.2.$$

\bigskip
\noindent
{\bf £3. A filtration for the basic $\F\-$locality\/}

\medskip
£3.1. Let $P$ be a finite $p\-$group, $\F$ a Frobenius $P\-$category and 
$(\tau^{_{\rm b}},\L^{^{\rm b}},\pi^{_{\rm b}})$ the corresponding {\it basic $\F\-$locality\/};
we already know that the {\it contravariant\/} functor
$$\Ker (\pi^{_{\rm b}}) : \L^{^{\rm b}}\too \Ab
\eqno £3.1.1\phantom{.}$$
factorizes troughout  the {\it exterior quotient\/} $\tilde\F$ of~$\F$ (cf. £2.7), so that it defines a 
{\it contravariant\/} functor 
$$\tilde{\frak k}_{\L^{^{\rm b}} } = \tilde\frak k^{^{\rm b}}_\F : \tilde\F\too \Ab
\eqno £3.1.2\phantom{.}$$
which, up to suitable identifications, maps any $\tilde\F\-$morphism 
$\tilde\varphi\,\colon R\to Q$ on the group homomorphism  (cf.~£2.17.2)
$$\tilde\frak k^{^{\rm b}}_\F (\tilde\varphi) = \sum_{(T,\eta)\in \frak O_Q}
\,\sum_{(U,\theta)\in \frak O_R}\,\sum_{f\in \frak I_{_{(T,\eta)}}^{^{(U,\theta)}}(\varphi)}\ab(\delta_f)\circ \ab^\frak c (\varepsilon_f)
\eqno £3.1.3\phantom{.}$$
from ${\displaystyle \bigoplus_{(T,\eta)\in \frak O_Q}}
\ab\Big(\bar N_{Q\times P}\big(\Delta_\eta (T)\big)\Big) $
to ${\displaystyle \bigoplus_{(U,\theta)\in \frak O_R}}
\ab\Big(\bar N_{R\times P}\big(\Delta_\theta (U)\big)\Big)$,
where we set $\frak I_{_{(T,\eta)}}^{^{(U,\theta)}}(\varphi) = \emptyset$ whenever
 condition~£2.16.2 is {\it not\/} fulfillied for any $u\in Q$ and any $n\in P\,.$

\medskip
£3.2. In particular, note that the homomorphism $\tilde\frak k^{^{\rm b}}_\F(\tilde\varphi)$ sends an
element of $\ab\Big(\bar N_{Q\times P}\big(\Delta_\eta (T)\big)\Big)$ to a sum of terms indexed by elements $(U,\theta)$ in $\frak O_R$ such that $U$ is contained in a $P\-$conjugated of $T\,;$ hence, for any subset $\N$ of~$\C_P$ which fulfills
\smallskip
\noindent
£3.2.1\quad {\it any $U\in \C_P$ which is contained in a $P\-$conjugated of 
$T\in \N$ belongs to $\N$,\/}
\smallskip
\noindent 
setting $\frak O_Q^\N = \bigsqcup_{T\in \N}\frak O_Q^T$
for any subgroup $Q$ of $P\,,$ it is quite clear that we get a {\it contravariant\/} subfunctor 
$\tilde\frak k^{^{\N}}_\F\! : \tilde\F\too \Ab$ of $\tilde\frak k^{^{\rm b}}_\F$ sending $Q$ to
$${\displaystyle \bigoplus_{(T,\eta)\in \frak O_Q^\N}}\ab\Big(\bar N_{Q\times P}\big(\Delta_\eta (T)\big)\Big)
\eqno £3.2.2$$
and we consider the corresponding {\it quotient $\F\-$locality\/} 
$\L^{^{\rm b}}/(\tilde\frak k^{^{\N}}_\F\circ \tilde\pi^{^{\rm b}})$ (cf.~£2.8)  --- denoted by 
$(\tau^{^{{\rm b},\N}},\L^{^{{\rm b},\N}},\pi^{^{{\rm b},\N}})$ --- of the {\it basic 
$\F\-$locality\/} $\L^{^{\rm b}}$ above.

\medskip
£3.3. It is quite clear that if $\M$ is another subset of $\C_P$ fulfilling condition~£3.2.1 and containing $\N\,,$  we  have a canonical functor 
$\frak l^{^{\M,\N}}_\F \,\colon \L^{^{{\rm b},\N}} \to  \L^{^{{\rm b},\M}}\,.$ From now on, we fix a proper subset $\N$ of $\C_P$ fulfilling condition~£3.2.1 and, in order to argue by induction on $\vert \C_P - \N\vert\,,$ we also fix a minimal element $U$ in $\C_P - \N\,,$ setting 
$\M = \N\cup \{U\}\,.$ Hence, it makes sense to consider the quotient {\it contravariant\/}  functor  
$$\tilde\frak k^{^U }_\F=\Ker (\bar\frak l^{^{\M,\N}}_\F) =  
\tilde\frak k^{^{\M}}_\F/\tilde\frak k^{^{\N}}_\F  :\tilde\F \too \Ab
\eqno £3.3.1\phantom{.}$$
which only depends on $U$ as we show in~£3.4 and~£3.5 below. More precisely, for any 
$m\in \Bbb N$ let us consider the   subfunctor $p^m\.\id\,\colon \Ab\to \Ab$ of the
 {\it identity functor $\id_{\Ab}$\/} sending any Abelian group $A$ to $p^m\.A$,  setting
$\frak s_m= p^m\.\id/p^{m+1}\.\id\,.$ Then, the key point to prove the main results announced in~£1.6 above is that, {\it for any $m\ge 0$ and any $n\ge 1\,,$ the $n\-$th 
{\it stable cohomology group\/} --- noted  $\Bbb H^{^n}_* (\tilde\F, \frak s_m\circ \tilde\frak k^{^U }_\F)$ {\rm (see [8,~A3.17])} --- of $\tilde\F$ over $\frak s_m\circ \tilde\frak k^{^U }_\F$ vanish\/}; that is to say, that the {\it differential subcomplex\/} in [11,~A2.2], where $\B = \tilde\F$ and $\frak a = \frak s_m\circ \tilde\frak k^{^U}_\F$, and where  we only consider the elements which are stable by $\tilde\F\-$isomorphisms, is exact.

\medskip
£3.4. This vanishing result will follow from Theorem~£3.11 below and from [11, Theorem A5.5]; 
that is to say, with the terminology introduced in [11, 45.1], in Theorem~£3.11 below we  prove
that, for any $m\in \Bbb N\,,$ the {\it contravariant\/} functor 
$\frak s_m\circ \tilde\frak k^{^U }_\F$ above admits indeed a {\it compatible complement\/}. From definition~£3.3.1 above  it is clear that, for any subgroup $Q$ of $P$  
$$\tilde\frak k^{^U}_\F\! (Q) = \bigoplus_{(T,\eta)\in \frak O^\M_Q - \frak O^\N_Q} 
\ab \big(\bar N_{Q\times P}\big(\Delta_{\eta}(U)\big)\big)
\eqno £3.4.1\phantom{.}$$
and then, for any $(T,\eta)\in \frak O^\M_Q - \frak O^\N_Q\,,$ we necessarily have $T = U\,;$ hence, we get
$$\tilde\frak k^{^U}_\F\! (Q) = \bigoplus_{(U,\eta)\in \frak O^U_Q}
\ab \big(\bar N_{Q\times P}\big(\Delta_{\eta}(U)\big)\big)
\eqno £3.4.2\phantom{.}$$
where $\frak O^U_Q\i \F (Q,U)$ is a set of representatives for  
$Q\backslash \F (Q,U)/N_P (U)\,;$ more precisely, the group $Q\times N_P (U)$ acts on 
$\F (Q,U)$ and if $\eta, \eta'\in \F(Q,U)$ are in the same  $Q\times N_P (U)\-$orbit then the conjugation by a suitable element $(u,n)$ in $Q\times N_P (U)$ induces a group isomorphism
$$\ab \big(\bar N_{Q\times P}\big(\Delta_{\eta}(U)\big)\big)\cong 
\ab \big(\bar N_{Q\times P}\big(\Delta_{\eta'}(U)\big)\big)
\eqno £3.4.3\phantom{.}$$
which is clearly independent of the choice of $(u,n)$ fulfilling $\eta' = (u,n)\.\eta\,.$ Consequently, from £3.4.2 we get a canonical isomorphism
$$\tilde\frak k^{^U}_\F \!(Q) \cong \bigg(\bigoplus_{\eta\in \F (Q,U)} 
\ab \Big(\bar N_{Q\times P}\big(\Delta_{\eta}(U)\big)\Big)\bigg)^{Q\times N_P (U)}
\eqno £3.4.4.$$

\medskip
£3.5. Moreover, for any $\tilde\F\-$morphism $\tilde\varphi\,\colon R\to Q\,,$ from £3.1.3 above we still get
$$\tilde\frak k^{^{U}}_\F (\tilde\varphi) = 
\sum_{(U,\eta)\in \frak O_Q^U}\,\sum_{(U,\theta)\in \frak O_R^U}\,\sum_{f\in \frak I_{_{(U,\eta)}}^{^{(U,\theta)}}(\varphi)}\ab(\delta_f)\circ \ab^\frak c (\varepsilon_f)
\eqno £3.5.1;$$
 in this case, it follows from~£2.16 that $\frak I_{_{(U,\eta)}}^{^{(U,\theta)}}(\varphi)$ is empty unless for suitable $u\in Q$ and $n\in P$ we have
$$\Delta_{\varphi\circ\theta} (U) = {}^{(u,n)}\Delta_\eta (U)
\eqno £3.5.2\phantom{.}$$ 
or, equivalently, $n$ belongs to $N_P (U)$ and $\varphi\circ\theta$ belongs to the classe of 
$\eta$ in~$Q\backslash \F (Q,U)/N_P (U)\,;$ in this case we have an injective $R\times P\-$set homomorphism
$$f :(R\times P)\big/\Delta_\theta (U) \too 
{\rm Res}_{\varphi\times {\rm id}_P}\big((Q\times P)\big/\Delta_\eta (U)\big)
\eqno £3.5.3\phantom{.}$$
sending the class of $(1,1)$ to the class of $(u,n)\,.$ Then,  denoting by
$$\varphi_\theta : \bar N_{R\times P}\big(\Delta_{\theta}(U)\big)\too 
\bar N_{Q\times P}\big(\Delta_{\varphi\circ\theta}(U)\big)
\eqno £3.5.4\phantom{.}$$
the group homomorphism induced by $\varphi\times {\rm id}_P$, and by
$$\kappa_{(u,n)}^{\eta,\varphi\circ\theta} :  \bar N_{Q\times P}\big(\Delta_{\eta}(U)\big)
\cong \bar N_{Q\times P}\big(\Delta_{\varphi\circ\theta}(U)\big)
\eqno £3.5.5,$$
 the conjugation by~$(u,n)\,,$ it is quite clear that the image of $\varphi_\theta$ stabilizes the image of $f$ and therefore that $f$ is the unique element of  
$\frak I_{_{(U,\eta)}}^{^{(U,\theta)}}(\varphi)\,,$ that $\delta_f$ is an isomorphism in £2.17.1
 and that we~get [9,~8.8]
$$\ab(\delta_f)\circ \ab^\frak c (\varepsilon_f) = \ab^\frak c \big((\kappa_{(u,n)}^{\eta,\varphi\circ\theta})^{-1}\circ \varphi_\theta) = \ab^\frak c (\varphi_\theta) \circ \ab (\kappa_{(u,n)}^{\eta,\varphi\circ\theta})
 \eqno £3.5.6.$$
Consequently, equality £3.5.1 becomes
$$\tilde\frak k^{^{U}}_\F (\tilde\varphi) = \sum_{(U,\theta)\in \frak O_R^U}
\ab^\frak c (\varphi_\theta) \circ \ab (\kappa_{(u,n)}^{\eta,\varphi\circ\theta})
\eqno £3.5.7\phantom{.}$$
where, for any $(U,\theta)\in \frak O_R^U\,,$ $(U,\eta)$ belongs to $\frak O_Q^U$ and $(u,n)$ fulfills~£3.5.2.

\medskip
£3.6. But, for our purposes,  we need a better description  as follows for the functor 
$\tilde\frak k^{^{U}}_\F$. It is quite clear that we have a functor $\frak n^{_U}_\F\colon \F\to \Gr$ mapping any subgroup $Q$
of $P$ on the direct product of groups
$$\frak n^{_U}_\F (Q) = {\displaystyle \prod_{\eta\in \F (Q,U)}} 
\bar N_{Q\times P}\big(\Delta_{\eta}(U)\big)
\eqno £3.6.1$$
and any $\F\-$morphism $\varphi\,\colon R\to Q$ on the direct product of group homomorphisms (cf.~£3.5.4)
$$\frak n^{_U}_\F (\varphi) =\!\!\! \prod_{\theta\in \F (R,U)}\! \!\varphi_\theta :\!\!\!\! {\displaystyle \prod_{\theta\in \F (R,U)}} \!\! \bar N_{R\times P}\big(\Delta_{\theta}(U)\big)\too \!\!\!
{\displaystyle \prod_{\eta\in \F (Q,U)}} \!\!
\bar N_{Q\times P}\big(\Delta_{\eta}(U)\big)
\eqno £3.6.2;$$
note that, for any $u\in Q$ denoting by $\kappa_{Q,u}\,\colon Q\cong Q$ the conjugation by $u\,,$ the action of $(u,1)\in Q\times N_P (U)$ on $\frak n^{_U}_\F\! (Q)$ coincides with  
$\frak n^U_\F (\kappa_{Q,u})\,.$ 
Similarly, as in~£3.4.3 above, for any $n\in N_P (U)$ the action of $(1,n)\in Q\times N_P (U)$ 
on~$\frak n^{_U}_\F (Q)$ induces obvious isomorphisms
$$\,\overline{\!(1,n)\!}\,^\eta :\bar N_{Q\times P}\big(\Delta_{\eta}(U)\big)\cong \bar N_{Q\times P}\big(\Delta_{\eta\circ\kappa_{U,n^{-1}}} (U)\big)
\eqno £3.6.3\phantom{.}$$
for any $\eta\in \F (Q,U)$; moreover, for any $\theta\in \F (R,U)\,,$ we obviously get
$$\,\overline{\!(1,n)\!}\,^{\varphi\circ \theta} \circ \varphi_\theta = \varphi_{\theta\circ\kappa_{U,n^{-1}}}
\eqno £3.6.4.$$
\eject

\medskip
£3.7.  Consequently, we also get the  functors (cf.~£2.14 and~£2.17)
$$\ab^\frak c \circ \frak n^{_U}_\F : \F\too \Ab^\circ \qq \ab\circ \frak n^{_U}_\F : \F\too \Ab
\eqno £3.7.1\phantom{.}$$ 
which send any subgroup $Q$ of $P$ to 
$$(\ab^\frak c \circ \frak n^{_U}_\F)(Q) = {\displaystyle \bigoplus_{\eta\in \F (Q,U)}} \ab\Big(\bar N_{Q\times P}\big(\Delta_{\eta}(U)\big)\Big) = (\ab\circ \frak n^{_U}_\F)(Q)
\eqno £3.7.2\phantom{.}$$
and we know that $Q\times N_P (U)$ acts on this Abelian $p\-$group; then, it is quite easy to check that
we  have a subfunctor of $\ab^\frak c \circ \frak n^{_U}_\F$ and a quotient functor of $\ab\circ \frak n^{_U}_\F\,,$
respectively  denoted by
$$\frak h^\circ (\ab^\frak c \circ \frak n^{_U}_\F) \,\colon\F\too \Ab^\circ\qq 
\frak h_\circ (\ab \circ \frak n^{_U}_\F) \,\colon\F \too \Ab
\eqno £3.7.3,$$
sending any subgroup $Q$ of $P$ to the subgroup $(\ab^\frak c \circ \frak n^{_U}_\F)(Q)^{Q\times N_P (U)}$ of {\it fixed\/} elements  and to the quotient $(\ab \circ \frak n^{_U}_\F)(Q)_{Q\times N_P (U)}$ of  {\it co-fixed\/} elements of~${\displaystyle \bigoplus_{\eta\in \F (Q,U)}} \ab\Big(\bar N_{Q\times P}\big(\Delta_{\eta}(U)\big)\Big)\,;$  actually, it is easily checked that both factorize through the {\it exterior quotient\/} $\tilde\F$ (cf.~£3.1) yielding respective functors
$$\tilde\frak h^\circ (\ab^\frak c \circ \frak n^{_U}_\F) : \tilde\F\too \Ab^\circ\qq 
\tilde\frak h_\circ (\ab \circ \frak n^{_U}_\F) : \tilde\F \too \Ab
\eqno £3.7.4.$$
In particular, it follows from~£3.4.4 that for  any subgroup $Q$ of $P$ we have
$$\tilde\frak h^\circ (\ab^\frak c \circ \frak n^{_U}_\F)(Q)\cong \tilde\frak k^{^U}_\F (Q)
\eqno £3.7.5.$$

\medskip
£3.8. Actually, we claim that for any $\tilde\F\-$morphism $\tilde\varphi\,\colon R\to Q$ we also have the commutative diagram
$$\matrix{\tilde\frak h^\circ (\ab^\frak c \circ \frak n^{_U}_\F)(Q)&\cong &\tilde\frak k^{^U}_\F\!(Q)\cr
\hskip-50pt{\scriptstyle \tilde\frak h^\circ (\ab^\frak c \circ \frak n^{_U}_\F)(\tilde\varphi)}\big\downarrow&\phantom{\Big\downarrow}&\big\downarrow{\scriptstyle \tilde\frak k^{^U}_\F \!(\tilde\varphi)}\hskip-20pt\cr
\tilde\frak h^\circ (\ab^\frak c \circ \frak n^{_U}_\F)(R)&\cong &\tilde\frak k^{^U}_\F(R)\cr}
\eqno £3.8.1$$
so that the functors $\tilde\frak h^\circ (\ab^\frak c \circ \frak n^{_U}_\F)$ and 
$\tilde\frak k^{^U}_\F \!$ 
from $\tilde\F$ to $\Ab^\circ$ are isomorphic. Indeed, 
$\tilde\frak h^\circ (\ab^\frak c \circ \frak n^{_U}_\F)(\tilde\varphi)$ sends any element 
$$a ={\displaystyle \sum_{\eta\in \F (Q,U)}}  a_\eta\in  \tilde\frak h^\circ 
(\ab^\frak c \circ \frak n^{_U}_\F)(Q)
\eqno £3.8.2,$$
where $ a_\eta$ belongs to $\ab\Big(\bar N_{Q\times P}\big(\Delta_{\eta}(U)\big)\Big)$ 
for any $\eta\in \F(Q,U)\,,$ to the element 
$${\displaystyle \sum_{\theta\in \F (R,U)}} \ab^\frak c (\varphi_{\theta})( a_{\varphi\circ\theta})\in 
\tilde\frak h^\circ (\ab^\frak c \circ \frak n^{_U}_\F)(R)
\eqno £3.8.3;$$ 
then, the commutativity of diagram~£3.8.1 follows from equality~£3.5.7.
\eject

\medskip
£3.9. Moreover,  for any subgroup $Q$ of $P$ we clearly have a canonical group isomorphism
$$\tr_Q : (\ab \circ \frak n^{_U}_\F)(Q)_{Q\times N_P (U)}\cong 
(\ab^\frak c \circ \frak n^{_U}_\F)(Q)^{Q\times N_P (U)}
\eqno £3.9.1\phantom{.}$$
which, for any $\eta\in \F (Q,U)\,,$ maps the class in $(\ab \circ \frak n^{_U}_\F)(Q)_{Q\times N_P (U)}$ of an element $ a_\eta$ of~$\ab\Big(\bar N_{Q\times P}\big(\Delta_{\eta}(U)\big)\Big)$ on the element (cf.~£3.6)
$${\rm tr}_{N_{Q\times P} (\Delta_{\eta}(U))}^{Q\times N_P (U)} ( a_\eta) = \sum_{(u,n)} (u,n)^\eta\. a_\eta
\eqno £3.9.2\phantom{.}$$
in $(\ab^\frak c \circ \frak n^{_U}_\F)(Q)^{Q\times N_P (U)}$ (cf.~£3.7.2), where $(u,n)$ runs over a set of representatives for the quotient set $\big(Q\times N_P (U)\big)\big/N_{Q\times P} (\Delta_{\eta}(U))\,.$ 

\medskip
£3.10. Explicitly, an element $(u,n)$ of $Q\times P$ belongs 
to~$N_{Q\times P} (\Delta_{\eta}(U))$ if and only if we have ${}^n U = U$ and 
$\eta ({}^n v) = {}^u \eta (v)$ for any $v\in U\,;$ in particular, $u$ normalizes $\eta (U)$ and, denoting by $\eta_*\,\colon U\cong \eta (U)$ the isomorphism induced by $\eta\,,$ this element belongs to the converse image $Q_\eta$ of 
$$\F_Q\big(\eta (U)\big)\cap \big(\eta_*\circ \F_P (U)\circ \eta_*^{-1}\big)
\eqno £3.10.1\phantom{.}$$
 in $N_Q\big(\eta (U)\big)\,;$ then, the conjugation by $\eta_*^{-1}$ induces a group homomorphism 
 $\nu_\eta\,\colon Q_\eta\to \F_P (U)\,;$ thus, setting 
$$\Delta^{\nu_\eta} (Q_\eta) = \{(u,\nu_\eta (u)\}_{u\in Q_\eta}\i Q\times \F_P (U)
\eqno £3.10.2,$$
we get the exact sequence
$$1\too \{1\}\times C_P (U)\too N_{Q\times P} \big(\Delta_\eta (U)\big)\too 
\Delta^{\nu_\eta} (Q_\eta) \too 1
\eqno £3.10.3\phantom{.}$$
and, in particular, denoting by $\tilde a_\eta$ the classe of  
$ a_\eta\in \ab\Big(\bar N_{Q\times P} \big(\Delta_{\eta}(U)\big)\Big)$ in 
$(\ab \circ \frak n^{_U}_\F)(Q)_{Q\times N_P (U)}$ and by $Q^{^\eta}\i Q$ a set of representatives 
for $Q/Q_\eta\,,$ since $\{1\}\times C_P (U)$ acts trivially on 
$\ab\Big(\bar N_{Q\times P} \big(\Delta_{\eta}(U)\big)\Big)\,,$  we still get 
$$\eqalign{\tr_Q (\tilde a_\eta) &= \sum_{u\in Q^{^\eta}}  \sum_{\nu\in \F_P (U)} (u,\nu)\.a_\eta\cr
& =  \sum_{\nu\in \F_P (U)} (1,\nu)\.\sum_{u\in Q^{^\eta}}  
 (\ab\circ \frak n^U_\F)(\kappa_{Q,u})(a_\eta)\cr}
\eqno £3.10.4.$$
Finally, for any $m\in \Bbb N$, for short we set
$$\tilde\frak r^{U,\circ}_{\F,m} = 
\frak s_m\!\circ \tilde\frak h^\circ (\ab^\frak c \! \circ \frak n^{_U}_\F) \!\!\qq \!\!
\tilde\frak r_{\F,\circ}^{U,m}  = \frak s_m\!\circ \tilde\frak h_\circ (\ab \circ \frak n^{_U}_\F)
\eqno £3.10.5;$$ 
moreover, for any subgroup $Q$ of $P\,,$ it is clear that $\tr_Q$ induces a group isomorphism
$$\tr_Q^m : \tilde\frak r_{\F,\circ}^{U,m} (Q) \cong \tilde\frak r^{U,\circ}_{\F,m} (Q)
\eqno £3.10.6.$$
\eject

\smallskip
\noindent
{\bf Theorem~£3.11.}{\footnote{\dag}{\cds In [9, Proposition 8.9] the statement and the proof are far from correction}}  
{\it With the notation above, the functor $\tilde\frak r^{U,\circ}_{\F,m} \,\colon \tilde\F\to \Ab^\circ$ admits a compatible complement $(\tilde\frak r^{U,\circ}_{\F,m})^\frak c \,\colon \tilde\F\to \Ab$ sending any $\tilde\F\-$morphism $\tilde\varphi\,\colon R\to Q$ to 
the group  homomorphism 
$$(\tilde\frak r^{U,\circ}_{\F,m})^\frak c (\tilde\varphi) = 
\tr_Q^m\circ\tilde\frak r_{\F,\circ}^{U,m} (\tilde\varphi)  \circ (\tr_R^m)^{-1}
\eqno £3.11.1.$$ \/}
\par
\noindent
{\bf Proof:} It is clear that equalities~£3.11.1 define a functor $\tilde\F\to \Ab$ sending any subgroup $Q$ of $P$ to $\frak s_m\big((\ab^\frak c \circ \frak n^{_U}_\F)(Q)^{Q\times N_P (U)}\big)\,;$ thus, it remains to prove that $(\tilde\frak r^{U,\circ}_{\F,m})^\frak c $ fulfills  the conditions~A5.1.2 and~A5.1.3 in [11].  With the notation in~£3.7 above, assuming that $a ={\displaystyle \sum_{\eta\in \F (Q,U)}}  a_\eta$ belongs to 
$p^m\.\tilde\frak h^\circ (\ab^\frak c \circ \frak n^{_U}_\F)(Q)$ and denoting by $\bar a^m$ its image in $\frak s_m\big((\ab^\frak c \circ \frak n^{_U}_\F)(Q)^{Q\times N_P (U)}\big)\,,$ for condition~A5.1.3 we have to 
compute~$\big((\tilde\frak r^{U,\circ}_{\F,m})^\frak c  (\tilde\varphi)\circ 
\tilde\frak r^{U,\circ}_{\F,m} (\tilde\varphi)\big) (\bar a^m)$ in $\tilde\frak r^{U,\circ}_{\F,m} (Q)\,;$ clearly, this element is the image in $\frak s_m\big((\ab^\frak c \circ \frak n^{_U}_\F)(Q)^{Q\times N_P (U)}\big)$ of 
$$\eqalign{\Big(\tr_Q^m&\circ (\tilde\frak h_\circ (\ab \circ \frak n^{_U}_\F))(\tilde\varphi)\circ (\tr_R^m)^{-1}\Big)\tilde\frak h^\circ (\ab^\frak c \circ \frak n^{_U}_\F) (\tilde\varphi)(a)\big)\cr
&= \Big(\tr_Q^m\circ \tilde\frak h_\circ (\ab \circ \frak n^{_U}_\F) (\tilde\varphi)\circ (\tr_R^m)^{-1}\Big)
\big({\displaystyle \sum_{\theta\in \F (R,U)}} \ab^\frak c (\varphi_{\theta})( a_{\varphi\circ\theta})\big)\cr}\eqno £3.11.2,$$ 
which is equal to zero whenever $\F (R,U)$ is empty.

\smallskip
Otherwise, for any element $(y,n)$ in $R\times N_P (U)\,,$ $(\varphi(y),n)$ belongs to $Q\times N_P (U)$
and therefore, with the obvious action of $R\times N_P (U)$ on $\F (R,U)\,,$ we~have 
$a_{\varphi\circ (y\.\theta\.n^{-1})} = (\varphi(y),n)\. a_{\varphi\circ\theta}\,;$ consequently, this element coincides with (cf.~£2.13 and~£3.9)
$$\sum_{\theta\in \frak O_R^U} {\rm tr}_{N_{Q\times P} 
(\Delta_{\varphi\circ \theta}(U))}^{Q\times N_P (U)} 
\Big(\big(\ab (\varphi_\theta)\circ \ab^\frak c (\varphi_{\theta})\big)( a_{\varphi\circ\theta})\Big)
\eqno £3.11.3\phantom{.}$$
and we know that for any $\theta\in  \frak O_R^U$ we have
$$\big(\ab (\varphi_\theta)\circ \ab^\frak c (\varphi_{\theta})\big)( a_{\varphi\circ\theta}) =
{\big\vert\bar N_{Q\times P} \big(\Delta_{\varphi\circ \theta} (U)\big)\big\vert \over 
\big\vert\bar N_{R\times P} \big(\Delta_\theta (U)\big)\big\vert}\. a_{\varphi\circ\theta}
\eqno £3.11.4;$$
thus, either $\big\vert\bar N_{Q\times P} \big(\Delta_{\varphi\circ \theta} (U)\big)\big\vert\not= 
\big\vert\bar N_{R\times P} \big(\Delta_\theta (U)\big)\big\vert$ and the term
$$ {\rm tr}_{N_{Q\times P} (\Delta_{\varphi\circ \theta}(U))}^{Q\times N_P (U)} 
\Big(\big(\ab (\varphi_\theta)\circ \ab^\frak c (\varphi_{\theta})\big)( a_{\varphi\circ\theta})\Big) 
\eqno £3.11.5\phantom{.}$$
belongs to $p^{m +1}\.\tilde\frak h^\circ (\ab^\frak c \circ \frak n^{_U}_\F)(Q)\,,$
or we have $\varphi (R_\theta) = Q_{\varphi\circ\theta}$ (cf.~£3.10.2).

\smallskip
But,  for any $\eta\in \F (Q,U)$ and any element $(u,n)$ in~$Q\times N_P (U)$ we still have 
$a_{u\.\eta\.n^{-1}} = (u,n)\. a_\eta\,;$ consequently, for any $\eta$ in the set $\tilde\varphi\circ\F (R,U)\,,$  setting 
$$\T_{Q\times N_P (U)} (\eta,\varphi\circ\theta) = \{(u,n)\in Q\times N_P (U)\mid \eta = (u,n)\.(\varphi\circ\theta)\}
\eqno £3.11.6,$$
in the second case we have
$$\eqalign{ {\rm tr}_{N_{Q\times P} (\Delta_{\varphi\circ \theta}(U))}^{Q\times N_P (U)} 
\Big(\big(\ab (\varphi_\theta)&\circ \ab^\frak c (\varphi_{\theta})\big)( a_{\varphi\circ\theta})\Big)\cr
& = \sum_{\eta\in \tilde\varphi\circ\F (R,U)} {\vert \T_{Q\times N_P (U)} (\eta,\varphi\circ\theta)\vert\over 
\big\vert N_{Q\times P} \big(\Delta_{\eta} (U)\big)\big\vert}\. a_\eta\cr}
\eqno £3.11.7.$$
Moreover, for any element $u$ in the {\it transporter\/} $\T_Q \big(\varphi (R),\eta (U)\big)$ (cf.~£2.2),
the following diagram  
$$\matrix{R &\buildrel \varphi'\over\cong& {}^{u\!^{-1}}\!\!\varphi (R)\cr
\hskip-10pt{\scriptstyle \theta'}\uparrow&& \cup\cr
U &\buildrel \eta\over\cong& \eta (U)}
\eqno £3.11.8\phantom{.}$$
determines a pair formed by $\varphi'\in \tilde\varphi$ and by $\theta'$ in the $\{1\}\times N_P (U)\-$orbit
of $\theta$ such that $\eta = \varphi'\circ\theta'$ and therefore it is quite clear that
$$\T_{Q\times N_P (U)} (\eta,\varphi\circ\theta) = \T_Q \big(\varphi (R),\eta (U)\big)\times N_P (U)
\eqno £3.11.9.$$
 Finally, note that the map
 $$\varphi (R)\backslash \T_Q \big(\varphi (R),\eta (U)\big)\too \varphi (R)\backslash Q / \eta (U)
 \eqno £3.11.10\phantom{.}$$
sending the class of $u\in \T_Q \big(\varphi (R),\eta (U)\big)$ to $\varphi (R) u \eta (U)$ is injective and
its image is the set of double classes of cardinal equal to $\vert \varphi (R)\vert\,,$ so that $p$ divides $\big\vert\varphi (R)\backslash \T_Q \big(\varphi (R),\eta (U)\big)\big\vert\,;$ in conclusion, $p$ also divides the quotient $\vert \T_{Q\times N_P (U)} (\eta,\varphi\circ\theta)\vert\big/ \big\vert N_{Q\times P} 
\big(\Delta_{\eta} (U)\big)\big\vert\.$ Consequently, in both cases  we obtain 
$$\big((\tilde\frak r^{U,\circ}_{\F,m})^\frak c (\tilde\varphi)\circ 
\tilde\frak r^{U,\circ}_{\F,m} (\tilde\varphi)\big) (\bar a^{m}) = 0\
\eqno £3.11.11.$$

\smallskip
In order to show condition~A5.1.3 in [11], for any pair of $\tilde\F\-$morphisms 
$\tilde\varphi\,\colon R\to Q$ and $\tilde\psi\,\colon T\to Q$ we have to prove that
$$\tilde\frak r^{U,\circ}_{\F,m} (\tilde\psi)\circ (\tilde\frak r^{U,\circ}_{\F,m})^\frak c (\tilde\varphi) = \sum_{w\in W} (\tilde\frak r^{U,\circ}_{\F,m})^\frak c (\tilde\psi_w)\circ 
\tilde\frak r^{U,\circ}_{\F,m} (\tilde\varphi_w)
\eqno £3.11.12\phantom{.}$$
where, choosing  a pair of representatives   $\varphi$ of $\tilde\varphi$ and $\psi$ of $\tilde\psi\,,$ and  a set of representatives $W\i Q$ for the set of double classes $\varphi (R)\backslash Q/\psi (T)\,,$  for any 
$w\in W$ we set  $S_w = \varphi (R)^w\cap \psi (T)$ and denote by
$$\varphi_w : S_w\too R\quad and\quad \psi_w : S_w\too T
\eqno £3.11.13\phantom{.}$$
the $\F\-$morphisms fulfilling $\varphi\big(\varphi_w (u)\big) = wuw^{-1}$ and 
$\psi\big(\psi_w (u)\big) = u$ for any element $u\in S_w\,.$
\eject

\smallskip
For any $\theta\in \F (R,U)\,,$ let $b_\theta$ be an element of 
$p^m\.\ab\Big(\bar N_{R\times P} \big(\Delta_\theta (U)\big)\Big)$ and denote by 
$\bar b_\theta$ the image of $b_\theta$ in $(\ab \circ \frak n^{_U}_\F)(R)_{R\times N_P (U)}$ (cf.~£3.7.2); thus, $\tr_R^m (\bar b_\theta)$ is an element of $\tilde\frak r^\circ_m (R)$ (cf.~£3.10.6) and we have to compute (cf.~£3.11.1)
$$\eqalign{\big(\tilde\frak r^{U,\circ}_{\F,m} (\tilde\psi)\circ (\tilde\frak r^{U,\circ}_{\F,m})^\frak c (\tilde\varphi)\big) &\big(\tr_R^m (\bar b_\theta)\big)\cr
& = \big(\tilde\frak r^{U,\circ}_{\F,m} (\tilde\psi)\circ \tr_Q^m\circ \tilde\frak r^{U,m}_{\F,\circ}(\tilde\varphi)\big) (\bar b_\theta)\cr} 
\eqno £3.11.14;$$
this element is clearly the image in $\tilde\frak r^{U,\circ}_{\F,m} (T)$
of the element (cf.~£3.9.2)
$$\eqalign{a = \tilde\frak h^\circ (\ab^\frak c \! &\circ \frak n^{_U}_\F)(\tilde\psi)
\bigg(  \tr_Q\Big(\,\overline{\tilde\frak h_\circ (\ab \circ \frak n^{_U}_\F)
(\tilde\varphi) (b_\theta)}\,\Big)\bigg)\cr
&= (\ab^\frak c \! \circ \frak n^{_U}_\F)(\psi)
\Big({\rm tr}_{N_{Q\times P} (\Delta_{\eta}(U))}^{Q\times N_P (U)}
 \big(\ab (\varphi_\theta)(b_\theta)\big)\Big)\cr}
\eqno £3.11.15\phantom{.}$$
where $\,\overline{\tilde\frak h_\circ (\ab \circ \frak n^{_U}_\F) (\tilde\varphi) (b_\theta)}\,$ denotes the image of $\tilde\frak h_\circ (\ab \circ \frak n^{_U}_\F) (\tilde\varphi) (b_\theta)$ in the quotient
$(\ab \circ \frak n^{_U}_\F)(Q)_{Q\times N_P (U)}$ and we set $\eta = \varphi\circ \theta\,.$

\smallskip
Then, as in~£3.10 above, denoting by $Q_\eta$  the converse image of the intersection $\F_Q\big(\eta (U)\big)\cap \big(\eta_*\circ \F_P (U)\circ \eta_*^{-1}\big)$ in $N_Q\big(\eta (U)\big)$ and by $Q^{^\eta}\i Q$ a set of representatives for $Q/Q_\eta\,,$ we have (cf.~£3.10.4)
$${\rm tr}_{N_{Q\times P} (\Delta_{\eta}(U))}^{Q\times N_P (U)}
\big( \ab (\varphi_\theta)(b_\theta)\big) = 
\sum_{u\in Q^{^\eta}} \sum_{\nu\in \F_P (U)} (u,\nu)\. \ab (\varphi_\theta)(b_\theta)
\eqno £3.11.16;$$
but, for any $u\in Q^\eta$ and any $\nu\in \F_P(U)\,,$ the element  
$(u,\nu)\. \ab (\varphi_\theta)(b_\theta)$ belongs to 
$p^m\.\ab\Big(\bar N_{Q\times P} \big(\Delta_{u\.\eta\circ\nu^{-1}} (U)\big)\Big)$ and therefore it follows from  definition~£3.6.2 that the element (cf.~£3.6)
$$\eqalign{a_{u,\nu} &=  (\ab^\frak c \! \circ \frak n^{_U}_\F)(\psi)
\big( (u,\nu)\. \ab (\varphi_\theta)(b_\theta)\big)\cr
&= (1,\nu)\.(\ab^\frak c \! \circ \frak n^{_U}_\F)(\kappa_{Q,u^{-1}}\circ \psi)
\big(\ab (\varphi_\theta)(b_\theta)\big)\cr}
\eqno £3.11.17\phantom{.}$$
is equal to zero unless  $(\kappa_{Q,u^{-1}}\circ \psi)(T)$ contains 
$\eta (U)\,,$ so that there is a unique $\zeta_u \in \F (T,U)$ fulfilling  
$\kappa_{Q,u}\circ\eta = \psi\circ\zeta_u\,;$ in this case, setting $\psi^u = \kappa_{Q,u^{-1}}\circ \psi$ we get
$$a_{u,\nu} =  (1,\nu)\.\ab^\frak c \big((\psi^u)_{\zeta_u}\big)
\big( \ab (\varphi_\theta)(b_\theta)\big)
\eqno £3.11.18;$$
let us denote by $\hat Q^\eta\i Q^\eta$ the subset of $u\in Q^\eta$ fulfilling this condition.

\smallskip
In this situation, note that we have the two injective group homomorphisms
$$\matrix{&&\bar N_{Q\times P} \big(\Delta_\eta (U)\big)\cr
&\hskip-30pt{\varphi_{\theta}\atop}\hskip-5pt\nearrow 
&&\hskip-35pt\nwarrow \hskip-5pt{(\psi^u)_{\zeta_u}\atop}\cr
\bar N_{R\times P} \big(\Delta_\theta (U)\big)\hskip-10pt
&&& \bar N_{T\times P} \big(\Delta_{\zeta_u} (U)\big)\cr}
\eqno £3.11.19;$$
thus, since $\ab^\frak c$ (cf.~£2.17) is a {\it Mackey complement\/} of $\ab\,,$  for any $u\in Q^\eta$ we need to consider the set of double classes 
$$\X_u =N_{\varphi (R)\times P} \big(\Delta_\eta (U)\big)\backslash
N_{Q\times P} \big(\Delta_\eta (U)\big)/N_{\psi^u (T)\times P} \big(\Delta_\eta (U)\big)
\eqno £3.11.20\phantom{.}$$
which, according to the exact sequence~£3.10.3, admits an obvious canonical bijection with the set of double classes 
$\big(\varphi (R)\cap Q_\eta\big)\big\backslash Q_\eta\big/\big(\psi^u (T)\cap Q_\eta\big)\,;$
hence, choosing a set $X_u\i Q_\eta$ of representatives for this last set of double classes, we get
$$\ab^\frak c \big( (\psi_u)_{\zeta_u}\big) \circ \ab (\varphi_\theta) = 
\sum_{x\in X_u}\ab (\psi_{\eta,u,x})\circ \ab^\frak c (\varphi_{\theta, u,x})
\eqno £3.11.21\phantom{.}$$
where  for any $x\in X_u$ we set  $S_{u,x} = \varphi (R)^x\cap \psi (T)^u$ and denote by
$$\eqalign{\varphi_{\theta, u,x} : \bar N_{S_{u,x}\times P} \big(\Delta_\eta (U)\big)
&\too \bar N_{R\times P} \big(\Delta_\theta (U)\big)\cr 
\psi_{\eta,u,x} : \bar N_{S_{u,x}\times P} \big(\Delta_\eta (U)\big)
&\too \bar N_{T\times P} \big(\Delta_{\zeta_u} (U)\big)\cr}
\eqno £3.11.22\phantom{.}$$
the $\F\-$morphisms fulfilling (cf.~£3.10.2)
$$\eqalign{\varphi_\theta\big(\varphi_{\theta, u,x} \,\overline{\!(s,n)\!}\,\big) 
&= \,\overline{\!\big({}^x s,{}^{\hat x} n)\big)\!} = (x,\hat x)\. \,\overline{\!(s,n)\!}\,\cr 
(\psi^u)_{\zeta_u}\big(\psi_{\eta,u,x} \,\overline{\!(s,n)\!}\,\big) 
&= \,\overline{\!(s,n)\!}\,\cr}
\eqno £3.11.23\phantom{.}$$ 
for any element $(s,n)\in N_{S_{u,x}\times P} \big(\Delta_\eta (U)\big)$ and for a choice of $\hat x\in Q_\eta$
lifting $\nu_\eta (x)$ (cf.~£3.10.2); note that the element $(x,\hat x)\in Q_\eta\times P$ normalizes $\Delta_\eta (U)\,.$
Hence, from £3.11.15, £3.11.18, £3.11.19 and~£3.11.21 we obtain
$$a = \!\!\!\sum_{\nu\in \F_P (U)} (1,\nu)\.\sum_{u\in \hat Q^\eta} \,\sum_{x\in X_u}
\!\big(\ab (\psi_{\eta,u,x})\circ \ab^\frak c (\varphi_{\theta, u,x})\big)(b_\theta)
\eqno £3.11.24.$$

\smallskip
On the other hand, we have to prove that the element (cf.~£3.11.12)
$$\bar c =\big(\sum_{w\in W} (\tilde\frak r^{U,\circ}_{\F,m})^\frak c (\tilde\psi_w)\circ 
\tilde\frak r^{U,\circ}_{\F,m} (\tilde\varphi_w)\big)\big(\tr_R^m (\bar b_\theta)\big)
\eqno £3.11.25\phantom{.}$$
is also the image of $a\,.$ But, according to~£3.9, $\tr_R^m (\bar b_\theta)$ is the image in 
$\tilde\frak r^{U,\circ}_{\F,m} (R)$~of
$${\rm tr}_{N_{R\times P} (\Delta_{\theta}(U))}^{R\times N_P (U)} ( b_\theta) = 
\sum_{y\in R^\theta}\, \sum_{\nu\in \F_P (U)} (y,\nu)\. b_\theta
\eqno £3.11.26\phantom{.}$$
in $p^m\.(\ab^\frak c \circ \frak n^{_U}_\F)(R)^{R\times N_P (U)}$ (cf.~£3.10.4) where,  denoting by 
$\theta_*\,\colon U\cong \theta (U)$ the isomorphism induced by $\theta$ and by $R_\theta$  the converse image  of the intersection 
$\F_R\big(\theta (U)\big)\cap \big(\theta_*\circ \F_P (U)\circ \theta_*^{-1}\big)$ in 
 $N_R\big(\theta (U)\big)\,,$ $R^\theta\i R$ is a set of representatives for $R/R_\theta\,.$ Note that,  according to~£3.10 we have $R_\theta = \varphi^{-1}(Q_\eta)\,.$
 \eject
 
 \smallskip
 In particular, for any $w\in W\,,$ the element
 $\tilde\frak r^{U,\circ}_{\F,m} (\tilde\varphi_w)\big(\tr_R^m (\bar b_\theta)\big)$
is the image in $\tilde\frak r^{U,\circ}_{\F,m} (S_w)$ of the element (cf.~£3.10.5)
$$d_w =\sum_{y\in R^\theta}\, \sum_{\nu\in \F_P (U)} (\ab^\frak c \circ \frak n^{_U}_\F)(\varphi_w) \big( (y,\nu)\. b_\theta\big)
\eqno £3.11.27;$$
but, as above, for any $y\in R^\theta$ and any $\nu\in \F_P(U)\,,$ the element  
$(y,\nu)\. b_\theta$ belongs to 
$p^m\.\ab\Big(\bar N_{R\times P} \big(\Delta_{y\.\theta\circ\nu^{-1}} (U)\big)\Big)$ and therefore it follows from  definition~£3.6.2 that the element
$$d_{w,y,\nu} = (\ab^\frak c \circ \frak n^{_U}_\F)(\varphi_w) \big( (y,\nu)\. b_\theta\big)
\eqno £3.11.28\phantom{.}$$
is equal to zero unless $\varphi_w(S_w)$ contains ${}^y \theta (U)\,,$ so that there is a unique 
$\xi_{w,y} \in \F (S_w,U)$ fulfilling  $\kappa_{R,y}\circ\theta = \varphi_w \circ\xi_{w,y}\,,$ which forces the equality $\kappa_{Q,w^{-1}\varphi (y)}\circ \eta = \iota_{S_w}^Q\circ \xi_{w,y}\,;$ in this case, we have
$$d_{w,y,\nu} =  (1,\nu)\.\ab^\frak c \big((\kappa_{R,y^{-1}}\circ\varphi_w )_{\xi_{w,y}}\big) (b_\theta)
\eqno £3.11.29;$$
let us denote by $\hat R^{\theta,w}\i R^\theta$ the subset of $y\in R^\theta$ fulfilling this condition; thus, we get
$$d_w = \sum_{\nu\in \F_P (U)} (1,\nu)\.\sum_{y\in \hat R^{\theta,w}} 
\ab^\frak c \big((\kappa_{R,y^{-1}}\circ\varphi_w )_{\xi_{w,y}}\big) (b_\theta)
\eqno £3.11.30.$$

\smallskip
Moreover, for any $y\in R^\theta$ fulfilling the above condition and any $s\in S_w\,,$ the product 
$\varphi_w(s)\. y$ still fulfills this condition and we have
$$\eqalign{\varphi_w \circ \xi_{w,\varphi_w(s)\. y} 
&= \kappa_{R,\varphi_w(s)\.y}\circ\theta = \kappa_{R,\varphi_w (s)}\circ \kappa_{R,y}\circ\theta\cr
&=\kappa_{R,\varphi_w (s)}\circ \varphi_w \circ\xi_{w,y} = \varphi_w\circ \kappa_{S_w,s}\circ \xi_{w,y}\cr}
\eqno £3.11.31,$$
so that we get $\xi_{w,\varphi_w(s)\. y}  = \kappa_{S_w,s}\circ \xi_{w,y}\,;$ in particular, $\varphi_w (S_w)$ has an action on
$\hat R^{\theta,w}$ and, choosing a set of representatives $\hat Y^{\theta,w}\i \hat R^{\theta,w}$
for the set of $\varphi_w (S_w)\-$orbits, the element $d_w$ above is also equal to
$$ \sum_{\nu\in \F_P (U)} (1,\nu)\.\sum_{y\in \hat Y^{\theta,w}} \, \sum_{s\in S_w^{\theta,y}}
\ab^\frak c \big((\kappa_{R,(\varphi_w (s)\.y)^{-1}}\circ
\varphi_w )_{\xi_{w,\varphi_w (s)\.y}}\big) (b_\theta)
\eqno £3.11.32$$
where for any $y\in \hat Y^{\theta,w}\,,$ setting $S_{w,\theta,y} = \varphi_w^{-1}\big((R_\theta)^y\big)\,,$ $S_w^{\theta,y}\i S_w$ is a set of representatives for~$S_w/S_{w,\theta,y}\,;$ but, it is quite clear that
$$\eqalign{\ab^\frak c \big(&(\kappa_{R,(\varphi_w (s)\.y)^{-1}}\circ
\varphi_w )_{\xi_{w,\varphi_w (s)\.y}}\big)\cr
&= \ab^\frak c \big((\kappa_{R,y^{-1}}\circ\varphi_w\circ 
\kappa_{S_w,s^{-1}} )_{\xi_{w,\varphi_w (s)\.y}}\big)\cr
&= \ab \big( (\kappa_{S_w,s})_{\xi_{w,y}}\big)\circ 
\ab^\frak c \big((\kappa_{R,y^{-1}}\circ\varphi_w )_{\xi_{w,y}}\big)\cr}
\eqno 3.11.33;$$
hence,  setting $\varphi_w^y = \kappa_{R,y^{-1}}\circ\varphi_w $ and   denoting by 
$\overline{\ab^\frak c \big((\varphi_w^y )_{\xi_{w,y}}\big) (b_\theta)}$ the image of 
$\ab^\frak c \big((\varphi_w^y )_{\xi_{w,y}}\big) (b_\theta)$ in the quotient
$(\ab \circ \frak n^{_U}_\F)(S_w)_{S_w\times N_P (U)}\,,$  according to~£3.10.4 we easily obtain
$$\eqalign{\sum_{\nu\in \F_P (U)} (1,\nu)\. \sum_{s\in S_w^{\theta,y}}
&\Big(\ab^\frak c \big((\kappa_{R,(\varphi_w (s)\.y)^{-1}}\circ
\varphi_w )_{\xi_{w,\varphi_w (s)\.y}}\big)\Big) (b_\theta)\cr
& = \tr_{S_w}^m \Big(\,\overline{\ab^\frak c \big((\varphi_w^y)_{\xi_{w,y}}\big) (b_\theta)}\,\Big)\cr}
\eqno £3.11.34.$$

\smallskip
Consequently, it follows from definition~£3.11.1 that we have (cf.~£3.11.25)
$$\bar c = \sum_{w\in W} \big(\tr_T^m \circ \tilde\frak r_\circ^m (\tilde\psi_w)\big)
\Big(\sum_{y\in \hat Y^{\theta,w}} \,\overline{\ab^\frak c \big((\varphi_w^y)_{\xi_{w,y}}\big) 
(b_\theta)}\,\Big)
\eqno £3.11.35;$$
this element is clearly the image in $\tilde\frak r_\circ^m (T)$ of the element
$$\hskip-2pt\eqalign{\sum_{w\in W} \sum_{y\in \hat Y^{\theta,w}} \!\!\!\! {\rm tr}_{N_{T\times P} 
(\Delta_{\zeta_{w,y}} (U))}^{T\times N_P (U)} &\bigg(\!\Big(\ab \big((\psi_w)_{\xi_{w,y}}\big)\!\circ\!
\ab^\frak c \big((\varphi_w^y)_{\xi_{w,y}}\big)\Big) (b_\theta)\!\bigg) \cr
&\Vert\cr
\sum_{w\in W} \sum_{y\in \hat Y^{\theta,w}}  \sum_{\nu\in \F_P (U)}\! \!\!\!(1,\nu)\.\!\!\!\!\sum_{z\in Z^{w,y}} \!\!&\bigg(\!\Big(\ab \big(({}^z\psi_w)_{\xi_{w,y}}\big)\!\circ\!
\ab^\frak c \big((\varphi_w^y)_{\xi_{w,y}}\big)\Big) (b_\theta)\!\bigg)\cr
&\Vert\cr
\sum_{\nu\in \F_P (U)}\! \!\!\!(1,\nu)\.\!\!\!\!\sum_{w\in W}\! \sum_{y\in \hat Y^{\theta,w}} \sum_{z\in Z^{w,y}} 
\!\!&\bigg(\!\Big(\ab \big(({}^z\psi_w)_{\xi_{w,y}}\big)\!\circ\!
\ab^\frak c \big((\varphi_w^y)_{\xi_{w,y}}\big)\Big) (b_\theta)\!\bigg)\cr}
\eqno £3.11.36$$
in $p^m\. (\ab\circ \frak n^U_\F)(T)^{T\times N_P (U)}$ where,  for any $w\in W$ and any 
$y\in \hat Y^{\theta,w}\,,$ setting~$\zeta_{w,y} = \psi_w\circ \xi_{w,y}$ and denoting by 
$(\zeta_{w,y})_*\,\colon U\cong \zeta_{w,y} (U)$ the isomorphism induced by $\zeta_{w,y}$ and by $T_{w,y}$  the converse image of the intersection 
$\F_T\big(\zeta_{w,y} (U)\big)\cap \big((\zeta_{w,y})_*\circ \F_P (U)\circ (\zeta_{w,y})_*^{-1}\big)$ in $N_T\big(\zeta_{w,y} (U)\big)\,,$  we choose as above a set of representatives $Z^{w,y}$ for 
$T/T_{w,y}$ and, for any $z\in Z^{w,y}\,,$ we set~${}^z \psi_w = \kappa_{T,z}\circ \psi_w\,.$
 
\smallskip
Finally, we claim that this element $\bar c$ coincides with $a$ in~£3.11.24 above; that is to say, considering the sets
$$X = \bigsqcup_{u\in \hat Q^\eta} \{u\}\times X_u \qq Z = \bigsqcup_{w\in W}\,\{w\}\times\!\! \bigsqcup_{y\in \hat Y^{\theta,w}} \{y\}\times Z^{w,y}
\eqno £3.11.37,$$
in $p^m\. (\ab\circ \frak n^U_\F)(T)^{T\times N_P (U)}$ we have to prove the equality
$$\eqalign{&\sum_{(u,x)\in X} \,\sum_{\nu\in \F_P(U)}
\!(1,\nu)\.\big(\ab (\psi_{\eta,u,x})\circ \ab^\frak c (\varphi_{\theta, u,x})\big)(b_\theta)\cr
&= \sum_{(w,y,z)\in Z}\,\sum_{\nu\in \F_P(U)}
\!\!(1,\nu)\.\bigg(\!\Big(\ab \big(({}^z\psi_w)_{\xi_{w,y}}\big)\!\circ\!
\ab^\frak c \big((\varphi_w^y)_{\xi_{w,y}}\big)\Big) (b_\theta)\!\bigg)\cr}
\eqno £3.11.38;$$
actually, we will define a bijection between $X$ and $Z$ such that the cor-responding terms in both sums coincide with each other.

\smallskip
Indeed, for any $w\in W\,,$ any $y\in \hat Y^{\theta, w}$ and any $z\in Z^{w,y}$ let us consider the element
$\varphi (y)^{-1}w\psi (z)^{-1}$ of $Q\,;$ this element certainly belongs to $Q_\eta u^{-1}$ for some $u$ in $Q^\eta$ so that 
we have $\varphi (y)^{-1}w\psi (z)^{-1} = vu^{-1}$ for some $v$ in~$Q_\eta\,;$ but, since $y$ belongs to $\hat Y^w\,,$ 
$\varphi_w (S_w)$ contains ${}^y \theta (U)$ and therefore ${}^w S_w$ contains ${}^{\varphi (y)} \eta (U)$ 
(cf.~£3.11.13); thus, ${}^{\psi (z)w^{-1}\varphi (y)}\eta (U)$ is contained in~$\psi (T)$ and, since 
$Q_\eta\i N_Q \big(\eta (U)\big)$ (cf.~£3.9), ${}^u \eta (U) = {}^{\psi (z)w^{-1}\varphi (y)v}\eta (U)$ is also contained in~$\psi (T)\,,$ so that~$u$ 
belongs to $\hat Q^\eta$ (cf.~£3.11.18). Moreover, the double class of $v$ in $\big(\varphi (R)\cap Q_\eta\big)\big\backslash Q_\eta\big/\big(\psi^u (T)\cap Q_\eta\big)$ determines an element $x$ in $X_u$ such that we have 
$v = \varphi (r)xu^{-1}\psi (t)u$ for some $r\in R_\theta$ and some $t\in T$ fulfilling $\psi (t)^u \in Q_\eta\,,$ so that we get
$$\varphi (y)^{-1}w\psi (z)^{-1} = \varphi (r)xu^{-1}\psi (t)
\eqno £3.11.39.$$
Thus, we obtain a map from $Z$ to $X$ sending $(w,y,z)$ to $(u,x)\,.$

\smallskip
Moreover, with the same notation, setting 
$$q = \psi (tz)^{-1}u = w^{-1}\varphi (yr)x
\eqno £3.11.40,$$
it is clear that the automorphism $\kappa_{Q,q}$ of $Q$ (cf.~£3.6) maps $S_{u,x}$ onto $S_w$  inducing a group isomorphism $\chi\,\colon S_{u,x}\cong S_w\,;$ hence, since we have (cf.~£3.10.3)
$$\kappa_{Q,q}\circ \eta = \kappa_{Q,w^{-1}\varphi(y)}\circ \eta \circ \nu_\eta \big(\varphi(r)x\big) = 
\iota_{S_w}^Q\circ \xi_{w,y}\circ  \nu_\eta \big(\varphi(r)x\big)
\eqno £3.11.41,$$
we get the group isomorphism (cf.~£3.5.4)
$$\chi_\eta : \bar N_{S_{u,x}\times P} \big(\Delta_{\eta} (U)\big) \cong 
\bar N_{S_w\times P} \big(\Delta_{\xi_{w,y}\circ \nu_\eta (\varphi(r)x)} (U))
\eqno £3.11.42.$$
Then, we claim that
$$\eqalign{\big(1,\nu_\eta(\varphi (r)x)^{-1}\big)\.\ab^\frak c (\varphi_{\theta, u,x}) &= 
\ab^\frak c (\chi_{\eta})\circ \ab^\frak c \big((\varphi_w^y)_{\xi_{w,y}}\big)\cr
\big(1,\nu_\eta (\varphi(r)x\psi^u (t)\big)\.\ab (\psi_{\eta,u,x}) &= \ab \big(({}^z\psi_w)_{\xi_{w,y}} \big)\circ \ab (\chi_\eta)\cr}
\eqno £3.11.43.$$

\smallskip
Indeed, for any $(s,n)\in N_{S_{u,x}\times P} \big(\Delta_\eta (U)\big)$ it is easily checked that we have (cf.~£3.11.23)
$$\eqalign{\big(\varphi_{\theta\circ \nu_\eta (\varphi(r)x) } &\circ (\varphi_w^y)_{\xi_{w,y}\circ \nu_\eta (\varphi(r)x)}\circ \chi_\eta\big)\,\overline{\!(s,n)\!}\,\cr
& = \big(\varphi_{\theta\circ \nu_\eta (\varphi(r)x)}\circ (\varphi_w^y)_{\xi_{w,y}\circ \nu_\eta (\varphi(r)x)}\big)
\,\overline{\!({}^q s, n)\!}\cr
& = \,\overline{\!({}^{\varphi (r)x} s,n)\!} = (\varphi(r)x,{\rm id}_P)\.\,\overline{\!(s,n)\!}\,\cr
&= (\varphi (r),\nu_\eta (x)^{-1})\. \big(\varphi_\theta \circ \varphi_{\theta, u,x}\big) 
\,\overline{\!(s,n)\!}\,\cr
&= \big(\varphi_{\theta\circ \nu_\eta(\varphi (r)x)} \circ \varphi_{\theta, u,x}\big) 
\,\overline{\!(s,n)\!}\,\cr}
\eqno £3.11.44\phantom{.}$$
since it follows from~£3.10.3 that $\kappa_{R,r}\circ \theta = \theta\circ \nu_\theta (r) = \theta\circ \nu_\eta \big(\varphi (r)\big)\,;$ thus, since the homomorphism $\varphi_{\theta\circ \nu_\eta (\varphi(r)x)}$ is injective, we get
$$(\varphi_w^y)_{\xi_{w,y}\circ \nu_\eta (\varphi(r)x)}\circ \chi_\eta = \varphi_{\theta, u,x}
\eqno £3.11.45\phantom{.}$$
\eject
\noindent
 and, according to~£3.6.4 above, we still get
 $$(\varphi_w^y)_{\xi_{w,y}}\circ \chi_\eta = \big(1,\nu_\eta(\varphi (r)x)\big)\.\varphi_{\theta, u,x}
 \eqno £3.11.46.$$

\smallskip
Similarly, since $q = \big(\psi (z)^{-1}u\big)\psi^u (t)^{-1}\,,$  we have (cf.~£3.11.23)
$$\eqalign{\big((\psi^u)_{{}^z\zeta_{w,y}\circ \nu_\eta (\varphi (r)x)} 
&\circ ({}^z\psi_w)_{\xi_{w,y}\circ \nu_\eta (\varphi (r)x)} \circ \chi_\eta\big)\,\overline{\!(s,n)\!}\,\cr
&= \big((\psi^u)_{{}^z\zeta_{w,y}\circ \nu_\eta (\varphi (r)x)} \circ ({}^z\psi_w)_{\xi_{w,y}} \big)\,\overline{\!({}^q s,n)\!}\,\cr
& = \,\overline{\!({}^{\psi^u (t^{-1})} s,n)\!}\, \cr
&= \big(\psi^u (t^{-1}),{\rm id}_P\big)\.\big((\psi^u)_{\zeta_u}\circ \psi_{\eta, u,x}\big) 
\,\overline{\!( s,n)\!}\,\cr
&= \big((\psi^u)_{\zeta^t_u}\circ \big(t^{-1},{\rm id}_P\big)\.\psi_{\eta, u,x}\big) 
\,\overline{\!( s,n)\!}\,\cr}
\eqno £3.11.47\phantom{.}$$
where, as in~£3.11.23 above, $\zeta_u\in \F (T,U)$ is the unique element fulfilling 
$\eta = \psi^u\circ\zeta_u\,;$ but, it is easily checked that 
$\zeta^t_u = {}^z\zeta_{w,y}\circ \nu_\eta (\varphi (r)x)\,;$ thus,  since the homomorphism 
$(\psi^u)_{\zeta^t_u}$ is injective, we get
$$ ({}^z\psi_w)_{\xi_{w,y}\circ \nu_\eta (\varphi (r)x)} \circ \chi_\eta
 = \big(t^{-1},{\rm id}_P\big)\. \psi_{\eta, u,x}
\eqno £3.11.48\phantom{.}$$
 and, according to~£3.6.4 above, we still get
 $$ ({}^z\psi_w)_{\xi_{w,y}} \circ \chi_\eta
 = \big(t^{-1},\nu_\eta (\varphi (r)x)\big)\. \psi_{\eta, u,x}
 \eqno £3.11.49;$$
 moreover, since $\psi^u (t)\in Q_\eta,$  in~£3.10.3 above $\big(t, \nu_\eta (\psi^u (t)\big)$ is the image of an element of  $N_{T\times P} \big(\Delta_{\zeta_u} (U)\big)\,;$ hence, it acts trivially over 
$\ab \Big(\bar N_{T\times P} \big(\Delta_{\zeta_u} (U)\big)\Big)$ and therefore we obtain
$$\ab\big(({}^z\psi_w)_{\xi_{w,y}}\big) \circ \ab\big(\chi_\eta\big) = 
\big(1,\nu_\eta (\varphi (r)x \psi^u (t)\big)\. \ab\big(\psi_{\eta, u,x}\big)
\eqno £3.11.50.$$

\smallskip
Finally, since $\ab^\frak c (\chi_\eta) = \ab (\chi_\eta)^{-1}\,,$ the composition of both equalities in~£3.11.43 yields
$$\eqalign{\big(1,\nu_\eta (\psi^u (t)\big)\. \ab (\psi_{\eta,u,x})&\circ \ab^\frak c (\varphi_{\theta\, u,x})\cr
& = \ab \big(({}^z\psi_w)_{\xi_{w,y}} \big)\circ  \ab^\frak c \big((\varphi_w^y)_{\xi_{w,y}}\big)\cr}
\eqno £3.11.51\phantom{.}$$
and therefore in~£3.11.38 we get
$$\eqalign{\sum_{\nu\in \F_P (U)} \!\!\!\Big(\ab (\psi_{\eta,u,x})&\circ 
\ab^\frak c (\varphi_{\theta\, u,x})\Big)(b_\theta)\cr
& = \!\!\!\sum_{\nu\in \F_P (U)} \Big(\ab \big(({}^z\psi_w)_{\xi_{w,y}} \big)
\circ  \ab^\frak c \big((\varphi_w^y)_{\xi_{w,y}}\big)\Big)(b_\theta)\cr}
\eqno £3.11.52.$$
\eject

\smallskip
Conversely, for any $u\in \hat Q^\eta$ and any $x\in X_u$ let us consider the element $w\in W$ determined by the double class of $xu^{-1}$ in $\varphi (R)\backslash Q/\psi (T)\,,$ so that we have
$$xu^{-1} = \varphi (y)^{-1}w \psi (z)
\eqno £3.11.53\phantom{.}$$
for suitable $y\in R$ and $z\in T\,;$ then, with the notation above, we claim that $\varphi_w (S_w)$ contains
${}^y\theta (U)$ or, equivalently, that ${}^w S_w = \varphi (R)\cap {}^w \psi (T)$ contains 
${}^{\varphi (y)} \eta (U)\,.$ Indeed, since $\theta (U)$ is contained in $R\,,$ it is clear that ${}^{\varphi (y)} \eta (U)$ is contained in $\varphi (R)\,;$ it remains to prove that $\eta(U)$ is contained in 
${}^{\varphi (y)^{-1}w} \psi (T)$ or, equivalently, in ${}^{xu^{-1}} \psi (T)\,;$ but, $x$ normalizes $\eta (U)$
and $\eta = \psi^u\circ \zeta_u\,,$ so that $\eta (U)$ is contained in $\psi^u (T)\,;$ this proves the claim.

\smallskip
Consequently, from the very definitions  of $\hat R^{\theta,w}\,,$ of $\hat Y^{\theta,w}$ and of 
$\hat S^{\theta,\hat y}_w$  above, we actually have $y = \varphi_w (s)\hat y r$ for a unique 
$\hat y\in \hat Y^{\theta,w}\,,$ a unique $s\in \hat S^{\theta,\hat y}_w$ and a unique $r\in R_\theta\,;$ now, the equality~£3.11.40 becomes
$$\eqalign{\varphi (r)x u^{-1} &= \varphi (\hat y)^{-1}(\varphi\circ \varphi_w)(s^{-1})w \psi (z)\cr
&= \varphi (\hat y)^{-1}w  s^{-1}\psi (z)\cr}
\eqno £3.11.54\phantom{.}$$
and, since $s\in S_w\i \psi (T)\,,$ there exist a unique $\hat z\in Z^{w,\hat y}$ and a unique $t\in T_{w,\hat y}$
fulfilling $s^{-1}\psi (z) = \psi (\hat z t^{-1})\,,$ so that equality~£3.11.55 becomes
$$\varphi (r)xu^{-1} \psi (t) = \varphi (\hat y)^{-1}w\psi(\hat z)
\eqno £3.11.55.$$
Thus, we obtain a map from $X$ to $Z$ sending $(u,x)$ to $(w, \hat y,\hat z)$ which is clearly the inverse of the map from $Z$ to $X$ defined above. We are done.

\medskip
£3.12. For the next result, we borrow the notation from~£A5 in [11]. Recall that in~£3.10.5 above, for any
$m\in \Bbb N$  we actually define the functors
$$\tilde\frak r^{U,\circ}_{\F,m} : \tilde\F\too \O\-\mod^\circ\qq 
\tilde\frak r^{U,m}_{\F,\circ} : \tilde\F\too \O\-\mod
\eqno £3.12.1.$$

\bigskip
\noindent
{\bf Corollary~£3.13.} {\it Let $\tilde\G$ be a subcategory of $\tilde\F$ having the same objects, only having $\tilde\G\-$isomorphisms and containing all the 
$\tilde\F_{\!P}\-$isomorphisms. Then, with the notation above, for any $m\in \Bbb N$ and any $n\ge 1$ we have 
$$\Bbb H^n_{\tilde\G} (\tilde\F,\tilde\frak r^{U,\circ}_{\F,m}) = \{0\}\,. $$\/}
\par
\noindent
{\bf Proof:} It is an immediate consequence of Theorems~£3.11 above and Theorem A5.5 in [11].

\bigskip
\noindent
{\bf £4. Existence and uniqueness of the perfect $\F\-$locality\/}

\medskip
£4.1. As in~£3.1 above, let $P$ be a finite $p\-$group, $\F$ a Frobenius $P\-$category and 
$(\tau^{_{\rm b}},\L^{^{\rm b}},\pi^{_{\rm b}})$ the corresponding {\it basic $\F\-$locality\/}. 
Recall that we have a {\it contravariant\/} functor [8, Proposition 13.14]
$${\frak c}^\frak f_\F : \F \too \Ab
\eqno £4.1.1\phantom{.}$$
\eject
\noindent
mapping any subgroup $Q$ of $P$ fully centralized in $\F$ on $C_P (Q)/F_{C_\F (Q)}$,
where $F_{C_\F (Q)}$ denotes the {\it $C_{\F} (Q)\-$focal subgroup\/} of $C_P (Q)$ [8, 13.1], and any $\F\-$morphism $\varphi\,\colon R\to Q$ between subgroups of~$P$ fully centralized
in $\F\,,$  on the group homomorphism
$$C_P (Q)/F_{C_\F (Q)}\too C_P (R)/F_{C_\F (R)}
\eqno £4.1.2\phantom{.}$$
induced by an $\F\-$morphism [8, 2.8.2]
$$\zeta : \varphi (R)\.C_P (Q)\too R\.C_P (R)
\eqno £4.1.3\phantom{.}$$ 
fulfilling $\zeta\big(\varphi (v)\big) = v$ for any $v\in R\,.$ Actually, it is easily checked that
this {\it contravariant\/} functor factorizes through the {\it exterior quotient\/} $\tilde\F$ inducing
a new {\it contravariant\/} functor 
$${\tilde\frak c}^\frak f_\F : \tilde\F \too \Ab
\eqno £4.1.4\phantom{.}$$

\bigskip
\noindent
{\bf Proposition £4.2.} {\it The structural functor $\tau^{_{\rm b}} : \T_P \to \L^{^{\rm b}} $ induces a  natural map $\hat\tau^{_{\rm b}}$  from ${\tilde\frak c}^\frak f_\F$ to 
$\tilde{\frak k}^{^{\rm b}}_\F$.\/}

\medskip
\noindent
{\bf Proof:} For any subgroup $Q$ of $P$, the functor $\tau^{_{\rm b}}$ induces a group homomorphism $\tau^{_{\rm b}}_{_Q}$ from $N_P (Q)$ to $\L^{^{\rm b}} (Q)$ which clearly maps $C_P (Q)$ in   $\big(\Ker (\pi^{_{\rm b}}) \big)(Q)$; we claim that this correspondence defines a  natural map (cf. £3.1.2)
$$\hat\tau^{_{\rm b}} : {\tilde\frak c}^\frak f_\F \too \tilde{\frak k}^{^{\rm b}}_\F
\eqno £4.2.1.$$
\smallskip
First of all, we claim that $\tau^{_{\rm b}}_{_Q}$ maps the {\it $C_{\F} (Q)\-$focal subgroup\/}  above on the trivial subgroup of $\L^{^{\rm b}} (Q)$; we may assume that $Q$ is fully centralized in $\F$ and then we know that $F_{C_\F (Q)}$ is generated by the elements $u^{-1}\theta (u)$ where $u$ runs over any subgroup $T$ of $C_P (Q)$ and $\theta$ runs over $\F (T.Q)$ stabilizing $T$ and acting trivially on $Q$ [8, 13.1]; but, according to £2.12 above, $\theta$ can be lifted to 
$\hat z\in N_G (T.Q)$ normalizing $T$ and centralizing $Q\,;$ hence, the element 
$u^{-1}\theta(u) = [u,\hat z^{-1}]$ belongs to $[C_G(Q),C_G(Q)]$ and
therefore it has indeed a trivial image in $\L^{^{\rm b}} (Q)\,;$ consequently,  the canonical homomorphism
$$C_P(Q) \i C_G (Q)\to {\rm Ker} (\pi_{_Q})
\eqno £4.2.2\phantom{.}$$
facorizes through a group homomorphism 
$\hat \tau_{_Q}^{_{\rm b}} : \tilde{\frak c}^\frak f_\F (Q)\to {\rm Ker}(\pi^{_{\rm b}}_{_Q})\,.$

\smallskip
In order to prove the naturality of this correspondence, let $x : R\to Q$ be an 
$\L^{^{\rm b}}\-$morphism  between subgroups of $P$ fully centralized in $\F$ and set 
$\varphi = \pi^{_{\rm b}}_{Q,R} (x)\,;$ it follows from  [8, 2.8.2] that there exists an 
$\F\-$morphism $\zeta$\ from  $\varphi (R)\.C_P (Q)$ to $R\.C_P (R)$ fulfilling $\zeta\big(\varphi (v)\big) = v$ for any $v\in R\,;$ then, $\zeta$ can be lifted to an $\L^{^{\rm b}}\-$morphism
$$y : \varphi (R)\.C_P (Q)\too R\.C_P (R)
\eqno £4.2.3\phantom{.}$$ 
\eject
\noindent 
fulfilling $\big( \pi^{_{\rm b}}_{R.C_P(R),\varphi(R)\.¢_P (Q)} (y)\big)\big(\varphi (v)\big) = v$ for any $v\in R\,;$ in particular, by the {\it divisibility\/}  of $\L^{^{\rm b}}\!,$ $y$ induces an 
$\L^{^{\rm b}}\-$isomorphism $y_{_R} : \varphi(R)\cong R$ and then, setting $x y_{_R} = z\,,$\  the $\L^{^{\rm b}}\-$morphism $z : \varphi (R)\to Q$  fulfills 
  $\pi^{_{\rm b}}_{Q,\varphi (R)}(z) = \iota_{_{\varphi(R)}}^{_Q}$ (§cf. £2.4); consequently, we easily get the following commutative diagram
$$\matrix{{\tilde\frak c}^\frak f_\F(Q) &\buildrel \hat \tau_{_Q}^{_{\rm b}}\over \too &{\rm Ker} (\pi_{_Q})\cr
\hskip-25pt{\scriptstyle {\tilde\frak c}^\frak f_\F (\tilde\varphi) }\Big\downarrow&& \Big\downarrow{\scriptstyle\tilde{\frak k}^{^{\rm b}}_\F (\tilde\varphi)}\hskip-20pt\cr
\tilde{\frak c}^\frak f_\F(R) &\buildrel \hat \tau_{_R}^{_{\rm b}}\over\too &{\rm Ker} (\pi_{_R})\cr}
\eqno £4.2.4.$$
We are done.

\medskip
£4.3. The image $\hat\tau^{_{\rm b}} (\tilde{\frak c}^\frak f_\F)$ of 
$\hat\tau^{_{\rm b}}$ is a subfunctor of $\tilde{\frak k}^{^{\rm b}}_\F$ and therefore, by~£2.8 above, it determines a quotient $\F\-$locality $\widetilde{\L^{^{\rm b}} } = \L^{^{\rm b}}\big/ 
\big(\hat\tau^{_{\rm b}} ({\tilde\frak c}^\frak f_\F)\circ\tilde \pi^{_{\rm b}}\big)$
of $\L^{^{\rm b}}\!$ (cf.~£2.3); we denote by 
$$\widetilde{\tau^{_{\rm b}}} : \T_P \too \widetilde{\L^{^{\rm b}}}\qq \widetilde{\pi^{_{\rm b}} }: \widetilde{\L^{^{\rm b}}}\too \F
\eqno £4.3.1\phantom{.}$$ 
the corresponding structural functors; the point is that {\it $\widetilde{\pi^{_{\rm b}}}$ admits an essentially unique section\/} as  proves the theorem below. First of all, we need the following lemma.

\bigskip
\noindent
{\bf Lemma £4.4.} {\it For any subgroup $Q$ of $P$ there is a group homomorphism
$\mu_Q\,\colon \F(Q)\to \widetilde{\L^{^{\rm b}}} (Q)$ fulfilling $\mu_Q\circ \kappa_Q = \widetilde{\tau^{_{\rm b}}_Q}\,.$\/}

\medskip
\noindent
{\bf Proof:} Since we can choose an $\F\-$isomorphism $\theta\,\colon Q\cong Q'$ such that       $Q'$ is fully normalized in $\F$ and $\theta$ can be lifted to 
$\widetilde{\L^{^{\rm b}}} (Q',Q)\,,$ we may assume that $Q$ is fully normalized in $\F\,.$ 

\smallskip
We apply [8, Lemma 18.8] to the groups $\F (Q)$ and $\widetilde{\L^{^{\rm b}}} (Q)\,,$ to the normal $p\-$subgroup ${\rm Ker} (\widetilde{\pi^{_{\rm b}}_Q })$ of $\widetilde{\L^{^{\rm b}}} (Q)$ and to the
group homomorphism ${\rm id}_{\F (Q)}\,.$ We consider the group homomorphism
$\widetilde{\tau^{_{\rm b}}_Q}\,\colon N_P (Q)\to \widetilde{\L^{^{\rm b}}} (Q)$ and, for any subgroup $R$ of $N_P (Q)$ and any $\alpha\in \F (Q)$ such that 
$\alpha\circ \F_{\!R} (Q)\circ \alpha^{-1}\i \F_{\!P} (Q)\,,$ it~follows from 
[8, Proposition 2.11]
that there exists $\zeta\in \F \big(N_P (Q),Q\. R\big)$ extending $\alpha\,;$ then, it follows from [8, 17.11.2] that there exists $x\in  \widetilde{\L^{^{\rm b}}} (Q)$ fulfilling 
$$\widetilde{\tau^{_{\rm b}}_Q} \big(\zeta (v)\big) ={}^x \widetilde{\tau^{_{\rm b}}_Q} (v)
\eqno £4.4.1\phantom{.}$$
for any $v\in R\,.$ That is to say, condition 18.8.1 in  [8, Lemma 18.8]  is fulfilled and therefore this lemma proves the existence of $\mu_Q$ as announced.

\bigskip
\noindent
{\bf Theorem £4.5.} {\it  With the notation above, the structural functor 
$\widetilde{\pi^{_{\rm b}}}$ admits a unique natural $\F\-$isomorphism class of $\F\-$locality functorial sections.\/}

\medskip
\noindent
{\bf Proof:} We consider the filtration of $\widetilde{\L^{^{\rm b}} }$ induced by the filtration
of the basic $\F\-$locality introduced in section 3 and then argue by induction. That is to say,
recall that we denote by $\C_P$ a set of representatives for the set of $P\-$conjugacy classes of subgroups $U$ of $P$ (cf. £2.13); now, for any subset $\N$ of $\C_P$
fulfilling condition £3.2.1, we consider the obvious functor (cf. £3.2)
$$\hat\tau^{_{\rm b}} ({\tilde\frak c}^\frak f_\F)\. \tilde\frak k^{^\N}_\F :  \tilde\F\too \Ab
\eqno £4.5.1\phantom{.}$$
sending any subgroup $Q$ of $P$ to 
$\hat\tau^{_{\rm b}}_Q \big({\tilde\frak c}^\frak f_\F (Q)\big).\tilde\frak k^{^\N}_\F (Q)\,,$
and  the quotient $\F\-$lo-cality $\widetilde{\L^{^{{\rm b},\N}}} = 
\L^{^{\rm b}} \big/\big(\hat\tau^{_{\rm b}} ({\tilde\frak c}^\frak f_\F)\. \tilde\frak k^{^\N}_\F\circ \tilde\pi^{_{\rm b}} \big)$ with the structural functors
$$\widetilde{\tau^{_{\rm b,\N}}} : \T_P\too \widetilde{\L^{^{{\rm b},\N}}}
\qq \widetilde{\pi^{_{\rm b,\N}}}  :  \widetilde{\L^{^{{\rm b},\N}}} \too \F
\eqno £4.5.2.$$
Note that if $\N = \emptyset$ then $\widetilde{\L^{^{{\rm b},\N}}} 
= \widetilde{\L^{^{\rm b}} }\,;$ hence,  arguing by induction on $\vert\C_P - \N\vert\,,$  it suffices to prove that $\widetilde{\pi^{_{\rm b,\N}}}$ admits a unique natural $\F\-$isomorphism class of $\F\-$locality functorial sections.

\smallskip
Moreover,  if $\N = \C_P$ then $ \tilde\frak k^{^\N}_\F =  \tilde\frak k^{^{\rm b}}_\F\,;$ therefore
$\widetilde{\L^{^{{\rm b},\N}}} = \F$ and $\widetilde{\pi^{^{{\rm b},\N}}} = \id_\F \,,$
so that we may assume that $\N\not= \C_P\,;$ then, we fix a minimal element $U$ 
in~$\C_P - \N\,,$ setting $\M = \N\cup \{U\}$ and 
$\tilde\frak k^{^U}_\F = \tilde\frak k^{^\M}_\F/ \tilde\frak k^{^\N}_\F\,.$ If $U\not= P$ then 
$\M\not= \C_P$ and, as a matter of fact, we have
$\hat\tau^{_{\rm b}} ({\tilde\frak c}^\frak f_\F) \cap  \tilde\frak k^{^\M}_\F = \{0\}\,,$ so that
$$\big(\hat\tau^{_{\rm b}} ({\tilde\frak c}^\frak f_\F)\. \tilde\frak k^{^\M}_\F\big) \big/
\big(\hat\tau^{_{\rm b}} ({\tilde\frak c}^\frak f_\F)\. \tilde\frak k^{^\N}_\F\big)\cong 
\tilde\frak k^{^U}_\F
\eqno £4.5.3;$$
in this case, for any $m\in \Bbb N$ we simply denote by $\tilde\frak l^{^{U,m}}_\F$ the converse image of~$p^m\. \tilde\frak k^{^U}_\F$ 
in~$\hat\tau^{_{\rm b}} ({\tilde\frak c}^\frak f_\F)\. \tilde\frak k^{^\M}_\F\,;$  set 
$\widetilde{\L^{^{{\rm b},U,m}}} = \L^{^{\rm b}}/\tilde\frak l^{^{U,m}}_\F$ and, coherently,
denote by $\widetilde{\pi^{^{{\rm b},U,m}}}$ and $\widetilde{\tau^{^{{\rm b},U,m}}}$  the corresponding structural functors. Note that, by £3.8 and
£3.10.5 above we get 
$$\tilde\frak l^{^{U,m}}_\F/\tilde\frak l^{^{U,m +1}}_\F\cong \tilde\frak r^{U,\circ}_{\F,m}
\eqno £4.5.4\phantom{.}$$
and in particular, by Corollary £3.13, for any $n\in \Bbb N$ we still get
$$\Bbb H^{^n}_* (\tilde\F,\tilde\frak l^{^{U,m}}_\F/\tilde\frak l^{^{U,m +1}}_\F) = \{0\}
\eqno £4.5.5.$$

\smallskip
If $U = P$ then $\M = \C_P\,,$ so that in this case 
$\tilde\frak k^{^\M}_\F = \tilde\frak k^{^{\rm b}}_\F$ and, denoting by $\frak d_P\,\colon \tilde\F\to \Ab$ the functor mapping $P$ on $Z(P)$ and any other subgroup of $P$ on $\{0\}\,,$ from £3.7 and £3.8 it is easily checked that 
$$ \tilde\frak k^{^{\rm b}}_\F\big/\big(\hat\tau^{_{\rm b}} ({\tilde\frak c}^\frak f_\F)\. 
\tilde\frak k^{^{\C_P - \{P\}}}_\F\big)
\cong \prod_{\tilde\sigma\in \tilde\F (P)} \frak d_P\Big/ \Delta (\frak d_P)
\eqno £4.5.6\phantom{.}$$
where $\Delta$ denotes de usual {\it diagonal map\/}; but, similarly we have
$$\tilde\frak r^{P,\circ}_{\F,m}\cong \prod_{\tilde\sigma\in \tilde\F (P)} \frak s_m\circ\frak d_P
\eqno £4.5.7\phantom{.}$$
\eject
\noindent
and, according to Corollary £3.13, we get 
$\Bbb H^{^n}_{\F_{\!P}} (\tilde\F, \frak s_m\circ\frak d_P) =\{0\}\,;$ moreover, since $p$ does not divide $\vert\tilde\F (P)\vert\,,$ we still have
$$ \prod_{\tilde\sigma\in \tilde\F (P)} \frak d_P\Big/ \Delta (\frak d_P) \cong 
 \prod_{\tilde\sigma\in \tilde\F (P) -\{{\rm \widetilde{id}}_P\}} \frak d_P
 \eqno £4.5.8;$$
 hence, still setting $\tilde\frak l^{^{P,m}}_\F = p^m\. \tilde\frak k^{^{\rm b}}_\F$ and 
 $\widetilde{\L^{^{{\rm b},P,m}}} = \L^{^{\rm b}}/\tilde\frak l^{^{P,m}}_\F\,,$ we still get
 $$\Bbb H^{^n}_* (\tilde\F,\tilde\frak l^{^{P,m}}_\F/\tilde\frak l^{^{P,m +1}}_\F) = \{0\}
\eqno £4.5.9.$$
 
\smallskip
Further, we denote by $\C_\F$ a set of representatives,  fully normalized in~$\F\,,$ for the 
$\F\-$isomorphism classes of subgroups of $P$ and, for any subgroup $Q$ in~$\C_\F\,,$  we choose a  group homomorphism $\mu_Q\,\colon \F (Q)\to \widetilde{\L^{^{\rm b}}} (Q)$ as in Lemma £4.4 above and, for any $m\in \Bbb N\,,$ simply denote by $\mu_Q^{_m}$ the corresponding group homomorphism from $\F (Q)$ to $\widetilde{\L^{^{{\rm b},U,m}}} (Q)\,.$ For any $\F\-$morphism 
$\varphi\,\colon R\to Q$ denote by $\F (Q)_\varphi$ and by 
$\widetilde{\L^{^{{\rm b},U,m}}} (Q)_\varphi$ the respective stabilizers of $\varphi (R)$ in 
$\F (Q)$ and in $\widetilde{\L^{^{{\rm b},U,m}}} (Q)\,;$  it is clear that we have a group homomorphism $a_\varphi\,\colon \F (Q)_\varphi \to \F (R)$ fulfilling 
$ \eta\circ \varphi = \varphi\circ a_\varphi (\eta)$
for any $\eta\in \F (Q)_\varphi\,;$ similarly, for any 
$x^{_m}\in  \widetilde{\L^{^{{\rm b},U,m}}} (Q,R)$ we have a group homomorphism 
$$a_{x^{_m}} : \widetilde{\L^{^{{\rm b},U,m}}} (Q)_\varphi \too
\widetilde{\L^{^{{\rm b},U,m}}} (R)
\eqno £4.5.10\phantom{.}$$
fulfilling $y^{_{m }}\!\. x^{_m} = x^{_m} \. a_{ x^{_m} } (y^{_{m}})$
for any $y^{_{m}}\!\in  \widetilde{\L^{^{{\rm b},U,m}}} (Q)_\varphi\,.$

\smallskip
For any  subgroups $Q$ and $R$ in $\C_\F$ and any $\F\-$morphism 
$\varphi\,\colon R\to Q\,,$ $\F_{\!P} (Q)$ and $\F_{\!P} (R)$ are respective 
Sylow $p\-$subgroups of $\F (Q)$ and $\F (R)$\break  [8, Proposition 2.11]; therefore,  there are $\alpha\in \F (Q)$ such that  $\F_{\!P} (Q)^\alpha$ contains a Sylow $\F_{\!P} (Q)^\alpha_\varphi$ $p\-$subgroup of~$\F (Q)_\varphi$ and  $\beta\in \F (R)$ such that  
$a_\varphi \big(\F_{\!P} (Q)^\alpha_\varphi\big)$ is contaioned in $\F_{\!P} (R)^\beta\,.$ Thus, we choose a set of representati-ves $\F_{\!Q,R}$ for the set of double classes 
$\F (Q)\backslash \F (Q,R)/\F (R)$ such that, for any $\varphi$ in $\F_{\!Q,R}\,,$ 
$\F_{\!P} (Q)$ contains a Sylow $p\-$subgroup of $\F (Q)_\varphi$ and 
$a_\varphi \big(\F_{\!P} (Q)_\varphi\big)$ is contaioned in $\F_{\!P} (R)\,;$ of course, we choose $\F_{Q,Q} = \{{\rm id}_Q\}\,.$

\smallskip
With all this notation and arguing by induction on $\vert\C_P - \N\vert$ and on~$m\,,$   we will prove that there is a functorial section 
$$\sigma^{_{m}} : \F\too \widetilde{\L^{^{{\rm b},U,m }}}
\eqno £4.5.11\phantom{.}$$
 such that, for any $Q\in C_\F$ and any $u\in Q\,,$ we have $\sigma^{_m} \big(\kappa_Q (u)\big) = \widetilde{\tau^{^{{\rm b},U,m }}_Q} (u)\,,$ and that, for any   groups $Q$ and $R$ in~$\C_\F$, and any $\F\-$morphism $\varphi\,\colon Q\to R$ in~$\F_{Q,R}\,,$ we have the commutative diagram
$$\matrix{\F (Q)_\varphi & \buildrel \mu_{_Q}^{_m}\over \too 
&\widetilde{\L^{^{{\rm b},U,m }}} (Q)_\varphi\cr
\hskip-20pt{\scriptstyle a_\varphi}\big\downarrow&
& \big\downarrow {\scriptstyle a_{\sigma^{_m} (\varphi)}}\hskip-20pt\cr
\F (R) & \buildrel \mu_{_R}^{_m}\over \too 
&\widetilde{\L^{^{{\rm b},U,m }}} (R)\cr}
\eqno £4.5.12.$$ Since we have $ \widetilde{\pi^{^{{\rm b},U,0}}} =  \widetilde{\pi^{^{{\rm b},\M}}}$ 
and $\vert\M\vert = \vert\N\vert + 1\,,$ by the induction hypothesis we actually may assume that $m\not= 0\,,$ that $ \widetilde{\pi^{^{{\rm b},U,m -  1}}}$ admits a functorial section 
$$\sigma^{_{m - 1}} : \F\too \widetilde{\L^{^{{\rm b},U,m - 1}}}
\eqno £4.5.13\phantom{.}$$
which fulfills the conditions above.

\smallskip
Then, for any $\varphi\in \F_{\!Q,R}$ it follows from [8, Proposition 2.11], applied to the  inverse $\varphi^*$ of the isomorphism $\varphi_*\,\colon R\cong \varphi (R)$ induced by 
$\varphi\,,$ that there exists an $\F\-$morphism $\zeta\,\colon N_P (Q)_\varphi\to N_P (R)$ fulfilling $\zeta\big(\varphi (v)\big) = v$ for any $v\in R\,,$ so that we easily get the following commutative diagram`
$$\matrix{N_P (Q)_\varphi &\buildrel \kappa_Q \over 
\too & \F (Q)_\varphi\cr
\hskip-15pt {\scriptstyle \zeta} \big\downarrow &
&\big\downarrow{\scriptstyle a_\varphi }\hskip-10pt  \cr
N_P (R) &\buildrel \kappa_R \over 
\too&\F (R)&\cr}
\eqno £4.5.14;$$
note that, if $Q = R$ and $\varphi = \kappa_Q (u)$ for some $u\in Q\,,$ we may assume that\break
$\zeta = \kappa_{N_P (Q)} (u)\,.$ In particular, since $\sigma^{_{m - 1}}$ fulfills the corresponding commutative diagram £4.5.12,  we still get the following commutative diagram 
$$\matrix{N_P (Q)_\varphi &\buildrel \widetilde{\tau^{^{{\rm b},U,m -1}}_Q} \over 
{\hbox to 40pt{\rightarrowfill}} & \widetilde{\L^{^{{\rm b},U,m - 1}}} (Q)_\varphi\cr
\hskip-10pt {\scriptstyle \zeta} \big\downarrow &
&\big\downarrow{\scriptstyle a_{\sigma^{_{m -1}} (\varphi)} }\hskip-20pt  \cr
N_P (R) &\buildrel \widetilde{\tau^{^{{\rm b},U,m -1}}_R} \over 
{\hbox to 40pt{\rightarrowfill}} & \widetilde{\L^{^{{\rm b},U,m - 1}}} (R)&\cr}
\eqno £4.5.15$$

\smallskip
The first step is, for any $\F\-$morphism $\varphi$ in $\F_{Q,R}\,,$ to choose a suitable lifting 
$\widehat{\sigma^{_{m -1}} (\varphi) }$ of $\sigma^{_{m - 1}}(\varphi)$ in  
$\widetilde{\L^{^{{\rm b},U,m }}} (Q,R)\,.$ We start by choosing a lifting 
$\widehat{\sigma^{_{m -1}} (\zeta) }$ of $\sigma^{_{m -1}} (\zeta)$ in the obvious stabilizer
$\widetilde{\L^{^{{\rm b},U,m }}} \big(N_P (R),N_P (Q)_\varphi\big)_{R,\varphi (R)}\,;$
thus, by the {\it coherence\/} of $\widetilde{\L^{^{{\rm b},U,m }}}$ (cf. (Q)), for any 
$u\in N_P (Q)_\varphi$ we have
$$\widehat{\sigma^{_{m -1}} (\zeta) }\.\widetilde{\tau^{^{{\rm b},U,m }}_{N_P (Q)_\varphi}} (u) = \widetilde{\tau^{^{{\rm b},U,m }}_{N_P (R)}}\big(\zeta (u)\big)\.\widehat{\sigma^{_{m -1}} (\zeta) }
\eqno £4.5.16;$$
moreover, by the {\it divisibility\/} of $\widetilde{\L^{^{{\rm b},U,m }}}$ (cf. £2.4), we find
$z_\varphi\in  \widetilde{\L^{^{{\rm b},U,m }}} \big(R,\varphi(R)\big)$ fulfilling
$$\widehat{\sigma^{_{m -1}} (\zeta) }\.\widetilde{\tau^{^{{\rm b},U,m }}_{N_P (Q)_\varphi,\varphi(R)}} (1) = \widetilde{\tau^{^{{\rm b},U,m }}_{N_P (R),R}} (1)\. z_\varphi
\eqno £4.5.17;$$
similarly, $\sigma^{_{m -1}} (\zeta)$ resrtricts to $\sigma^{_{m -1}} (\varphi^*)\in \widetilde{\L^{^{{\rm b},U,m -1}}} \big(R,\varphi(R)\big)\,,$ so that it is easily checked that $z_\varphi$ lifts $\sigma^{_{m -1}} (\varphi^*)$ to $ \widetilde{\L^{^{{\rm b},U,m }}} \big(R,\varphi(R)\big)$
and therefore
 $\widetilde{\sigma^{_{m -1}} (\varphi) } = \widetilde{\tau^{^{{\rm b},U,m }}_{Q\. \varphi (R)}} (1)\; z_\varphi^{-1}$ lifts $\sigma^{_{m -1}} (\varphi)$ to $ \widetilde{\L^{^{{\rm b},U,m }}} (Q,R)\,.$
 \eject

\smallskip
Then, from £4.5.16 and £4.5.17 above,  for any 
$u\in N_P (Q)_\varphi$ we get
$$\eqalign{\widehat{\sigma^{_{m -1}} (\zeta) }\.
\widetilde{\tau^{^{{\rm b},U,m }}_{N_P (Q)_\varphi}} (u)&\.
.\widetilde{\tau^{^{{\rm b},U,m }}_{N_P (Q)_\varphi,\varphi(R)}} (1) = 
\widetilde{\tau^{^{{\rm b},U,m }}_{N_P (R),R}} (1)\. 
z_\varphi \. \widetilde{\tau^{^{{\rm b},U,m }}_{\varphi (R)}} (u)\cr
\Vert\cr
\widetilde{\tau^{^{{\rm b},U,m }}_{N_P (R)}}\big(\zeta (u)\big)\.\widehat{\sigma^{_{m -1}} (\zeta) }&\. .\widetilde{\tau^{^{{\rm b},U,m }}_{N_P (Q)_\varphi,\varphi(R)}} (1)\cr
&= \widetilde{\tau^{^{{\rm b},U,m }}_{N_P (R),R}} (1)\.
\widetilde{\tau^{^{{\rm b},U,m }}_R}\big(\zeta (u)\big)\. z_\varphi\cr}
\eqno £4.5.18\phantom{.}$$
and therefore we still get $z_\varphi\. \widetilde{\tau^{^{{\rm b},U,m }}_{\varphi (R)}} (u) = \widetilde{\tau^{^{{\rm b},U,m }}_R}\big(\zeta (u)\big)\. z_\varphi\,,$ so that
$$ \widetilde{\tau^{^{{\rm b},U,m }}_Q} (u)\. \widetilde{\sigma^{_{m -1}} (\varphi) } = \widetilde{\sigma^{_{m -1}} (\varphi) }\. \widetilde{\tau^{^{{\rm b},U,m }}_R}\big(\zeta (u)\big)
\eqno £4.5.19\phantom{.}$$
or, equivalently, we have $a_{ \widetilde{\sigma^{_{m -1}} (\varphi) } } \big( \widetilde{\tau^{^{{\rm b},U,m }}_Q} (u)\big) = \widetilde{\tau^{^{{\rm b},U,m }}_R}\big(\zeta (u)\big)\,.$

\smallskip
At this point, we will apply the uniqueness part of [8, Lemma 18.8] to the groups 
$\F (Q)_\varphi$ and $ \widetilde{\L^{^{{\rm b},U,m }}} (R)$  and to the composition of group homomorphisms
$$a_{\sigma^{^{m -1}}(\varphi)}\circ \mu_Q^{_{m -1}} : \F (Q)_\varphi\too  \widetilde{\L^{^{{\rm b},U,m -1}}} (Q)_\varphi\too \widetilde{\L^{^{{\rm b},U,m -1}}} (R)
\eqno £4.5.20,$$
together with the composition of group homomorphisms
$$\widetilde{\tau^{^{{\rm b},U,m }}_R}\circ \zeta : N_P (Q)_\varphi\too N_P (R)\too 
 \widetilde{\L^{^{{\rm b},U,m}}} (R)
 \eqno £4.5.21.$$
 Now, according to the commutative diagrams £4.5.12 for $m -1$ and £4.5.14, and to equality £4.5.18 above, the two group homomorphisms
 $$\eqalign{a_{\widetilde{\sigma^{^{m -1}}(\varphi)}}\circ \mu_Q^{_{m}} &: \F (Q)_\varphi\too  \widetilde{\L^{^{{\rm b},U,m}}} (Q)_\varphi\too \widetilde{\L^{^{{\rm b},U,m}}} (R)\cr
\mu_R^{_{m}}\circ  a_\varphi &: \F (Q)_\varphi\too  \F (R)\too \widetilde{\L^{^{{\rm b},U,m}}} (R)\cr}
 \eqno £4.5.22$$
 both fulfill the conclusion of [8, Lemma 18.8]; consequently, according to this lemma, there is
 $k_\varphi$ in the kernel of the canonical homomorphism from
  $ \widetilde{\L^{^{{\rm b},U,m }}} (R)$
 to $ \widetilde{\L^{^{{\rm b},U,m -1}}} (R)$ such that, denoting by 
 ${\rm int}_{\widetilde{\L^{^{{\rm b},U,m}}} (R)} (k_\varphi)$ the conjugation by $k_\varphi$ in
  $\widetilde{\L^{^{{\rm b},U,m}}} (R)\,,$ we have
  $${\rm int}_{\widetilde{\L^{^{{\rm b},U,m}}} (R)} (k_\varphi)\circ a_{\widetilde{\sigma^{^{m -1}}(\varphi)}}\circ \mu_Q^{_{m}} = \mu_R^{_{m}}\circ  a_\varphi
  \eqno £4.5.23;$$
  but, it is easily checked that 
  $${\rm int}_{\widetilde{\L^{^{{\rm b},U,m}}} (R)} (k_\varphi)\circ a_{\widetilde{\sigma^{^{m -1}}(\varphi)}} = a_{\widetilde{\sigma^{^{m -1}} (\varphi)}\. k_\varphi^{-1}}
 \eqno £4.5.24.$$
 \eject
 \noindent

 \smallskip
 Finally, we choose $\widehat{\sigma^{^{m -1}} (\varphi)} 
 = \widetilde{\sigma^{^{m -1}} (\varphi)}\. k_\varphi^{-1}\,,$  lifting indeed $\sigma^{_{m -1}} (\varphi)$
 to $\widetilde{\L^{^{{\rm b},U,m}}} (Q,R)$ and, according to equalities £4.5.23 and £4.5.24, fulfilling the following commutative 
 diagram
$$\matrix{\F (Q)_\varphi & \buildrel \mu_{_Q}^{_m}\over \too 
&\widetilde{\L^{^{{\rm b},U,m }}} (Q)_\varphi\cr
\hskip-20pt{\scriptstyle a_\varphi}\big\downarrow&
& \big\downarrow {\scriptstyle a_{\widehat{\sigma^{_{m -1}} (\varphi)}}}\hskip-30pt\cr
\F (R) & \buildrel \mu_{_R}^{_m}\over \too 
&\widetilde{\L^{^{{\rm b},U,m }}} (R)\cr}
\eqno £4.5.25;$$
note that, if $Q = R$ and $\varphi = \kappa_Q (u)$ for some $u\in Q\,,$ this choice is com-patible  wtih 
$\widehat{\sigma^{^{m -1}} \big( \kappa_Q (u)\big)} = \widetilde{\tau^{^{{\rm b},U,m }}_Q} (u)\,.$ 
In particular, considering the action of $\F (Q)\times \F (R)\,,$ by composition on the left- and
on the right-hand, on $\F (Q,R)$ and on $\widetilde{\L^{^{{\rm b},U,m }}} (Q,R)$ {\it via\/} $\mu^{_m}_Q$ and $\mu^{_m}_R\,,$ we have the inclusion of stabilizers
$$\big(\F (Q)\times \F (R)\big)_\varphi \i \big(\F (Q)\times \F (R)\big)_{\sigma^{_{m -1}} (\varphi)}
\eqno £4.5.26;$$
indeed, it is quite clear that $(\alpha,\beta)\in \big(\F (Q)\times \F (R)\big)_\varphi $ forces
$\alpha\in \F (Q)_\varphi\,;$  then, since $\alpha\circ \varphi = \varphi\circ a_\varphi (\alpha)\,,$
we get $\beta =  a_\varphi (\alpha)$ and the inclusion above follows from the commutativity of diagram £4.5.26.

\smallskip
This allows us to choose a family of liftings $\{\widehat{\sigma^{_{m -1}} (\varphi)}\}_\varphi\,,$ where $\varphi$ runs over the set of  $\F\-$morphisms, which is compatible with 
$\F\-$isomorphisms; precisely,  for any  pair of   subgroups $Q$ and $R$ in~$\C_\F\,,$ and any $\varphi\in \F_{Q,R}\,,$ we choose a lifting  $\widehat{\sigma^{_{m -1}} (\varphi)}$ of 
$\sigma^{_{m -1}} (\varphi)$ in  $\widetilde{\L^{^{{\rm b},U,m }}} (Q,R)$ as above.
 Then, any subgroup $Q$ of $P$ determines a unique  $\hat Q$ in $\C_\F$ which is 
$\F\-$isomorphic to~$Q$ and we choose an $\F\-$isomorphism 
$\omega_Q\,\colon \hat Q\cong Q$ and a lifting 
$x_Q\in \widetilde{\L^{^{{\rm b},U,m }}} (Q,\hat Q)$ of` $\omega_Q\,;$  in particular, we choose
$\omega_{\hat Q} = {\rm id}_{\hat Q}$ and 
$x_{\hat Q} =  \widetilde{\tau^{^{{\rm b},U,m }}_{\hat Q}} (1)\,.$ Thus, any $\F\-$morphism 
$\varphi\,\colon R\to Q$ determines subgroups $\hat Q$ and $\hat R$ in~$\C_\F$ and an element $\hat\varphi$ in  $\F_{\hat Q,\hat R}$ in such a way that there are 
$\alpha_\varphi\in \F (\hat Q)$ and  $\beta_\varphi\in \F (\hat R)$ fulfilling
 $$\varphi = \omega_Q\circ \alpha_\varphi\circ \hat\varphi\circ\beta^{-1}_\varphi	\circ \omega_R^{-1π}
 \eqno £4.5.27\phantom{.}$$
 and we define
 $$\widehat{\sigma^{_{m -1}} (\varphi)} = x_Q\. \mu_{\hat Q}^{_m} (\alpha_\varphi)\.
 \widehat{\sigma^{_{m -1}} (\hat\varphi)} \.  \mu_{\hat R}^{_m} (\beta_\varphi)^{-1} \. x_R^{-1}
 \eqno £4.5.28;$$
once again,, if $Q = R$ and $\varphi = \kappa_Q (u)$ for some $u\in Q\,,$ we actually get
$\widehat{\sigma^{^{m -1}} \big( \kappa_Q (u)\big)} = \widetilde{\tau^{^{{\rm b},U,m }}_Q} (u)\,.$
 This definition does not depend on the choice of  $(\alpha_\varphi,\beta_\varphi)$ since for another choice $(\alpha',\beta')$ we clearly have $\alpha' = \alpha_\varphi\circ \alpha''$ and 
 $\beta' = \beta_\varphi\circ \beta''$ for a suitable $(\alpha'',\beta'')$  in~$\big(\F (\hat Q)\times \F (\hat R)\big)_{\hat\varphi}$ and it suffices to apply inclusion £4.5.26.
 \eject

 \smallskip 
 Moreover, for any pair of  $\F\-$isomorphisms $\zeta\,\colon Q\cong Q'$ and 
 $\xi\,\colon R\cong R'\,,$  considering $\varphi' = \zeta\circ \varphi\circ \xi^{-1}$ we claim that  
 $$\widehat{\sigma^{_{m -1}} (\varphi')} = \widehat{\sigma^{_{m -1}} (\zeta)} \. \widehat{\sigma^{_{m -1}} (\varphi) } \. \widehat{\sigma^{^{m -1}} (\xi)}^{-1}
\eqno £4.5.29;$$
indeed, it is clear that $Q'$ also determines $\hat Q$ in  $\C_\F$ and therefore, if we have
$\zeta = \omega_Q\circ \alpha_\zeta\circ \omega_{Q'}^{-1}$ then we obtain 
$\widehat{\sigma^{_{m -1}} (\zeta)} = x_{Q'}\. \mu_{\hat Q}^{_m} (\alpha_\zeta)\. x_Q^{-1}\,;$
similarly,   if we have $\xi = \omega_R\circ \beta_\xi\circ \omega_{R'}^{-1}$  we also obtain
$\widehat{\sigma^{_{m -1}} (\xi)}^{-1} = x_{R}\. \mu_{\hat R}^{_m} (\beta_\xi)^{-1}\. x_{R'}^{-1}\,;$
further, $\varphi'$ also determines $\hat\varphi$ in $\F_{\hat Q,\hat R}\,;$
consequently, we get
$$\eqalign{\widehat{\sigma^{_{m -1}} (\zeta)} &\. \widehat{\sigma^{_{m -1}} (\varphi) } \. \widehat{\sigma^{^{m -1}} (\xi)}^{-1} \cr
&= (x_{Q'}\. \mu_{\hat Q}^{_m} (\alpha_\zeta)\. x_Q^{-1})\.   \widehat{\sigma^{_{m -1}} (\varphi)}\. (x_{R}\. \mu_{\hat R}^{_m} (\beta_\xi)^{-1}\. x_{R'}^{-1})\cr
&= x_{Q'} \.  \mu_{\hat Q}^{_m} (\alpha_\zeta \circ \alpha_\varphi)\. \hat\varphi\. 
\mu_{\hat R}^{_m} (\beta_\varphi^{-1}\circ \beta_\xi^{-1})\. x_{R'}^{-1} = \widehat{\sigma^{_{m -1}} (\varphi')}\cr}
\eqno £4.5.30.$$

Recall that we have the exact sequence of {\it contravariant\/} functors from $\tilde\F$ to $\Ab$ (cf.~£2.7 and £2.8)
 $$0\too \tilde\frak l^{^{U,m - 1}}_\F/\tilde\frak l^{^{U,m}}_\F\too
 \widetilde{\Ker} (\widetilde{\bar\pi^{^{{\rm b},U,m}}})\too \widetilde{\Ker} (\widetilde{\bar\pi^{^{{\rm b},U,m - 1}}})
  \too 0
 \eqno £4.5.31;$$
 hence, for another $\F\-$morphism $\psi: T\to R$ we clearly have
$$\widehat{\sigma^{_{m -1}} (\varphi)}\. \widehat{\sigma^{_{m -1}} (\psi)}
 = \widehat{\sigma^{_{m -1}} (\varphi\circ\psi)}\. \gamma_{\psi,\varphi}^{_m}
 \eqno £4.5.32\phantom{.}$$
 for some $\gamma_{\psi,\varphi}^{_m}$ in 
  $(\tilde\frak l^{^{U,m - 1}}_\F/\tilde\frak l^{^{U,m}}_\F) (T)\,.$
 That is to say, borrowing notation and terminology from  [8,~A2.8], we get a correspondence sending any {\it $\F\-$chain\/}  $\frak q\,\colon \Delta_2\to \F$  to the element 
 $ \gamma^{_m}_{\frak q (0\bullet 1),\frak q(1\bullet 2)}$ 
 in $(\tilde\frak l^{^{U,m - 1}}_\F/\tilde\frak l^{^{U,m}}_\F) \big(\frak q(0)\big)$ and, setting
$$\Bbb C^n \big(\tilde\F, \tilde\frak l^{^{U,m - 1}}_\F/\tilde\frak l^{^{U,m}}_\F\big) = 
\prod_{\tilde\frak q\in \Fct(\Delta_n,\tilde\F)} (\tilde\frak l^{^{U,m - 1}}_\F/\tilde\frak l^{^{U,m}}_\F) \big(\tilde\frak q(0)\big)
\eqno £4.5.33\phantom{.}$$
 for any $n\in \Bbb N\,,$  we claim that this correspondence determines an 
 {\it  stable\/} element $\gamma^{_m}$ of 
 $ \Bbb C^2 \big(\tilde\F,\tilde\frak l^{^{U,m - 1}}_\F/\tilde\frak l^{^{U,m}}_\F\big) $ 
 [8,~A3.17].

\smallskip
 Indeed, for another $\F-$isomorphic {\it $\F\-$chain\/} $\frak q'\,\colon \Delta_2\to \F$ and a {\it natural $\F\-$isomorphism\/} $\nu\,\colon \frak q\cong\frak q'\,,$
setting  
$$\eqalign{T = \frak q (0)  ,\ T' = \frak q'(0) ,\ R = \frak q(1)
& ,\ R'  = \frak q' (1) ,\ Q = \frak q (2) ,\ Q' = \frak q' (2)\cr
\psi = \frak q (0\bullet 1)\!\!\quad,\quad \!\!\varphi = \frak q(1\bullet 2)\!\!\quad &,\quad  \!\!\psi' = \frak q'(0\bullet 1)\!\!\quad,\quad \!\!\varphi' = \frak q'(1\bullet 2)\cr
\nu_0 = \eta\quad,\quad  
\nu_1 = &\ \xi  \qq 
\nu_2 = \zeta\cr}
\eqno £4.5.34,$$
 from~£4.5.30 we have
 $$\eqalign{\widehat{\sigma^{_{m -1}} (\varphi')} &=  \widehat{\sigma^{_{m -1}} (\zeta)}\. \widehat{\sigma^{_{m -1}}(\varphi)} \.  \widehat{\sigma^{_{m -1}} (\xi)}^{-1}\cr
 \widehat{\sigma^{_{m -1}} (\psi')} &=  \widehat{\sigma^{_{m -1}} (\xi)}\. 
 \widehat{\sigma^{_{m -1}} (\psi)} \.  \widehat{\sigma^{_{m -1}} (\eta)}^{-1} \cr 
\widehat{\sigma^{_{m -1}} (\varphi'\circ\psi')} 
&=  \widehat{\sigma^{_{m -1}} (\zeta)} \. \widehat{\sigma^{_{m -1}} (\varphi\circ\psi)} \. 
 \widehat{\sigma^{_{m -1}} (\eta)}^{-1}\cr}
 \eqno £4.5.35\phantom{.}$$
 \eject
\noindent
 and therefore we get
$$\eqalign{&\widehat{\sigma^{_{m -1}} (\varphi' \circ\psi')} \.  \gamma^{_m}_{\varphi',\psi'} 
= \widehat{\sigma^{_{m -1}} (\varphi')} \.  \widehat{\sigma^{_{m -1}} (\psi')} \cr
& = \big(\widehat{\sigma^{_{m -1}} (\zeta)} \. \widehat{\sigma^{_{m -1}} (\varphi)} \. 
\widehat{\sigma^{_{m -1}} (\xi)^{-1}}\big)\. 
\big(\widehat{\sigma^{_{m -1}} (\xi)} \. \widehat{\sigma^{_{m -1}} (\psi)} \. 
 \widehat{\sigma^{_{m -1}} (\eta)}^{-1} \big)\cr
& =  \widehat{\sigma^{_{m -1}} (\zeta)} \. \big(\widehat{\sigma^{_{m -1}} (\varphi \circ\psi) }\.  \gamma^{_m}_{\varphi,\psi} \big) \.  \widehat{\sigma^{_{m -1}} (\eta)}^{-1}\cr
&=  \widehat{\sigma^{_{m -1}} (\varphi'\circ\psi')} \. 
\big((\tilde\frak l^{^{U,m - 1}}_\F/\tilde\frak l^{^{U,m}}_\F)
( \widehat{\sigma^{_{m -1}} (\eta)}^{-1})\big)( \gamma^{_m}_{\varphi,\psi} ) \cr}
\eqno £4.5.36,$$ 
so that, by the divisibility of $ \widetilde{\L^{^{{\rm b},U,m }}}\,,$ we have
$$\gamma^{_m}_{\varphi',\psi'}  = \big((\tilde\frak l^{^{U,m - 1}}_\F/\tilde\frak l^{^{U,m}}_\F)
( \widehat{\sigma^{_{m -1}} (\eta)}^{-1})\big)( \gamma^{_m}_{\varphi,\psi})
\eqno £4.5.37;$$ 
this proves  that the correspondence $\gamma^{_m}$ sending $(\varphi,\psi)$ 
to~$\gamma^{m}_{\varphi,\psi}$ is {\it stable\/} and, in particular, that 
$\gamma^{m}_{\varphi,\psi}$ only depends on the corresponding $\tilde\F\-$morphisms 
$\tilde\varphi$ and~$\tilde\psi\,;$ thus we set $\gamma^{_m}_{\tilde\varphi,\tilde\psi}= \gamma^{_m}_{\varphi,\psi}\,.$

\smallskip
On the other hand, considering the usual differential map
$$d_2 : \Bbb C^2 \big(\tilde\F,\tilde\frak l^{^{U,m - 1}}_\F/\tilde\frak l^{^{U,m}}_\F\big)\too 
\Bbb C^3 \big(\tilde\F,\tilde\frak l^{^{U,m - 1}}_\F/\tilde\frak l^{^{U,m}}_\F\big)
\eqno £4.5.38,$$
we claim that $d_2 (\gamma^{_m}) = 0\,;$ indeed, for a third $\F\-$morphism $\varepsilon\,\colon W\to T$ we get
$$\eqalign{\big(\widehat{\sigma^{_{m -1}} (\varphi)} &\.\widehat{\sigma^{_{m -1}} (\psi)}  \big)\. \widehat{\sigma^{^{m -1}} (\varepsilon)} 
= (\widehat{\sigma^{_{m -1}} (\varphi\circ\psi)}\. \gamma^{_m}_{\tilde\varphi,\tilde\psi})\. \widehat{\sigma^{_{m -1}} (\varepsilon)} \cr
&= \big(\widehat{\sigma^{_{m -1}} (\varphi\circ\psi)} \. \widehat{\sigma^{_{m -1}} (\varepsilon)}\big)\.\big((\tilde\frak l^{^{U,m - 1}}_\F/\tilde\frak l^{^{U,m}}_\F)(\tilde\varepsilon)\big)(\gamma^{_m}_{\tilde\varphi,\tilde\psi})\cr
&= \widehat{\sigma^{_{m -1}} (\varphi\circ\psi\circ\varepsilon)} \.\gamma^{_m}_{\tilde\varphi\circ\tilde\psi,\tilde\varepsilon}\.\big((\tilde\frak l^{^{U,m - 1}}_\F/\tilde\frak l^{^{U,m}}_\F)(\tilde\varepsilon)\big) (\gamma^{_m}_{\tilde\varphi,\tilde\psi})\cr
\widehat{\sigma^{_{m -1}} (\varphi)} &\. \big(\widehat{\sigma^{_{m -1}} (\psi)} \. \widehat{\sigma^{_{m -1}} (\varepsilon)} \big)
 = \widehat{\sigma^{_{m -1}} (\varphi)} \.\big(\widehat{\sigma^{_{m -1}} (\psi\circ\varepsilon)} \. 
 \gamma^{_m}_{\tilde\psi,\tilde\varepsilon}\big)\cr 
& = \widehat{\sigma^{_{m -1}} (\varphi\circ\psi\circ\varepsilon)} \. \gamma^{_m}_{\tilde\varphi,\tilde\psi\circ\tilde\varepsilon}\. \gamma^{_m}_{\tilde\psi,\tilde\varepsilon}\cr}
\eqno £4.5.39\phantom{.}$$
and   the {\it divisibility\/} of $\widetilde{\L^{^{{\rm b},U,m}} }$ forces
$$\gamma^{_m}_{\tilde\varphi\circ\tilde\psi,\tilde\varepsilon}
\.\big((\tilde\frak l^{^{U,m - 1}}_\F/\tilde\frak l^{^{U,m}}_\F)(\tilde\varepsilon)\big)
(\gamma^{_m}_{\tilde\varphi,\tilde\psi}) = \gamma^{_m}_{\tilde\varphi,\tilde\psi\circ\tilde\varepsilon} \.\gamma^{_m}_{\tilde\psi,\tilde\varepsilon}
\eqno £4.5.40;$$
since ${\rm Ker}(\widetilde{\pi_W^{^{{\rm b},U,m  }}})$ is abelian, with the additive notation we obtain
$$0 = \big((\tilde\frak l^{^{U,m - 1}}_\F/\tilde\frak l^{^{U,m}}_\F) (\tilde\varepsilon)\big)
(\gamma^{_m}_{\tilde\varphi,\tilde\psi}) - \gamma^{_m}_{\tilde\varphi,\tilde\psi\circ\tilde\varepsilon} + \gamma^{_m}_{\tilde\varphi\circ\tilde\psi,\tilde\varepsilon} - 
\gamma^{_m}_{\tilde\psi,\tilde\varepsilon}
\eqno £4.5.41,$$
proving our claim.

\smallskip
At this point, it follows from equalities £4.5.5 and £4.5.9  that  
$\gamma^{_m} = d_1 (\beta^{^m})$ for some {\it stable\/} element 
$\beta^{^m} = (\beta^{^m}_{\tilde\frak r})_{\tilde\frak r\in \Fct(\Delta_1,\tilde\F^{^{\frak X}})}$ in 
$\Bbb C^1 \big(\tilde\F,\tilde\frak l^{^{U,m - 1}}_\F/\tilde\frak l^{^{U,m}}_\F\big)\,;$ 
that is to say, with the notation above we get
$$\gamma^{_m}_{\tilde\varphi,\tilde\psi} = 
\big((\tilde\frak l^{^{U,m - 1}}_\F/\tilde\frak l^{^{U,m}}_\F)(\tilde\psi)\big)
(\beta^{^m}_{\tilde\varphi})\.(\beta^{^m}_{\tilde\varphi\circ\tilde\psi})^{-1}\.\beta^{^m}_{\tilde\psi}
\eqno £4.5.42\phantom{.}$$
\eject
\noindent
where we  identify any $\tilde\F\-$morphism with the obvious 
{\it $\tilde\F\-$chain\/} $\Delta_1\to \tilde\F\,;$ hence, from equality~£4.5.32 we obtain
$$\eqalign{\big(\widehat{\sigma^{_{m -1}} (\varphi)} &\. (\beta^{^m}_{\tilde\varphi})^{-1}\big)
\.\big(\widehat{\sigma^{_{m -1}} (\psi )} \.( \beta^{^m}_{\tilde\psi})^{-1}\big)\cr
&=   \big( (\widehat{\sigma^{_{m -1}} (\varphi)} \.\widehat{\sigma^{_{m -1}} (\psi)} \big) \. \Big(\beta^{^m}_{\tilde\psi} \. \big((\tilde\frak l^{^{U,m - 1}}_\F/\tilde\frak l^{^{U,m}}_\F)(\tilde\psi)\big)(\beta^{^m}_{\tilde\varphi}) \Big)^{-1}\cr
&= \widehat{\sigma^{_{m -1}} (\varphi\circ\psi)} \.(\beta^{^m}_{\tilde\varphi\circ\tilde\psi})^{-1}\cr}
\eqno £4.5.43,$$
which amounts to saying that the correspondence $\sigma^{_m}$  sending $\varphi\in \F (Q,R)$ 
to~$\widehat{\sigma^{_{m -1}} (\varphi)} \. (\beta^{^m}_{\tilde\varphi})^{-1}\in \widetilde{\L^{^{{\rm b}, U,m}}}(Q,R)$ defines  a {\it functorial\/}  section of 
$\widetilde{\pi^{^{{\rm b},U,m}}}\,;$ note that, if $Q = R$ and $\varphi = \kappa_Q (u)$ for some $u\in Q\,,$ we have $\tilde\varphi = \widetilde{{\rm id}}_Q$ and 
$\beta^{^m}_{\tilde\varphi} = 1\,,$ so that
 $\sigma^{_m} \big( \kappa_Q (u)\big) = \widetilde{\tau^{^{{\rm b},U,m }}_Q}(u) \,.$ It remains to prove that this  {\it functorial\/} section fulfills the commutativity of the corresponding diagram £4.5.12; since we already have the commutativity of diagram £4.5.25, it suffices to get the commutativity of the following diagram
$$\matrix{\F (R) & \buildrel \mu_{_R}^{_m}\over \too 
&\widetilde{\L^{^{{\rm b},U,m }}} (R)\cr
\hskip-20pt{\scriptstyle {\rm id}_{\F (R)}}\big\downarrow&
& \big\downarrow {\scriptstyle a_{(\beta^{_m}_{\tilde\varphi})^{-1}}}\hskip-30pt\cr
\F (R) & \buildrel \mu_{_R}^{_m}\over \too 
&\widetilde{\L^{^{{\rm b},U,m }}} (R)\cr}
\eqno £4.5.44$$
which follows from the fact that $\beta^{_m}$ is {\it stable\/} and therfore 
$(\beta^{_m}_{\tilde\varphi})^{-1}$ fixes the image of $\mu_{_R}^{_m}\,.$

\smallskip
We can  modify this correspondence in order to get an {\it $\F\-$locality functorial section\/}; indeed, for any $\F_{\!P}\-$morphism $\kappa_{Q,R} (u) \colon R\to Q$ where 
$u$ belongs to $\T_P (Q,R)\,,$ the $\widetilde{\L^{^{{\rm b}, U,m}}}(Q,R)\-$morphisms 
$\sigma^{_m} \big(\kappa_{Q,R} (u)\big)$ and $\widetilde{\tau^{^{{\rm b},U,m}}_{_{Q,R}}}(u)$  both lift $\kappa_{Q,R} (u)\in \F (Q,R)\,;$ thus,  the 
{\it divisibility\/} of $ \widetilde{\L^{^{{\rm b}, U,m}}}$ guarantees the existence and the uniqueness of $\delta_{\kappa_{Q,R} (u)}\in {\rm Ker}(\widetilde{\pi^{^{{\rm b},U,m}}_R})$ fulfilling
$$\widetilde{\tau^{^{{\rm b},U,m}}_{_{Q,R}}}(u) = \sigma^{_m} \big(\kappa_{Q,R} (u)\big).\delta_{\kappa_{Q,R} (u)}
\eqno £4.5.45\phantom{.}$$
and, since we have 
$\sigma^{_m} \big( \kappa_Q (w)\big) = \widetilde{\tau^{^{{\rm b},U,m }}_Q}(w)$ for any 
$w\in Q\,,$  it is quite clear that $\delta_{\kappa_{Q,R} (u)}$ only depends on the class of
$\kappa_{Q,R} (u)$ in $\tilde\F (Q,R)$

\smallskip
 For a second $\F_{\! P}\-$morphism $\kappa_{R,T} (v)\,\colon T\to R\,,$ setting 
 $\xi = \kappa_{R,T} (u)$ and $\eta = \kappa_{R,T} (v)$ we get
$$\eqalign{\sigma^{_m} (\xi\circ \eta)\. \delta_{\tilde\xi\circ\tilde\eta} 
&= \widetilde{\tau^{^{{\rm b},U,m}}_{_{Q,T}}}(uv) = \widetilde{\tau^{^{{\rm b},U,m}}_{_{Q,R}}} (u) \. \widetilde{\tau^{^{{\rm b},U,m}}_{_{R,T}}}(v)\cr
&= \sigma^{_m} (\xi)\. \delta_{\tilde\xi}\.\sigma^{_m} (\eta).  \delta_{\tilde\eta}\cr
&= \sigma^{_m} (\xi\circ \eta)\.\big(\widetilde\Ker (\widetilde{\pi^{^{{\rm b},U,m}}})(\tilde\eta)\big)(\delta_{\tilde\xi})\. \delta_{\tilde\eta}\cr}
\eqno £4.5.46;$$
then, once again  the {\it divisibility\/} of $ \widetilde{\L^{^{{\rm b}, U,m}}}$ forces
$$\delta_{\tilde\xi\circ \tilde\eta} = \big(\widetilde\Ker (\widetilde{\pi^{^{{\rm b},U,m}}})
(\tilde\eta)\big)(\delta_{\tilde\xi})\. \delta_{\tilde\eta}
\eqno £4.5.47\phantom{.}$$
and, since ${\rm Ker}(\widetilde{\pi^{^{{\rm b},U,m}}_T})$ is abelian, with the additive notation we obtain
$$0 =\big(\widetilde\Ker (\widetilde{\pi^{^{{\rm b},U,m}}}) (\tilde\eta)\big)(\delta_{\tilde\xi}) -  \delta_{\tilde\xi\circ \tilde\eta}+ \delta_{\tilde\eta}
\eqno £4.5.48.$$

\smallskip
That is to say,  denoting by $\frak i\,\colon \tilde\F_{\! P}\i \tilde\F$ the obvious {\it inclusion functor\/}, the correspondence $\delta$ sending any $\tilde\F_{\! P}\-$morphism 
$\tilde\xi\,\colon R\to Q$ to~$\delta_{\tilde\xi}$ defines a {\it $1\-$cocycle\/} in 
$\Bbb C^1\big(\tilde\F_{\! P},\widetilde\Ker (\widetilde{\pi^{^{{\rm b},U,m}}}) \circ\frak i\big)\,;$ but, since the category $\tilde\F_{\! P}$ has a final object, we actually have [8,~Corollary~A4.8]
$$\Bbb H^1\big(\tilde\F_{\! P},\widetilde\Ker (\widetilde{\pi^{^{{\rm b},U,m}}})
\circ\frak i\big) = \{0\}
\eqno £4.5.49;$$
consequently, we obtain $\delta = d_0 (w)$ for some element 
$w = (w_Q)_{Q\i P}$ in 
$$\Bbb C^0\big(\tilde\F_{\! P},\widetilde\Ker (\widetilde{\pi^{^{{\rm b},U,m}}}) \circ\frak i\big) = \Bbb C^0\big(\tilde\F,\widetilde\Ker (\widetilde{\pi^{^{{\rm b},U,m}}}))\big)
\eqno £4.5.50.$$
In conclusion, equality~£4.5.45 becomes
$$\eqalign{\widetilde{\tau^{^{{\rm b},U,m}}_{_{Q,R}}}(u) &= 
\sigma^{_m} \big(\kappa_{Q,R} (u)\big) \.\big(\widetilde\Ker (\widetilde{\pi^{^{{\rm b},U,m}}})
(\widetilde{\kappa_{Q,R} (u)}\big)(w_Q)\.w_R^{-1}\cr
& = w_Q\.\sigma^{_m} \big(\kappa_{Q,R} (u)\big)\.w_R^{-1}\cr}
\eqno £4.5.51\phantom{.}$$
and therefore  the new correspondence  sending $\varphi\in \F (Q,R)$ to  
$w_Q\.\sigma^{_m} \big(\varphi)\.w_R^{-1}$ defines  a {\it $\F\-$locality functorial section\/}   
of $\widetilde{\pi^{^{{\rm b},U,m}}}\,.$ From now on, we still denote by 
$\sigma^{_m}$ this  {\it $\F\-$locality functorial section\/}   
of $\widetilde{\pi^{^{{\rm b},U,m}}}\,.$

\smallskip
Let $\sigma'^{_m}\,\colon \F\to \widetilde{\L^{^{{\rm b},U,m}}}$ be another $\F\-$locality functorial section\    of $\widetilde{\pi^{^{{\rm b},U,m}}}\,;$ arguing by induction on 
$\vert\C_P - \N\vert$ and on $m\,,$ and up to natural $\F\-$isomor-phisms, 
we clearly may assume that  $\sigma'^{_m}$ also lifts $\sigma^{_{m-1}}\,;$ in this case, for any $\F\-$morphism $\varphi\,\colon R\to Q\,,$ we have 
$\sigma'^{_m} (\varphi) = \sigma^{_m} (\varphi) \.  \varepsilon^{_m}_\varphi$  for some 
$\varepsilon^{_m}$ in~$(\tilde\frak l^{^{U,m - 1}}_\F/\tilde\frak l^{^{U,m}}_\F) (R)\,;$ that is to say,  as above we  get a correspondence sending any {\it $\F\-$chain\/} 
 $\frak q\,\colon \Delta_1\to \F$   to~$ \varepsilon^{_m}_{\frak q (0\bullet 1),}$ 
 in $(\tilde\frak l^{^{U,m - 1}}_\F/\tilde\frak l^{^{U,m}}_\F) \big(\frak q(0)\big)$ and  we claim that this correspondence determines an  {\it  $\F_{\!P}\-$stable\/} element $\varepsilon^{_m}$ of 
 $ \Bbb C^1 \big(\tilde\F,\tilde\frak l^{^{U,m - 1}}_\F \! /\tilde\frak l^{^{U,m}}_\F\big) $ 
 [8,~A3.17].

\smallskip
 Indeed, for another $\F_{\!P}\-$isomorphic {\it $\F\-$chain\/} $\frak q'\,\colon \Delta_1\to \F$ and a {\it natural $\F_{\!P}\-$isomorphism\/} $\nu\,\colon \frak q\cong\frak q'\,,$
as in~£4.5.34 above setting  
$$\eqalign{R = \frak q(0)\  ,\ R'  = \frak q' (0)\ & ,\ Q = \frak q (1)\ ,\ Q' = \frak q' (1)\cr
\varphi = \frak q(0\bullet 1)&\quad,\quad \!\!\varphi' = \frak q'(0\bullet 1)\cr 
\nu_0 = \kappa_{R',R}(v)  &\qq 
\nu_1 = \kappa_{Q',Q}(u)\cr}
\eqno £4.5.52,$$
from £4.5.29 we get
$$\eqalign{ \sigma'^{_m} (\varphi')  &=  \kappa_{Q',Q}(u)\. \sigma^{_m} (\varphi) \.  \varepsilon^{_m}_\varphi \. \kappa_{R',R}(v)^{-1}\cr
& = \sigma^{_m} (\varphi') \. \big((\tilde\frak l^{^{U,m - 1}}_\F/\tilde\frak l^{^{U,m}}_\F)(\widetilde{ \kappa_{R',R}(v)}^{-1})\big)(\varepsilon^{_m}_\varphi)\cr
 \sigma'^{_m} (\varphi')  &=  \sigma^{_m} (\varphi') \. \varepsilon^{_m}_{\varphi'}\cr}
\eqno £4.5.53$$
and the divisibility of  $\widetilde{\L^{^{{\rm b},U,m}}}$  forces
$$\varepsilon^{_m}_{\varphi'} =  \big((\tilde\frak l^{^{U,m - 1}}_\F/\tilde\frak l^{^{U,m}}_\F)(\widetilde{ \kappa_{R',R}(v)}^{-1}) \big)(\varepsilon^{_m}_\varphi)
\eqno £4.5.54;$$
this proves  that the correspondence $\varepsilon^{_m}$ sending $\varphi$ 
to~$\varepsilon^{_m}_\varphi$ is {\it $\F_{\!P}\-$stable\/} and, in particular, that 
$\varepsilon^{_m}_\varphi$ only depends on the corresponding $\tilde\F\-$morphism 
$\tilde\varphi\,,$  thus we set $\varepsilon^{_m}_{\tilde\varphi}= \varepsilon^{_m}_\varphi,.$

\smallskip
On the other hand, considering the usual differential map
$$d_1 : \Bbb C^1 \big(\tilde\F,\tilde\frak l^{^{U,m - 1}}_\F/\tilde\frak l^{^{U,m}}_\F\big)\too 
\Bbb C^2 \big(\tilde\F,\tilde\frak l^{^{U,m - 1}}_\F/\tilde\frak l^{^{U,m}}_\F\big)
\eqno £4.5.55,$$
we claim that $d_1 (\varepsilon^{_m}) = 0\,;$ indeed, for a second $\F\-$morphism 
$\psi\,\colon T\to R$ we get
$$\eqalign{\sigma'^{_m} (\varphi) \. \sigma'^{_m} (\psi) 
&= \sigma^{_m} (\varphi)\. \varepsilon^{_m}_\varphi\. \sigma^{_m} (\psi)\. \varepsilon^{_m}_\psi\cr
&= \sigma^{_m} (\varphi\circ \psi)\. \big((\tilde\frak l^{^{U,m - 1}}_\F/\tilde\frak l^{^{U,m}}_\F)(\tilde\psi)\big)(\varepsilon^{_m}_{\tilde\varphi}) \.  \varepsilon^{_m}_\psi \cr
\sigma'^{_m} (\varphi) \. \sigma'^{_m} (\psi) 
 &=  \sigma^{_m} (\varphi\circ \psi)\. \varepsilon^{_m}_{\varphi\circ \psi}\cr }
\eqno £4.5.56\phantom{.}$$
and   the {\it divisibility\/} of $\widetilde{\L^{^{{\rm b},U,m}} }$ forces
$$\big((\tilde\frak l^{^{U,m - 1}}_\F/\tilde\frak l^{^{U,m}}_\F)(\tilde\psi)\big)(\varepsilon^{_m}_{\tilde\varphi}) \.  \varepsilon^{_m}_{\tilde\psi } =  \varepsilon^{_m}_{\tilde\varphi\circ \tilde\psi}
\eqno £4.5.57;$$
since ${\rm Ker}(\widetilde{\pi_T^{^{{\rm b},U,m  }}})$ is Abelian, with the additive notation we obtain
$$0 = \big((\tilde\frak l^{^{U,m - 1}}_\F/\tilde\frak l^{^{U,m}}_\F)(\tilde\psi)\big)(\varepsilon^{_m}_{\tilde\varphi}) -\varepsilon^{_m}_{\tilde\varphi\circ \tilde\psi}
+  \varepsilon^{_m}_{\tilde\psi }
\eqno £4.5.58,$$
proving our claim.

\smallskip
At this point, it follows from equalities £4.5.5 and £4.5.9  that  $\varepsilon^{_m} = d_0 (y)$ for some {\it stable\/} element $y = (y_Q)_{Q\i P}$ in 
$\Bbb C^0 \big(\tilde\F,\tilde\frak l^{^{U,m - 1}}_\F/\tilde\frak l^{^{U,m}}_\F\big)\,;$ 
that is to say, with the notation above we get
$$\varepsilon^{_m}_{\tilde\varphi} = 
\big((\tilde\frak l^{^{U,m - 1}}_\F/\tilde\frak l^{^{U,m}}_\F)(\tilde\varphi)\big)(y_Q)\.y_R^{-1}
\eqno £4.5.59;$$
 hence,  we obtain
$$\eqalign{\sigma'^{_m} (\varphi) = \sigma^{_m} (\varphi) \. \big((\tilde\frak l^{^{U,m - 1}}_\F/\tilde\frak l^{^{U,m}}_\F)(\tilde\varphi)\big)(y_Q)\.y_R^{-1} =  y_Q\. \sigma^{_m} (\varphi)\. y_R^{-1}\cr}
\eqno £4.5.60,$$
which amounts to saying that  $\sigma'^{_m}$ is naturally $\F\-$isomorphic to 
$\sigma^{_m}\,.$ We are done
\eject

\bigskip
\noindent
{\bf Corollary £4.6.} {\it  There exists a perfect $\F\-$locality $\P\,.$\/} 

\par
\noindent
{\bf Proof:} Denote by $\bar\P$ the converse image in $\L^{^{\rm b}}$ of the image of $\F$
 in $\widetilde{\L^{^{\rm b}}}$ by a section of $\widetilde{\pi^{_b}}\,;$ since 
 $\hat\tau (\frak c^\frak f_\F)$ is contained in the image of $\tau^{_{\rm b}}\!\!\,,$ we still have a functor $\tau^{_{\rm b}}\,\colon \T_P\to \bar\P\,;$ thus, together with the restriction of 
 $\pi^{_{\rm b}}$ to $\bar\P\,,$ $\bar\P$ becomes an {\it $\F\-$locality\/} and, since 
 $\L^{^{\rm b}}$ is {\it coherent\/}, $\bar\P$ is coherent too.
 
 \smallskip
 We claim that  $\bar\P^{^{\rm sc}}$ is a {\it perfect $\F^{^{\rm sc}}\-$locality\/}; indeed, for any 
 $\F\-$selfcentrali-zing subgroup $Q$ of $P$ fully normalized in $\F\,,$ since $C_P (Q)/F_{C_\F (Q)} = Z(Q)$  we have a {\it group extension\/} (cf.~£4.3)
 $$1\too Z(Q)\too \bar\P (Q)\too \F (Q)\too 1
 \eqno £4.6.1\phantom{.}$$
 together with an injective group homomorphism 
 $\tau^{_{\rm b}}_Q\,\colon N_P(Q)\to \bar\P (Q)\,;$ consequently, it follows from [8, 18.5] that
 $\bar\P (Q)$ is the {\it $\F\-$localizer\/} of $Q\,;$ thus, by the very definition in [8, 17.4 and 17.13], $\bar\P^{^{\rm sc}}$ is a {\it perfect $\F^{^{\rm sc}}\-$locality\/}.
 
 \smallskip
 But, in [8, Ch. 20] we prove that any  {\it perfect $\F^{^{\rm sc}}\-$locality\/} can be extended to a unique  {\it perfect $\F\-$locality\/} $\P\,.$ We are done.
 
 \medskip
£4.7. The uniqueness of the perfect $\F\-$locality is an easy consequence of the following theorem; the proof of this result follows the same pattern than the proof of Theorem £4.5, but
we firstly  need the following lemmas.

\bigskip
\noindent
{\bf Lemma £4.8.} {\it Let $(\tau,\P,\pi)$ be  a perfect $\F\-$locality and 
$\hat\varphi ∞\colon Q\to P$ be a $\P\-$morphism such that $\pi_{\hat\varphi} (Q)$ is fully normalized in $\F\,.$ Then there is a  $\P\-$morphism $\hat\zeta\,\colon N_P(Q)\to P$ such that $\hat\varphi = \hat\zeta\. \tau_{N_P(Q),Q} (1)\,.$\/}
\medskip
\noindent
{\bf Proof:} Denoting by $\varphi$ the image of $\hat\varphi$ in $\F (P,Q)\,,$ it follows from
[8, 2.8.2]  that there is an $\F\-$morphism $\zeta\,\colon N_P (Q)\to P$ extending $\varphi\,;$
then, lifting $\zeta$ to $\hat\zeta$ in $\P \big(P,N_P (R)\big)\,,$ it is clear that the 
$\P\-$morphisms $\hat\zeta\. \tau_{N_P(Q),Q} (1)$ and $\hat\varphi$ have the same image\
$\varphi$ in $\F (P,Q)$ and therefore, by the very definition of $\P$ in [8,~17.13], there is 
$z\in C_P (Q)$ such that $\hat\zeta\. \tau_{N_P(Q),Q} (1)\. \tau_Q (z) = \hat\varphi\,;$
bur, it is clear that
$$ \tau_{N_P(Q),Q} (1)\. \tau_Q (z) =  \tau_{N_P(Q),Q} (z) =  \tau_{N_P(Q)} (z)\. \tau_{N_P(Q),Q} (1)
\eqno £4.8.1;$$
consequently, $\hat\zeta\.  \tau_{N_P(Q)} (z)$ extends $\hat\varphi$ in $\P\,.$ We are done

\bigskip
\noindent
{\bf Lemma £4.9.} {\it Let $(\tau,\P,\pi)$ be  a perfect $\F\-$locality. For any subgroup $Q$ of $P$ there is a group homomorphism $\hat\mu_Q\,\colon \P(Q)\to \L^{^{\rm b}} (Q)$ fulfilling 
$\hat\mu_Q\circ \tau_Q = \tau^{_{\rm b}}_Q\,.$}

\medskip
\noindent
{\bf Proof:} Since we can choose an $\F\-$isomorphism $\theta\,\colon Q\cong Q'$ such that       $Q'$ is fully normalized in $\F$ and $\theta$ can be lifted to $\P (Q',Q)$ and to
$\L^{^{\rm b}} (Q',Q)\,,$ we may assume that $Q$ is fully normalized in $\F\,.$ 

\smallskip
We apply [8, Lemma 18.8] to the groups $\P (Q)$ and $\L^{^{\rm b}} (Q)\,,$ to the normal 
$p\-$subgroup ${\rm Ker} (\pi^{_{\rm b}}_Q )$ of $\L^{^{\rm b}} (Q)$ and to the
group homomorphism $\tau_Q$ from $\P (Q)$ to 
$\F (Q)\cong \L^{^{\rm b}} (Q)/{\rm Ker} (\pi^{_{\rm b}}_Q )\,.$ We consider the group homomorphism $\tau^{_{\rm b}}_Q\,\colon N_P (Q)\to \L^{^{\rm b}} (Q)$ and, for any subgroup $R$ of $N_P (Q)$ and any $\hat\alpha\in \P (Q)$ such that 
$\hat\alpha\. \tau_Q (R)\. \hat\alpha^{-1}\i \tau_Q \big(N_P (Q)\big)\,,$ it~follows from 
[8, 2.10.1] that there exists $\zeta\in \F \big(N_P (Q),Q\. R\big)$ extending the image of 
$\hat\alpha$ in $\F (Q)\,;$ then, it follows from [8, 17.11.2] that there exists 
$x\in  \L^{^{\rm b}} (Q)$ fulfilling 
$${\tau^{_{\rm b}}_Q} \big(\zeta (v)\big) ={}^x {\tau^{_{\rm b}}_Q} (v)
\eqno £4.9.1\phantom{.}$$
for any $v\in Q .R\,.$ That is to say, condition 18.8.1 in  [8, Lemma 18.8]  is fulfilled and therefore this lemma proves the existence of $\hat\mu_Q$ as announced.

.

\bigskip
\noindent
{\bf Theorem £4.10.} {\it For any perfect $\F\-$locality $\P$ there exists a unique natural 
$\F\-$isomorphism class of $\F\-$locality functors to $\L^{^{\rm b}}\!\!\,.$\/}

\medskip
\noindent
{\bf Proof:}  Let $\P$ be a perfect $\F\-$locality with the structural functors
$$\tau : \T_P \too \P \qq \pi : \P \too \F
\eqno £4.10.1\phantom{.}$$
and for any subgroups $Q$ of $P$ and $R$ of $Q$ we set $i_R^Q = \tau_{Q,R} (1)\,.$
We consider the filtration of the basic $\F\-$locality introduced in section~3 and then argue by induction. That is to say, recall that we denote by $\C_P$ a set of representatives for the set of $P\-$conjugacy classes of subgroups $U$ of $P$ (cf. £2.13); now, for any subset $\N$ of $\C_P$ fulfilling condition £3.2.1, we have the functor 
$\tilde\frak k^{^\N}_\F :  \tilde\F\to \Ab$  (cf. £3.2) and we consider  the quotient $\F\-$locality 
$\L^{^{{\rm b},\N}} = \L^{^{\rm b}} /( \tilde\frak k^{^\N}_\F\circ \tilde\pi^{_{\rm b}} )$ with the structural functors
$$\tau^{_{\rm b,\N}} : \T_P\too \L^{^{{\rm b},\N}}
\qq \pi^{_{\rm b,\N}}  :  \L^{^{{\rm b},\N}} \too \F
\eqno £4.10.2.$$
Note that if $\N = \emptyset$ then $\L^{^{{\rm b},\N}} = \L^{^{\rm b}} \,;$ hence,  arguing by induction on $\vert\C_P - \N\vert\,,$  it suffices to prove the existence of a unique natural $\F\-$isomorphism class of $\F\-$locality functors from $\P$ to $\L^{^{{\rm b},\N}} $.

\smallskip
Moreover,  if $\N = \C_P$ then $ \tilde\frak k^{^\N}_\F =  \tilde\frak k^{^{\rm b}}_\F$ and therefore $\L^{^{{\rm b},\N}} = \F\,,$ so that we may assume that $\N\not= \C_P\,;$ in this situation, we fix a minimal element $U$ in $\C_P - \N\,,$ setting $\M = \N\cup \{U\}$ and 
$\tilde\frak k^{^U}_\F = \tilde\frak k^{^\M}_\F/ \tilde\frak k^{^\N}_\F\,;$
 for any $m\in \Bbb N$ we simply denote by~$\tilde\frak l^{^{U,m}}_\F$ the converse image of~$p^m\. \tilde\frak k^{^U}_\F$  in~$\tilde\frak k^{^\M}_\F\,;$  set 
$\L^{^{{\rm b},U,m}} = \L^{^{\rm b}}/\tilde\frak l^{^{U,m}}_\F$ and, coherently, denote by 
$\pi^{_{\rm b,U,m}}$ and $\tau^{_{\rm b,U,m}}$  the corresponding structural functors. Note that, by £3.8 and
£3.10.5 above we get 
$$\tilde\frak l^{^{U,m}}_\F/\tilde\frak l^{^{U,m +1}}_\F\cong \tilde\frak r^{U,\circ}_{\F,m}
\eqno £4.10.3\phantom{.}$$
and in particular, by Corollary £3.13, for any $n\in \Bbb N$ we still get
$$\Bbb H^{^n}_* (\tilde\F,\tilde\frak l^{^{U,m}}_\F/\tilde\frak l^{^{U,m +1}}_\F) = \{0\}
\eqno £4.10.4.$$
\eject

\smallskip
As above, we denote by $\C_\F$ a set of representatives,  fully normalized in~$\F\,,$ for the 
$\F\-$isomorphism classes of subgroups of $P$ and, for any subgroup $Q$ in~$\C_\F\,,$  we choose a  group homomorphism $\hat\mu_Q\,\colon \P (Q)\to {\L^{^{\rm b}}} (Q)$ as in Lemma £4.9 above and, for any $m\in \Bbb N\,,$ simply denote by $\hat\mu_Q^{_m}$ the corresponding group homomorphism from $\P (Q)$ to ${\L^{^{{\rm b},U,m}}} (Q)\,.$ For any 
$\F\-$morphism $\varphi\,\colon R\to Q$ denote by $\P (Q)_{\varphi}$ and by 
${\L^{^{{\rm b},U,m}}} (Q)_\varphi$ the respective stabilizers of $\varphi (R)$ in 
$\P (Q)$ and in ${\L^{^{{\rm b},U,m}}} (Q)\,.$ As above, for any $\hat\varphi\in \P (Q,R)$ and any
$x^{_m}\in  {\L^{^{{\rm b},U,m}}} (Q,R)$ we have  group homomorphisms
$$a_{\hat\varphi} : \P (Q)_\varphi\too \P (R) \qq   a_{x^{_m}} : {\L^{^{{\rm b},U,m}}} (Q)_\varphi \too {\L^{^{{\rm b},U,m}}} (R)
\eqno £4.10.5\.$$
For any  subgroups $Q$ and $R$ in $\C_\F\,,$  we choose as above a set of representati-ves 
$\P_{\!Q,R}$ for the set of double classes $\P (Q)\backslash \P (Q,R)/\P (R)$ such that, for any 
$\hat\varphi$ in $\P_{\!Q,R}\,,$ denoting by $\varphi$ its image in $\F (Q,R)\,,$
$\F_{\!P} (Q)$ contains a Sylow $p\-$subgroup of $\F (Q)_\varphi$ and 
$a_\varphi \big(\F_{\!P} (Q)_\varphi\big)$ is contaioned in $\F_{\!P} (R)\,;$ of course, we choose 
$\P_{Q,Q} =\{\tau_Q (1)\}\,.$

\smallskip
With all this notation and arguing by induction on $\vert\C_P - \N\vert$ and on~$m\,,$   we will prove that there is a functor
$$\lambda^{^{m}} : \P\too {\L^{^{{\rm b},U,m }}}
\eqno £4.10.6\phantom{.}$$
 such that, for any $Q\in C_\F$ and any $u\in Q\,,$ we have 
 $\lambda^{^{m}} \big(\tau_Q (u)\big) = {\tau^{^{{\rm b},U,m }}_Q} (u)\,,$ and that, for any   groups $Q$ and $R$ in~$\C_\F$, and any  $\hat\varphi$ in~$\P_{Q,R}\,,$ denoting by $\varphi$ its image in $\F (Q,R)\,,$ we have the commutative diagram
$$\matrix{\P (Q)_\varphi & \buildrel \hat\mu_{_Q}^{_m}\over \too 
&{\L^{^{{\rm b},U,m }}} (Q)_\varphi\cr
\hskip-20pt{\scriptstyle a_{\hat\varphi}}\big\downarrow&
& \big\downarrow {\scriptstyle a_{\lambda^{^{m}} (\hat\varphi)}}\hskip-20pt\cr
\P(R) & \buildrel  \hat\mu_{_R}^{_m}\over \too  &{\L^{^{{\rm b},U,m }}} (R)\cr}
\eqno £4.10.7.$$
 Since we have $ {\pi^{_{{\rm b},U,0}}} =  {\pi^{_{{\rm b},\M}}}$ 
and $\vert\M\vert = \vert\N\vert + 1\,,$ by the induction hypothesis we actually may assume that $m\not= 0$ and that we have a functor 
$$\lambda^{^{m - 1}} : \P \too {\L^{^{{\rm b},U,m - 1}}}
\eqno £4.10.8\phantom{.}$$
which fulfills the conditions above.

\smallskip
As above, for any $\hat\varphi \in \P_{\!Q,R}\,,$  denoting by $\varphi$ its image in $\F (Q,R)\,,$ it follows from [8, Proposition~2.11], applied to the  inverse $\varphi^*$ of the isomorphism $\varphi_*\,\colon R\cong \varphi (R)$ induced by 
$\varphi\,,$ that there exists an $\F\-$morphism $\zeta\,\colon N_P (Q)_\varphi\to N_P (R)$ fulfilling $\zeta\big(\varphi (v)\big) = v$ for any $v\in R\,,$ so that we easily get the following commutative diagram`
$$\matrix{N_P (Q)_\varphi &\buildrel \tau_Q \over 
\too & \P (Q)_\varphi\cr
\hskip-15pt {\scriptstyle\ \zeta} \big\downarrow &
&\big\downarrow{\scriptstyle a_{\hat\varphi} }\hskip-10pt  \cr
N_P (R) &\buildrel \tau_R \over \too&\P (R)&\cr}
\eqno £4.10.9;$$
note that, if $Q = R$ and $\hat\varphi = \tau_Q (u)$ for some $u\in Q\,,$ we may assume that\break  $\zeta = \kappa_{N_P (Q)} (u)\,.$ In particular, since $\lambda^{^{m - 1}}$ fulfills the corresponding commutative diagram £4.10.7,  we still get the following commutative diagram 
$$\matrix{N_P (Q)_\varphi &\buildrel {\tau^{^{{\rm b},U,m -1}}_Q} \over 
{\hbox to 40pt{\rightarrowfill}} & {\L^{^{{\rm b},U,m - 1}}} (Q)_\varphi\cr
\hskip-10pt {\scriptstyle \zeta} \big\downarrow &
&\big\downarrow{\scriptstyle a_{\lambda^{_{m -1}} (\hat\varphi)} }\hskip-20pt  \cr
N_P (R) &\buildrel {\tau^{^{{\rm b},U,m -1}}_R} \over 
{\hbox to 40pt{\rightarrowfill}} & {\L^{^{{\rm b},U,m - 1}}} (R)&\cr}
\eqno £4.10.10$$

\smallskip
With the notation above ,the first step is to choose a suitable lifting $\widehat{\lambda^{^{m -1}} (\hat\varphi) }$ of  $\lambda^{_{m - 1}}(\hat\varphi)$ in  ${\L^{^{{\rm b},U,m }}} (Q,R)\,.$ Choosing a lifting 
$\hat\zeta$ of $\zeta$ in  the obvious stabilizer 
$\P  \big(N_P (R),N_P (Q)_\varphi\big)_{R,\varphi (R)}\,,$ we start by choosing a lifting  $\widehat{\lambda^{^{m -1}} (\hat\zeta) }$ of $\lambda^{^{m -1}} (\hat\zeta)$ in 
${\L^{^{{\rm b},U,m }}} \big(N_P (R),N_P (Q)_\varphi\big)_{R,\varphi (R)}\,;$
thus, by the {\it coherence\/} of ${\L^{^{{\rm b},U,m }}}$ (cf. (Q)), for any 
$u\in N_P (Q)_\varphi$ we have
$$\widehat{\lambda^{^{m -1}} (\hat\zeta) }\.{\tau^{^{{\rm b},U,m }}_{N_P (Q)_\varphi}} (u) = {\tau^{^{{\rm b},U,m }}_{N_P (R)}}\big(\zeta (u)\big)\.\widehat{\lambda^{^{m -1}} (\hat\zeta) }
\eqno £4.10.11;$$
moreover, by the {\it divisibility\/} of ${\L^{^{{\rm b},U,m }}}$ (cf. £2.4), we find
$z_{\hat\varphi} \in  {\L^{^{{\rm b},U,m }}} \big(R,\varphi(R)\big)$ fulfilling
$$\widehat{\lambda^{^{m -1}} (\hat\zeta) }\.{\tau^{^{{\rm b},U,m }}_{N_P (Q)_\varphi,\varphi(R)}} (1) = {\tau^{^{{\rm b},U,m }}_{N_P (R),R}} (1)\. z_{\hat\varphi}
\eqno £4.10.12;$$
similarly, denoting by $\hat\varphi^*\,\colon \varphi (R)\cong R$ the  
$\P\-$isomorphism determined by~$\hat\varphi\,,$ $\lambda^{^{m -1}} (\hat\zeta)$ restricts to 
$\lambda^{^{m -1}} (\hat\varphi^*)$ and it is easily checked that $z_{\hat\varphi}$ lifts 
$\lambda^{^{m -1}} (\hat\varphi^*)$ to~$ {\L^{^{{\rm b},U,m }}} \big(R,\varphi(R)\big)\,;$
consequently, $\widetilde{\lambda^{^{m -1}} (\hat\varphi) } = 
{\tau^{^{{\rm b},U,m }}_{Q\. {\hat\varphi} (R)}} (1)\. z_\varphi^{-1}$ lifts 
$\lambda^{^{m -1}} (\hat\varphi)$\break to $ {\L^{^{{\rm b},U,m }}} (Q,R)\,.$.

\smallskip
Then, from £4.10.11 and £4.10.12 above,  for any 
$u\in N_P (Q)_\varphi$ we get
$$\eqalign{\widehat{\lambda^{^{m -1}} (\hat\zeta) }\.
{\tau^{^{{\rm b},U,m }}_{N_P (Q)_\varphi}} (u)&\.
.{\tau^{^{{\rm b},U,m }}_{N_P (Q)_\varphi,\varphi(R)}} (1) = 
{\tau^{^{{\rm b},U,m }}_{N_P (R),R}} (1)\. 
z_{\hat\varphi} \. {\tau^{^{{\rm b},U,m }}_{\varphi (R)}} (u)\cr
\Vert\cr
{\tau^{^{{\rm b},U,m }}_{N_P (R)}}\big(\zeta (u)\big)\.\widehat{\lambda^{^{m -1}} (\hat\zeta) }&\. .{\tau^{^{{\rm b},U,m }}_{N_P (Q)_\varphi,\varphi(R)}} (1)\cr
&= {\tau^{^{{\rm b},U,m }}_{N_P (R),R}} (1)\.
{\tau^{^{{\rm b},U,m }}_R}\big(\zeta (u)\big)\. z_{\hat\varphi}\cr}
\eqno £4.10.13\phantom{.}$$
and therefore we still get $z_{\hat\varphi}\. {\tau^{^{{\rm b},U,m }}_{\varphi (R)}} (u) = {\tau^{^{{\rm b},U,m }}_R}\big(\zeta (u)\big)\. z_{\hat\varphi}\,,$ so that
$$ {\tau^{^{{\rm b},U,m }}_Q} (u)\. \widetilde{\lambda^{^{m -1}} (\hat\varphi) } = 
{\lambda^{^{m -1}} (\hat\varphi) }\. {\tau^{^{{\rm b},U,m }}_R}\big(\zeta (u)\big)
\eqno £4.10.14\phantom{.}$$
or, equivalently, we have $a_{ \widetilde{\lambda^{^{m -1}} (\hat\varphi) } } \big( {\tau^{^{{\rm b},U,m }}_Q} (u)\big) = {\tau^{^{{\rm b},U,m }}_R}\big(\zeta (u)\big)\,.$

\smallskip
At this point, we will apply the uniqueness part of [8, Lemma 18.8] to the groups 
$\P (Q)_\varphi$ and $ {\L^{^{{\rm b},U,m }}} (R)\,,$ to the kernel of the canonical homomorphism
from $ {\L^{^{{\rm b},U,m }}} (R)$ to $ {\L^{^{{\rm b},U,m -1}}} (R)\,,$  and to the composition of group homomorphisms
$$a_{\lambda^{^{m -1}}(\hat\varphi)}\circ \hat\mu_Q^{_{m -1}} : \P (Q)_\varphi\too  {\L^{^{{\rm b},U,m -1}}} (Q)_\varphi\too {\L^{^{{\rm b},U,m -1}}} (R)
\eqno £4.10.15,$$
together with the composition of group homomorphisms
$${\tau^{^{{\rm b},U,m }}_R}\circ \zeta : N_P (Q)_\varphi\too N_P (R)\too 
 {\L^{^{{\rm b},U,m}}} (R)
 \eqno £4.10.16.$$
 Now, according to the commutative diagrams £4.10.7 for $m -1$ and £4.10.9, and to equality £4.10.14 above, the two group homomorphisms
 $$\eqalign{a_{\widetilde{\lambda^{^{m -1}}(\hat\varphi)}}\circ \hat\mu_Q^{_{m}} 
 &: \P (Q)_\varphi\too  {\L^{^{{\rm b},U,m}}} (Q)_\varphi\too {\L^{^{{\rm b},U,m}}} (R)\cr
\hat\mu_R^{_{m}}\circ  a_{\hat\varphi} &: \P (Q)_\varphi\too  \P (R)\too {\L^{^{{\rm b},U,m}}} (R)\cr}
 \eqno £4.10.17,$$
 both fulfill the conclusion of [8, Lemma 18.8]; consequently, according to this lemma, there is
 $k_{\hat\varphi}$ in the kernel of the canonical homomorphism from
  $ {\L^{^{{\rm b},U,m }}} (R)$ to $ {\L^{^{{\rm b},U,m -1}}} (R)$ such that, denoting by 
 ${\rm int}_{{\L^{^{{\rm b},U,m}}} (R)} (k_{\hat\varphi})$ the conjugation by $k_{\hat\varphi}$ in
  ${\L^{^{{\rm b},U,m}}} (R)\,,$ we have
  $${\rm int}_{{\L^{^{{\rm b},U,m}}} (R)} (k_{\hat\varphi)}\circ a_{\widetilde{\lambda^{^{m -1}}(\hat\varphi)}}\circ \hat\mu_Q^{_{m}} = \hat\mu_R^{_{m}}\circ  a_{\hat\varphi}
  \eqno £4.10.18;$$
  but, it is easily checked that 
  $${\rm int}_{{\L^{^{{\rm b},U,m}}} (R)} (k_{\hat\varphi})\circ a_{\widetilde{\lambda^{^{m -1}}(\hat\varphi)}} = a_{\widetilde{\lambda^{^{m -1}} }(\hat\varphi)}\. k_{\hat\varphi}^{-1}
\eqno £4.10.19.$$

 \smallskip
 Finally, we choose $\widehat{\lambda^{^{m -1}} (\hat\varphi)} 
 = \widetilde{\lambda^{^{m -1}} (\hat\varphi)}\. k_{\hat\varphi}^{-1}\,,$  lifting indeed 
 $\sigma^{_{m -1}} (\varphi)$
 to ${\L^{^{{\rm b},U,m}}} (Q,R)$ and, according to equalities £4.10.18 and £4.10.19, fulfilling the following commutative 
 diagram
$$\matrix{\P (Q)_\varphi & \buildrel \hat\mu_{_Q}^{_m}\over \too 
&{\L^{^{{\rm b},U,m }}} (Q)_\varphi\cr
\hskip-20pt{\scriptstyle a_{\hat\varphi}}\big\downarrow&
& \big\downarrow {\scriptstyle a_{\widehat{\lambda^{^{m -1}} (\hat\varphi)}}}\hskip-30pt\cr
\P (R) & \buildrel \hat\mu_{_R}^{_m}\over \too &{\L^{^{{\rm b},U,m }}} (R)\cr}
\eqno £4.10.20;$$
note that, if $Q = R$ and $\hat\varphi = \tau_Q (u)$ for some $u\in Q\,,$ this choice is com-patible  wtih $\widehat{\lambda^{^{m -1}} \big( \tau_Q (u)\big)} = {\tau^{^{{\rm b},U,m }}_Q} (u)\,.$ 
In particular, considering the action of $\P (Q)\times \P (R)\,,$ by composition on the left- and
on the right-hand, on $\P (Q,R)$ and on ${\L^{^{{\rm b},U,m }}} (Q,R)$ {\it via\/} $\hat\mu^{_m}_Q$ and $\hat\mu^{_m}_R\,,$ we have the inclusion of stabilizers
$$\big(\P (Q)\times \P (R)\big)_{\hat\varphi} \i \big(\P (Q)\times \P (R)\big)_{\lambda^{^{m -1}} (\hat\varphi)}
\eqno £4.10.21;$$
\eject
\noindent
indeed, it is quite clear that $(\hat\alpha,\hat\beta)\in \big(\P (Q)\times \P (R)\big)_{\hat\varphi} $ forces $\hat\alpha\in \P (Q)_\varphi\,;$  then, since $\hat\alpha\. \hat\varphi = \hat\varphi\.a_{\hat\varphi} (\hat\alpha)\,,$ we get $\hat\beta =  a_{\hat\varphi} (\hat\alpha)$ by the {\it divisibility\/} of $\P\,,$ and the inclusion above follows from the commutativity of diagram £4.10.20.

\smallskip
This allows us to choose a family of liftings   
$\{\widehat{\lambda^{_{m -1}} (\hat\varphi)}\}_{\hat\varphi}\,,$ where $\hat\varphi$ runs over the set of  $\P\-$morphisms, which is compatible with $\P\-$isomorphisms; precisely,  for any  pair of   subgroups $Q$ and $R$ in~$\C_\F\,,$ and any $\hat\varphi\in \P_{Q,R}\,,$ we choose  as above a lifting  $\widehat{\lambda^{^{m -1}} (\hat\varphi)}$ of $\lambda^{^{m -1}} (\hat\varphi)$ in  
${\L^{^{{\rm b},U,m }}} (Q,R)\,.$ Then, any subgroup $Q$ of $P$ determines a unique  $\bar Q$ in $\C_\F$ which is  $\F\-$isomorphic to~$Q$ and we choose a $\P\-$isomorphism 
$\hat\omega_Q\,\colon \bar Q\cong Q$ and a lifting  $x_Q\in {\L^{^{{\rm b},U,m }}} (Q,\hat Q)$ of the image $\omega_Q\in \F (Q,\bar Q)$ of $\hat\omega_Q\,;$  in particular, we choose
$\omega_{\bar Q} = \tau_{\bar Q} (1)$ and 
$x_{\bar Q} =  {\tau^{^{{\rm b},U,m }}_{\bar Q}} (1)\,.$ Thus, any $\P\-$morphism 
$\hat\varphi\,\colon R\to Q$ determines subgroups $\bar Q$ and $\bar R$ in~$\C_\F$ and an element $\bar{\hat\varphi}$ in  $\P_{\bar Q,\bar R}$ in such a way that there are 
$\hat\alpha_{\hat\varphi}  \in \P (\bar Q)$ and  $\hat\beta_{\hat\varphi}  \in \P (\bar R)$ fulfilling
 $$\hat\varphi = \hat\omega_Q\. \hat\alpha_{\hat\varphi} \.  \bar{\hat\varphi} \. \hat\beta^{-1}_{\hat\varphi} \. \hat \omega_R^{-1}
 \eqno £4.10.22\phantom{.}$$
 and we define
 $$\widehat{\lambda^{^{m -1}} ({\hat\varphi)}} = x_Q\.\hat \mu_{\bar Q}^{_m} (\hat\alpha_{\hat\varphi})\.
 \widehat{\lambda^{^{m -1}} (\bar{\hat\varphi)}} \.  \hat\mu_{\hat R}^{_m} (\hat\beta_{\hat\varphi})^{-1} \. x_R^{-1}
 \eqno £4.10.23;$$
once again,, if $Q = R$ and $\hat\varphi = \tau_Q (u)$ for some $u\in Q\,,$ we actually get
$\widehat{\lambda^{^{m -1}} \big( \tau_Q (u)\big)} = {\tau^{^{{\rm b},U,m }}_Q} (u)\,.$
 This definition does not depend on the choice of  
 $(\hat\alpha_{\hat\varphi}, \hat\beta_{\hat\varphi})$ since for another choice 
 $(\hat\alpha',\hat\beta')$ we clearly have $\hat\alpha' = \hat\alpha_{\hat\varphi }\. \hat\alpha''$ and  $\hat\beta' = \hat\beta_{\hat\varphi} \. \hat\beta''$ for a suitable $(\hat\alpha'', \hat\beta'')$  in~$\big(\P (\bar Q)\times \P (\bar R)\big)_{\hat\varphi}$ and it suffices to apply inclusion £4.10.21.

 \smallskip 
 Moreover, for any pair of  $\P\-$isomorphisms $\hat\zeta\,\colon Q\cong Q'$ and 
 $\hat\xi\,\colon R\cong R'\,,$  considering $\hat\varphi' = \hat\zeta\. \hat\varphi\. \hat\xi^{-1}$ we claim that  
 $$\widehat{\lambda^{^{m -1}} (\hat\varphi')} = \widehat{\lambda^{^{m -1}} (\hat\zeta)} \. \widehat{\lambda^{^{m -1}} (\hat\varphi) } \. \widehat{\lambda^{^{m -1}} (\hat\xi)}^{-1}
\eqno £4.10.24;$$
indeed, it is clear that $Q'$ also determines $\bar Q$ in  $\C_\F$ and therefore, if we have
$\hat\zeta = \hat\omega_Q\.\hat \alpha_{\hat\zeta}\. \hat\omega_{Q'}^{-1}$ then we obtain 
$\widehat{\lambda^{^{m -1}} (\hat\zeta)} = 
x_{Q'}\. \hat\mu_{\bar Q}^{_m} (\hat\alpha_{\hat\zeta})\. x_Q^{-1}\,;$
similarly,   if we have $\hat\xi = \hat\omega_R\. \hat\beta_{\hat\xi}\.\hat\omega_{R'}^{-1}$  we also obtain $\widehat{\lambda^{^{m -1}} (\hat\xi)}^{-1} = 
x_{R}\. \hat\mu_{\hat R}^{_m} (\hat\beta_{\hat\xi})^{-1}\. x_{R'}^{-1}\,;$
further, $\hat\varphi'$ also determines $\bar{\hat\varphi}$ in $\P_{\bar Q,\bar R}\,;$
consequently, we get
$$\eqalign{\widehat{\lambda^{^{m -1}} (\hat\zeta)} &\. \widehat{\lambda^{^{m -1}} (\hat\varphi) } \. \widehat{\lambda^{^{m -1}} (\hat\xi)}^{-1} \cr
&= (x_{Q'}\.\hat \mu_{\bar Q}^{_m} (\hat\alpha_{\hat\zeta})\. x_Q^{-1})\.
\widehat{\lambda^{^{m -1}} (\hat\varphi)}\. (x_{R}\.\hat \mu_{\hat R}^{_m} (\hat\beta_{\hat\xi})^{-1}\. x_{R'}^{-1})\cr
&= x_{Q'} \. \hat \mu_{\hat Q}^{_m} (\hat\alpha_{\hat\zeta} \. \hat\alpha_{\hat\varphi})\. \bar{\hat\varphi}\. \hat\mu_{\hat R}^{_m} (\hat\beta_{\hat\varphi}^{-1}\. \hat\beta_{\hat\xi}^{-1})\. x_{R'}^{-1} = \widehat{\lambda^{^{m -1}} (\hat\varphi')}\cr}
\eqno £4.10.25.$$
   
 \smallskip
Recall that we have the exact sequence of {\it contravariant\/} functors from $\tilde\F$ to $\Ab$ (cf.~£2.7 and £2.8)
 $$0\too \tilde\frak l^{^{U,m - 1}}_\F/\tilde\frak l^{^{U,m}}_\F\too
 \widetilde{\Ker} ({\pi^{_{{\rm b},U,m}}})\too \widetilde{\Ker} ({\pi^{_{{\rm b},U,m - 1}}}) \too 0
 \eqno £4.10.26;$$
 hence, for another $\P\-$morphism $\hat\psi: T\to R$ we clearly have
$$\widehat{\lambda^{^{m -1}} (\hat\varphi)}\. \widehat{\lambda^{^{m -1}} (\hat\psi)}
 = \widehat{\lambda^{^{m -1}} (\hat\varphi\. \hat\psi)} \. \gamma_{\hat\psi,\hat\varphi}^{_m}
 \eqno £4.10.27\phantom{.}$$
 for some $\gamma_{\psi,\varphi}^{_m}$ in 
  $(\tilde\frak l^{^{U,m - 1}}_\F/\tilde\frak l^{^{U,m}}_\F) (T)\,.$
 That is to say, borrowing notation and terminology from  [8,~A2.8], we get a correspondence sending any {\it $\P\-$chain\/}  $\frak q\,\colon \Delta_2\to \P$  to the element 
 $ \gamma^{_m}_{\frak q (0\bullet 1),\frak q(1\bullet 2)}$ 
 in $(\tilde\frak l^{^{U,m - 1}}_\F/\tilde\frak l^{^{U,m}}_\F) \big(\frak q(0)\big)$ and, setting
$$\Bbb C^n \big(\tilde\P, \tilde\frak l^{^{U,m - 1}}_\F/\tilde\frak l^{^{U,m}}_\F\big) = 
\prod_{\tilde\frak q\in \Fct(\Delta_n,\tilde\P)} (\tilde\frak l^{^{U,m - 1}}_\F/\tilde\frak l^{^{U,m}}_\F) \big(\tilde\frak q(0)\big)
\eqno £4.10.28\phantom{.}$$
 for any $n\in \Bbb N\,,$  we claim that this correspondence determines an 
 {\it  stable\/} element $\gamma^{_m}$ of 
 $ \Bbb C^2 \big(\tilde\P,\tilde\frak l^{^{U,m - 1}}_\F/\tilde\frak l^{^{U,m}}_\F\big) $ 
 [8,~A3.17]; note that $\tilde\P \cong \tilde\F\,.$

\smallskip
 Indeed, for another $\P-$isomorphic {\it $\P\-$chain\/} $\frak q'\,\colon \Delta_2\to \P$ and a {\it natural $\P\-$isomorphism\/} $\nu\,\colon \frak q\cong\frak q'\,,$
setting  
$$\eqalign{T = \frak q (0)  ,\ T' \!= \frak q'(0) ,\ R = \frak q(1)
& ,\ R' \! = \frak q' (1) ,\ Q = \frak q (2) ,\ Q' \!= \frak q' (2)\cr
\hat\psi = \frak q (0\bullet 1)\!\!\quad,\quad\!\! \hat\varphi = \frak q(1\bullet 2)\!\!\quad &,\quad  \!\!\hat\psi' = \frak q'(0\bullet 1)\!\!\quad,\quad \!\!\hat\varphi' = \frak q'(1\bullet 2)\cr
\nu_0 = \hat\eta\quad,\quad  \nu_1 = &\  \hat\xi  \qq \nu_2 =  \hat\zeta\cr}
\eqno £4.10.29,$$
 from~£4.5.30 we have
 $$\eqalign{\widehat{\lambda^{^{m -1}} (\hat\varphi')} &=  \widehat{\lambda^{^{m -1}} (\hat\zeta)}\. \widehat{\lambda^{^{m -1}}(\hat\varphi)} \.  \widehat{\lambda^{^{m -1}} (\hat\xi)}^{-1}\cr
 \widehat{\lambda^{^{m -1}} (\hat\psi')} &=  \widehat{\lambda^{^{m -1}} (\hat\xi)}\. 
 \widehat{\lambda^{^{m -1}} (\hat\psi)} \.  \widehat{\lambda^{^{m -1}} (\hat\eta)}^{-1} \cr 
\widehat{\lambda^{^{m -1}} (\hat\varphi'\.. \hat\psi')} 
&=  \widehat{\lambda^{^{m -1}} (\hat\zeta)} \. \widehat{\lambda^{^{m -1}} (\hat\varphi\. \hat\psi)} \. 
 \widehat{\lambda^{^{m -1}} (\hat\eta)}^{-1}\cr}
 \eqno £4.10.30\phantom{.}$$
 and therefore we get
$$\eqalign{&\widehat{\lambda^{^{m -1}} (\hat\varphi' \. \hat\psi')} \.  \gamma^{_m}_{\varphi',\psi'} 
= \widehat{\lambda^{^{m -1}} (\hat\varphi')} \.  \widehat{\lambda^{^{m -1}} (\hat\psi')} \cr
& = \big(\widehat{\lambda^{^{m -1}} (\hat\zeta)} \. \widehat{\lambda^{^{m -1}} (\hat\varphi)} \. 
\widehat{\lambda^{^{m -1}} (\hat\xi)^{-1}}\big)\. 
\big(\widehat{\lambda^{^{m -1}} (\hat\xi)} \. \widehat{\lambda^{^{m -1}} (\hat\psi)} \. 
 \widehat{\lambda^{^{m -1}} (\hat\eta)}^{-1} \big)\cr
& =  \widehat{\lambda^{^{m -1}} (\hat\zeta)} \. \big(\widehat{\lambda^{^{m -1}} (\hat\varphi \. \hat\psi) }\.  \gamma^{_m}_{\hat\varphi,\hat\psi} \big) \.  \widehat{\lambda^{^{m -1}} (\hat\eta)}^{-1}\cr
&=  \widehat{\lambda^{^{m -1}} (\hat\varphi'\. \hat\psi')} \. 
\big((\tilde\frak l^{^{U,m - 1}}_\F/\tilde\frak l^{^{U,m}}_\F)
( \widehat{\lambda^{^{m -1}} (\hat\eta)}^{-1})\big)( \gamma^{_m}_{\hat\varphi,\hat\psi} ) \cr}
\eqno £4.10.31,$$ 
so that, by the divisibility of $ {\L^{^{{\rm b},U,m }}}\,,$ we have
$$\gamma^{_m}_{\hat\varphi',\hat\psi'}  = 
\big((\tilde\frak l^{^{U,m - 1}}_\F/\tilde\frak l^{^{U,m}}_\F)
( \widehat{\lambda^{^{m -1}} (\hat\eta)}^{-1})\big)( \gamma^{_m}_{\hat\varphi,\hat\psi})
\eqno £4.10.32;$$ 
\eject
\noindent
this proves  that the correspondence $\gamma^{_m}$ sending $(\hat\varphi,\hat\psi)$ 
to~$\gamma^{m}_{\hat\varphi,\hat\psi}$ is {\it stable\/} and, in particular, that 
$\gamma^{m}_{\hat\varphi,\hat\psi}$ only depends on the corresponding $\tilde\P\-$morphisms 
$\skew2\tilde{\hat\varphi}$ and~$\skew2\tilde{\hat\psi}\,;$ thus we set $\gamma^{_m}_{\tilde\varphi,\tilde\psi}= \gamma^{_m}_{\hat\varphi,\hat\psi}$ where $\varphi$ and $\psi$ are the correponding $\F\-$morphisms.

\smallskip
On the other hand, considering the usual differential map
$$d_2 : \Bbb C^2 \big(\tilde\P,\tilde\frak l^{^{U,m - 1}}_\F/\tilde\frak l^{^{U,m}}_\F\big)\too 
\Bbb C^3 \big(\tilde\P,\tilde\frak l^{^{U,m - 1}}_\F/\tilde\frak l^{^{U,m}}_\F\big)
\eqno £4.10.33,$$
we claim that $d_2 (\gamma^{_m}) = 0\,;$ indeed, for a third $\F\-$morphism $\varepsilon\,\colon W\to T$ we get
$$\eqalign{\big(\widehat{\lambda^{^{m -1}} (\hat\varphi)} &\.\widehat{\lambda^{^{m -1}} (\hat\psi)}  \big)\. \widehat{\lambda^{^{m -1}} (\hat\varepsilon)} 
= (\widehat{\lambda^{^{m -1}} (\hat\varphi\. \hat\psi)}\. \gamma^{_m}_{\tilde\varphi,\tilde\psi})\. \widehat{\lambda^{^{m -1}} (\hat\varepsilon)} \cr
&= \big(\widehat{\lambda^{^{m -1}} (\hat\varphi\. \hat\psi)} \. \widehat{\lambda^{^{m -1}} (\hat\varepsilon)}\big)\.\big((\tilde\frak l^{^{U,m - 1}}_\F/\tilde\frak l^{^{U,m}}_\F)(\tilde\varepsilon)\big)(\gamma^{_m}_{\tilde\varphi,\tilde\psi})\cr
&= \widehat{\lambda^{^{m -1}} (\hat\varphi\. \hat\psi\. \hat\varepsilon)} \.\gamma^{_m}_{\tilde\varphi\.\tilde\psi,\tilde\varepsilon}\.\big((\tilde\frak l^{^{U,m - 1}}_\F/\tilde\frak l^{^{U,m}}_\F)(\tilde\varepsilon)\big) (\gamma^{_m}_{\tilde\varphi,\tilde\psi})\cr
\widehat{\lambda^{^{m -1}} (\hat\varphi)} &\. \big(\widehat{\lambda^{^{m -1}} (\hat\psi)} \. \widehat{\lambda^{^{m -1}} (\hat\varepsilon)} \big)
 = \widehat{\lambda^{^{m -1}} (\hat\varphi)} \.\big(\widehat{\lambda^{^{m -1}} (\hat\psi\. \hat\varepsilon)} \.  \gamma^{_m}_{\tilde\psi,\tilde\varepsilon}\big)\cr 
& = \widehat{\lambda^{^{m -1}} (\hat\varphi\. \hat\psi\. \hat\varepsilon)} \. \gamma^{_m}_{\tilde\varphi,\tilde\psi\. \tilde\varepsilon}\. \gamma^{_m}_{\tilde\psi,\tilde\varepsilon}\cr}
\eqno £4.10.34\phantom{.}$$
and   the {\it divisibility\/} of ${\L^{^{{\rm b},U,m}} }$ forces
$$\gamma^{_m}_{\tilde\varphi\. \tilde\psi,\tilde\varepsilon}
\.\big((\tilde\frak l^{^{U,m - 1}}_\F/\tilde\frak l^{^{U,m}}_\F)(\tilde\varepsilon)\big)
(\gamma^{_m}_{\tilde\varphi,\tilde\psi}) = \gamma^{_m}_{\tilde\varphi,\tilde\psi\. \tilde\varepsilon} \.\gamma^{_m}_{\tilde\psi,\tilde\varepsilon}
\eqno £4.10.35;$$
since ${\rm Ker}(\widetilde{\pi_W^{^{{\rm b},U,m  }}})$ is abelian, with the additive notation we obtain
$$0 = \big((\tilde\frak l^{^{U,m - 1}}_\F/\tilde\frak l^{^{U,m}}_\F) (\tilde\varepsilon)\big)
(\gamma^{_m}_{\tilde\varphi,\tilde\psi}) - \gamma^{_m}_{\tilde\varphi,\tilde\psi\. \tilde\varepsilon} + \gamma^{_m}_{\tilde\varphi\. \tilde\psi,\tilde\varepsilon} - 
\gamma^{_m}_{\tilde\psi,\tilde\varepsilon}
\eqno £4.10.36,$$
proving our claim.

\smallskip
At this point, it follows from  equality £4.10.4  that  $\gamma^{_m} = d_1 (\beta^{^{m}})$ for some {\it stable\/} element $\beta^{^{m}} = 
(\beta^{^{m}}_{\tilde\frak r})_{\tilde\frak r\in \Fct(\Delta_1,\tilde\P)}$ in 
$\Bbb C^1 \big(\tilde\P,\tilde\frak l^{^{U,m - 1}}_\F/\tilde\frak l^{^{U,m}}_\F\big)\,;$ 
that is to say, with the notation above we get
$$\gamma^{_m}_{\tilde\varphi,\tilde\psi} = 
\big((\tilde\frak l^{^{U,m - 1}}_\F/\tilde\frak l^{^{U,m}}_\F)(\tilde\psi)\big)
(\beta^{^m}_{\tilde\varphi})\.(\beta^{^{m}}_{\tilde\varphi\. \tilde\psi})^{-1}\.\beta^{^m}_{\tilde\psi}
\eqno £4.10.37;$$
 hence, from equality~£4.10.27 we obtain
$$\eqalign{\big(\widehat{\lambda^{^{m -1}} (\hat\varphi)} &\. (\beta^{^ m}_{\tilde\varphi})^{-1}\big)
\.\big(\widehat{\lambda^{^{m -1}} (\hat\psi )}  \.( \beta^{^m}_{\tilde\psi})^{-1}\big)\cr
&=   \big( (\widehat{\lambda^{^{m -1}} (\hat\varphi)} \.\widehat{\lambda^{^{m -1}} (\hat\psi)} \big) \. \Big(\beta^{^m}_{\tilde\psi} \. \big((\tilde\frak l^{^{U,m - 1}}_\F/\tilde\frak l^{^{U,m}}_\F)(\tilde\psi)\big)(\beta^{^m}_{\tilde\varphi}) \Big)^{-1}\cr
&= \widehat{\lambda^{^{m -1}} (\hat\varphi\. \hat\psi)} \.(\beta^{^m}_{\tilde\varphi\circ\tilde\psi})^{-1}\cr}
\eqno £4.10.38,$$
which amounts to saying that the correspondence $\lambda^{^m}$  sending 
$\hat\varphi\in \P (Q,R)$ to~$\widehat{\lambda^{^{m -1}} (\hat\varphi)} \. (\beta^{^m}_{\tilde\varphi})^{-1}\in {\L^{^{{\rm b}, U,m}}}(Q,R)$ defines  the announced functor; note that, if 
$Q = R$ and $\hat\varphi = \tau_Q (u)$ for some $u\in Q\,,$ we have $\tilde\varphi = \widetilde{{\rm id}}_Q$ and $\beta^{^m}_{\tilde\varphi} = 1\,,$\break so that
 $\lambda^{^m} \big( \tau_Q (u)\big) = {\tau^{^{{\rm b},U,m }}_Q}(u) \,.$ It remains to prove that this  {\it functorial\/} section fulfills the commutativity of the corresponding diagram £4.10.7; since we already have the commutativity of diagram £4.10.20, it suffices to get the commutativity of the following diagram
$$\matrix{\P (R) & \buildrel \hat\mu_{_R}^{_m}\over \too &{\L^{^{{\rm b},U,m }}} (R)\cr
\hskip-20pt{\scriptstyle {\rm id}_{\P (R)}}\big\downarrow&
& \big\downarrow {\scriptstyle a_{(\beta^{^m}_{\tilde\varphi})^{-1}}}\hskip-30pt\cr
\P (R) & \buildrel \hat\mu_{_R}^{_m}\over \too &{\L^{^{{\rm b},U,m }}} (R)\cr}
\eqno £4.10.39$$
which follows from the fact that $\beta^{^m}$ is {\it stable\/} and therfore 
$(\beta^{^m}_{\tilde\varphi})^{-1}$ fixes the image of $\hat\mu_{_R}^{_m}\,.$

\smallskip
We can  modify this correspondence in order to get an {\it $\F\-$locality functor\/}; indeed, for any 
$\P\-$morphism $\tau_{Q,R} (u) \colon R\to Q$ where $u$ belongs to $\T_P (Q,R)\,,$ the 
${\L^{^{{\rm b}, U,m}}}(Q,R)\-$morphisms 
$\lambda^{^m} \big(\tau_{Q,R} (u)\big)$ and ${\tau^{^{{\rm b},U,m}}_{_{Q,R}}}(u)\,,$  both lift 
$\kappa_{Q,R} (u)$ in $\F (Q,R)\,;$ thus,  the {\it divisibility\/} of $ {\L^{^{{\rm b}, U,m}}}$ guarantees the existence and the uniqueness of $\delta_{\kappa_{Q,R} (u)}\in {\rm Ker}({\pi^{^{{\rm b},U,m}}_R})$ fulfilling
$${\tau^{^{{\rm b},U,m}}_{_{Q,R}}}(u) = \lambda^{^m} \big(\tau_{Q,R} (u)\big).\delta_{\kappa_{Q,R} (u)}
\eqno £4.10.40\phantom{.}$$
and, since we have $\lambda^{^m} \big( \tau_Q (w)\big) = {\tau^{^{{\rm b},U,m }}_Q}(w)$ for any 
$w\in Q\,,$  it is quite clear that $\delta_{\kappa_{Q,R} (u)}$ only depends on the class of
$\kappa_{Q,R} (u)$ in $\tilde\F (Q,R)$

\smallskip
 For a second $\P\-$morphism $\tau_{R,T} (v)\,\colon T\to R\,,$ setting 
 $\hat\xi = \tau_{R,T} (u)$ and $\hat\eta = \tau_{R,T} (v)$ we get
$$\eqalign{\lambda^{^m} (\hat\xi\. \hat \eta)\. \delta_{\tilde\xi\circ\tilde\eta} 
&= {\tau^{^{{\rm b},U,m}}_{_{Q,T}}}(uv) = {\tau^{^{{\rm b},U,m}}_{_{Q,R}}} (u) \. {\tau^{^{{\rm b},U,m}}_{_{R,T}}}(v)\cr
&= \lambda^{^m} (\hat\xi)\. \delta_{\tilde\xi}\.\lambda^{^m} (\hat\eta).  \delta_{\tilde\eta}\cr
&= \lambda^{^m} (\hat\xi\. \hat\eta)\.\big(\widetilde\Ker ({\pi^{^{{\rm b},U,m}}})(\tilde\eta)\big)(\delta_{\tilde\xi})\. \delta_{\tilde\eta}\cr}
\eqno £4.10.41;$$
then, once again  the {\it divisibility\/} of $ {\L^{^{{\rm b}, U,m}}}$ forces
$$\delta_{\tilde\xi\circ \tilde\eta} = \big(\widetilde\Ker ({\pi^{_{{\rm b},U,m}}})
(\tilde\eta)\big)(\delta_{\tilde\xi})\. \delta_{\tilde\eta}
\eqno £4.10.42\phantom{.}$$
and, since ${\rm Ker}({\pi^{^{{\rm b},U,m}}_T})$ is abelian, with the additive notation we obtain
$$0 =\big(\widetilde\Ker ({\pi^{^{{\rm b},U,m}}}) (\tilde\eta)\big)(\delta_{\tilde\xi}) -  \delta_{\tilde\xi\. \tilde\eta}+ \delta_{\tilde\eta}
\eqno £4.10.43.$$

\smallskip
That is to say,  denoting by $\frak i\,\colon \tilde\F_{\! P}\i \tilde\F$ the obvious {\it inclusion functor\/}, the correspondence $\delta$ sending any $\tilde\F_{\! P}\-$morphism 
$\tilde\xi\,\colon R\to Q$ to~$\delta_{\tilde\xi}$ defines a {\it $1\-$cocycle\/} in 
$\Bbb C^1\big(\tilde\F_{\! P},\widetilde\Ker ({\pi^{^{{\rm b},U,m}}}) \circ\frak i\big)\,;$ but, since the category $\tilde\F_{\! P}$ has a final object, we actually have [8,~Corollary~A4.8]
$$\Bbb H^1\big(\tilde\F_{\! P},\widetilde\Ker ({\pi^{^{{\rm b},U,m}}})
\circ\frak i\big) = \{0\}
\eqno £4.10.44;$$
consequently, we obtain $\delta = d_0 (w)$ for some element 
$w = (w_Q)_{Q\i P}$ in 
$$\Bbb C^0\big(\tilde\F_{\! P},\widetilde\Ker ({\pi^{_{{\rm b},U,m}}}) \circ\frak i\big) = \Bbb C^0\big(\tilde\F,\widetilde\Ker ({\pi^{^{{\rm b},U,m}}}))\big)
\eqno £4.10.45.$$
In conclusion, equality~£4.10.40 becomes
$$\eqalign{{\tau^{^{{\rm b},U,m}}_{_{Q,R}}}(u) &= 
\lambda^{^m} \big(\tau_{Q,R} (u)\big) \.\big(\widetilde\Ker ({\pi^{^{{\rm b},U,m}}})
(\widetilde{\tau_{Q,R} (u)}\big)(w_Q)\.w_R^{-1}\cr
& = w_Q\.\lambda^{^m} \big(\tau_{Q,R} (u)\big)\.w_R^{-1}\cr}
\eqno £4.10.46\phantom{.}$$
and therefore  the new correspondence  sending $\hat\varphi\in \P (Q,R)$ to  
$w_Q\.\lambda^{^m} \big(\hat\varphi)\.w_R^{-1}$ defines  a {\it $\F\-$locality functor\/}.   
 From now on, we still denote by $\lambda^{^m}$ this  {\it $\F\-$lo-cality functor}.

\smallskip
Let $\lambda'^{^m}\,\colon \P\to {\L^{^{{\rm b},U,m}}}$ be another $\F\-$locality functor; arguing by induction on $\vert\C_P - \N\vert$ and on $m\,,$ and up to natural $\F\-$isomorphisms, 
we clearly may assume that  $\lambda'^{^m}$ also lifts $\lambda^{^{m-1}}\,;$ in this case, for any $\P\-$morphism $\hat\varphi\,\colon R\to Q\,,$ we have 
$\lambda'^{^m} (\varphi) = \lambda^{^m} (\hat\varphi) \.  \varepsilon^{_m}_{\hat\varphi}$  for some $\varepsilon^{_m}_{\hat\varphi}$ in~$(\tilde\frak l^{^{U,m - 1}}_\F/\tilde\frak l^{^{U,m}}_\F) (R)\,;$ that is to say,  as above we  get a correspondence sending any {\it $\P\-$chain\/} 
 $\frak q\,\colon \Delta_1\to \P$   to~$ \varepsilon^{_m}_{\frak q (0\bullet 1),}$ 
 in $(\tilde\frak l^{^{U,m - 1}}_\F/\tilde\frak l^{^{U,m}}_\F) \big(\frak q(0)\big)$ and  we claim that this correspondence determines a  {\it  $\P\-$stable\/} element $\varepsilon^{_m}$ of 
 $ \Bbb C^1 \big(\tilde\P,\tilde\frak l^{^{U,m - 1}}_\F \! /\tilde\frak l^{^{U,m}}_\F\big) $ 
 [8,~A3.17].

\smallskip
 Indeed, for another $\P-$isomorphic {\it $\P\-$chain\/} $\frak q'\,\colon \Delta_1\to \P$ and a {\it natural $\P\-$isomorphism\/} $\nu\,\colon \frak q\cong\frak q'\,,$
as in~£4.10.29 above setting  
$$\eqalign{R = \frak q(0)\  ,\ R'  = \frak q' (0)\  ,\  &Q = \frak q (1)\ ,\ Q' = \frak q' (1)\cr
\hat\varphi = \frak q(0\bullet 1)\quad  ,\quad \!\! & \hat\varphi' = \frak q'(0\bullet 1)\cr 
\nu_0 = \hat\xi  \qq &
\nu_1 = \hat\zeta\cr}
\eqno £4.10.47,$$
from £4.10.24 we get
$$\eqalign{ \lambda'^{m} (\hat\varphi')  &=  \hat\zeta\. \lambda^{^m} (\varphi) \.  \varepsilon^{_m}_\varphi \. \hat\xi^{-1}\cr
& = \lambda^{^m} (\hat\varphi') \. \big((\tilde\frak l^{^{U,m - 1}}_\F/\tilde\frak l^{^{U,m}}_\F)(\widetilde{ \xi}^{-1})\big)(\varepsilon^{_m}_{\hat\varphi})\cr
 \lambda'^{^m} (\hat\varphi')  &=  \lambda^{^m} (\hat\varphi') \. \varepsilon^{_m}_{\hat\varphi'}\cr}
\eqno £4.10.48\phantom{.}$$
and the divisibility of  ${\L^{^{{\rm b},U,m}}}$  forces
$$\varepsilon^{_m}_{\hat\varphi'} =  \big((\tilde\frak l^{^{U,m - 1}}_\F/\tilde\frak l^{^{U,m}}_\F)(\widetilde{\xi}^{-1}) \big)(\varepsilon^{_m}_{\hat\varphi})
\eqno £4.10.49;$$
\eject
\noindent
this proves  that the correspondence $\varepsilon^{_m}$ sending $\hat\varphi$ 
to~$\varepsilon^{_m}_{\hat\varphi}$ is {\it $\P\-$stable\/} and, in particular, that 
$\varepsilon^{_m}_{\hat\varphi}$ only depends on the corresponding $\tilde\F\-$morphism 
$\tilde\varphi\,,$  thus we set $\varepsilon^{_m}_{\tilde\varphi}= \varepsilon^{_m}_{\hat\varphi}\,.$

\smallskip
On the other hand, considering the usual differential map
$$d_1 : \Bbb C^1 \big(\tilde\P,\tilde\frak l^{^{U,m - 1}}_\F/\tilde\frak l^{^{U,m}}_\F\big)\too 
\Bbb C^2 \big(\tilde\P,\tilde\frak l^{^{U,m - 1}}_\F/\tilde\frak l^{^{U,m}}_\F\big)
\eqno £4.10.50,$$
we claim that $d_1 (\varepsilon^{_m}) = 0\,;$ indeed, for a second $\P\-$morphism 
$\hat\psi\,\colon T\to R$ we get
$$\eqalign{\lambda'^{^m} (\hat\varphi) \. \lambda'^{^m} (\hat\psi) 
&= \lambda^{^m} (\hat\varphi)\. \varepsilon^{_m}_{\tilde\varphi}\. \lambda^{^m} (\tilde\psi)\. \varepsilon^{_m}_{\hat\psi}\cr
&= \lambda^{^m} (\hat\varphi\. \hat\psi)\. \big((\tilde\frak l^{^{U,m - 1}}_\F/\tilde\frak l^{^{U,m}}_\F)(\tilde\psi)\big)(\varepsilon^{_m}_{\tilde\varphi}) \.  \varepsilon^{_m}_{\tilde\psi} \cr
\lambda'^{^m} (\hat\varphi) \. \lambda'^{^m} (\hat\psi) 
 &=  \lambda^{^m} (\hat\varphi\. \hat\psi)\. \varepsilon^{_m}_{\tilde\varphi\. \tilde\psi}\cr }
\eqno £4.10.51\phantom{.}$$
and   the {\it divisibility\/} of ${\L^{^{{\rm b},U,m}} }$ forces
$$\big((\tilde\frak l^{^{U,m - 1}}_\F/\tilde\frak l^{^{U,m}}_\F)(\tilde\psi)\big)(\varepsilon^{_m}_{\tilde\varphi}) \.  \varepsilon^{_m}_{\tilde\psi } =  \varepsilon^{_m}_{\tilde\varphi\. \tilde\psi}
\eqno £4.10.52;$$
since ${\rm Ker}({\pi_T^{^{{\rm b},U,m  }}})$ is Abelian, with the additive notation we obtain
$$0 = \big((\tilde\frak l^{^{U,m - 1}}_\F/\tilde\frak l^{^{U,m}}_\F)(\tilde\psi)\big)(\varepsilon^{_m}_{\tilde\varphi}) -\varepsilon^{_m}_{\tilde\varphi\. \tilde\psi}
+  \varepsilon^{_m}_{\tilde\psi }
\eqno £4.10.53,$$
proving our claim.

\smallskip
At this point, it follows from equality £4.10.4   that  $\varepsilon^{_m} = d_0 (\hat\nu)$ for some 
{\it stable\/} element $\hat\nu = (\hat\nu_Q)_{Q\i P}$ in 
$\Bbb C^0 \big(\tilde\P,\tilde\frak l^{^{U,m - 1}}_\F/\tilde\frak l^{^{U,m}}_\F\big)\,;$ 
that is to say, with the notation above we get
$$\varepsilon^{_m}_{\tilde\varphi} = 
\big((\tilde\frak l^{^{U,m - 1}}_\F/\tilde\frak l^{^{U,m}}_\F)(\tilde\varphi)\big)(\hat\nu_Q)\.\hat\nu_R^{-1}
\eqno £4.10.54;$$
 hence,  we obtain
$$\eqalign{\lambda'^{^m} (\hat\varphi) = 
\lambda^{^m} (\hat\varphi) \. \big((\tilde\frak l^{^{U,m - 1}}_\F/\tilde\frak l^{^{U,m}}_\F)(\tilde\varphi)\big)(\hat\nu_Q)\.\hat\nu_R^{-1} =  \hat\nu_Q\. \lambda^{^m} (\hat\varphi)\. \hat\nu_R^{-1}\cr}
\eqno £4.10.55,$$
which amounts to saying that  $\lambda'^{^m}$ is naturally $\F\-$isomorphic to 
$\lambda^{^m}\,.$ We are done

\bigskip
\noindent
{\bf Corollary £4.11.}  {\it There exists a unique perfect $\F\-$locality $\P$ up to natural $\F\-$isomorphisms.\/}
\medskip
\noindent
{\bf Proof:} The existence has been proved in Corollary £4.6 above and the uniqueness is an easy consequence of Theorem £4.10.

\vfill\eject
\bigskip
\bigskip
\centerline{\large References}
\bigskip
\noindent
[1]\phantom{.} Dave Benson, personal letter 1994
\smallskip\noindent
[2]\phantom{.} Carles Broto, Ran Levi and Bob Oliver,  {\it The homotopy theory
of fusion systems\/}, Journal of Amer. Math. Soc. 16(2003), 779-856.
\smallskip\noindent
[3]\phantom{.} Andrew Chermak. {\it Fusion systems and localities\/}, 
Acta Mathematica, 211(2013), 47-139.
\smallskip\noindent
[4]\phantom{.} Stefan Jackowski and James McClure, {\it Homotopy
decomposition of classifying spaces via elementary abelian subgroups\/},
Topology, 31(1992), 113-132.
\smallskip\noindent
[5]\phantom{.} George Glauberman \& Justin Lynd, {\it Control of fixed points and existence and
uniqueness of centric systems\/}, arxiv.org/abs/1506.01307.
\smallskip\noindent
[6] Bob Oliver. {\it Existence and Uniqueness of Linking Systems: Chermak's proof via obstruction theory\/}, 
Acta Mathematica, 211(2013), 141-175.
\smallskip\noindent
[7]\phantom{.}  Llu\'\i s Puig, {\it Brauer-Frobenius categories\/}, Manuscript notes 1993
\smallskip\noindent
[8] Llu\'\i s Puig, {\it ``Frobenius categories versus Brauer blocks''\/}, Progress in Math. 
274(2009), Birkh\"auser, Basel.
\smallskip\noindent
[9] Llu\'\i s Puig, {\it Existence, uniqueness and functoriality
of the perfect locality over a Frobenius $P\-$category\/}, arxiv.org/abs/1207.0066,
 Algebra Colloquium, 23(2016) 541-622.
\smallskip\noindent
[10] Llu\'\i s Puig, {\it A correction to the uniqueness of a partial perfect locality over a Frobenius $P$-category\/}, arxiv.org/abs/1706.04349, Algebra Colloquium. 26(2019) 541-559.
\smallskip\noindent
[11] Llu\'\i s Puig, {\it Categorizations of limits of  Grothendieck groups over a Frobenius P-category\/}, submitted to Algebra Colloquium.

\end